\documentclass{amsproc}

\usepackage{a4wide}
\usepackage{bbm}
\usepackage{enumerate}

\newtheorem{theorem}{Theorem}[section]
\newtheorem{lemma}[theorem]{Lemma}
\newtheorem{Prop}[theorem]{Proposition}
\newtheorem{Cor}[theorem]{Corollary}
\theoremstyle{definition}
\newtheorem{definition}[theorem]{Definition}
\newtheorem{example}[theorem]{Example}

\theoremstyle{remark}
\newtheorem{remark}[theorem]{Remark}

\numberwithin{equation}{section}

\newcommand{\bt}{\begin{theorem}\ \ }  
\newcommand{\et}{\end{theorem}}  
\newcommand{\bp}{\begin{Prop}\ \ }  
\newcommand{\ep}{\end{Prop}}  
\newcommand{\bc}{\begin{Cor}\ \ }  
\newcommand{\ec}{\end{Cor}}  
\newcommand{\bl}{\begin{lemma}\ \ }  
\newcommand{\el}{\end{lemma}}  
\newcommand{\bd}{\begin{definition}\ \ }  
\newcommand{\ed}{\end{definition}}  
\newcommand{\pf}{\begin{proof}}  
\newcommand{\epf}{\end{proof}}  
\newcommand{\br}{\begin{remark}\ \ }
\newcommand{\er}{\end{remark}}
\newcommand{\brsn}{\begin{remarks*}\ \ }
\newcommand{\ersn}{\end{remarks*}}

\renewcommand{\a}{\alpha}
\renewcommand{\b}{\beta}
\renewcommand{\o}{\omega}
\renewcommand{\l}{\lambda}  
\renewcommand{\L}{\Lambda}
\newcommand{\ga}{\gamma}
\renewcommand{\t}{\tau}  
\newcommand{\g}{\mathfrak{g}}
\newcommand{\h}{\mathfrak{h}}
\newcommand{\z}{\mathfrak{z}}
\renewcommand{\k}{\mathfrak{k}}
\newcommand{\gl}{\mathfrak{l}} 
\newcommand{\Sp}{{\rm Sp}}
\newcommand{\GL}{{\rm GL}}
\newcommand{\SL}{{\rm SL}}
\newcommand{\SO}{{\rm SO}}
\newcommand{\Og}{{\rm O}}
\newcommand{\U}{{\rm U}}
\newcommand{\SU}{{\rm SU}}
\newcommand{\Spin}{{\rm Spin}}
\newcommand{\G}{{\rm G}}
\newcommand{\Stab}{{\rm Stab}}

\newcommand\re[1]{(\ref{#1})}
\newcommand{\arr}{\begin{array}{rlll}}
\newcommand{\ea}{\end{array}}
\newcommand{\bea}{\begin{eqnarray}}
\newcommand{\eea}{\end{eqnarray}}
\newcommand{\bean}{\begin{eqnarray*}}
\newcommand{\eean}{\end{eqnarray*}}

\newcommand{\op}{\oplus}
\newcommand{\ot}{\otimes}
\newcommand{\ccdot}{\!\cdot\!}

\newcommand{\n}{\nabla}

\newcommand{\bR}{\mathbb{R}}
\newcommand{\bC}{\mathbb{C}}

\newcommand{\bO}{\mathbb{O}}

\usepackage{
amsfonts,
amssymb,
amsmath,
amsthm,
amsgen
}
\newcommand{\cal}{\mathcal}

\newcommand{\R}{\ensuremath{\mathbb{R}}}

\newcommand{\bcase}{\begin{case}}
\newcommand{\ecase}{\end{case}}

\newcommand{\bclaim}{\begin{claim}}
\newcommand{\eclaim}{\end{claim}}

\newcommand{\bstep}{\begin{step}}
\newcommand{\estep}{\end{step}}

\newcommand{\bhlem}{\begin{hlem}}
\newcommand{\ehlem}{\end{hlem}}

\newcommand{\bleer}{\begin{leer}}
\newcommand{\eleer}{\end{leer}}
\newcommand{\bde}{\begin{definition}}
\newcommand{\ede}{\end{definition}}

\newcommand{\ol}{\overline}

\newcommand{\bs}{\begin{Prop}}
\newcommand{\es}{\end{Prop}}
\newcommand{\btheo}{\begin{theorem}}
\newcommand{\etheo}{\end{theorem}}
\newcommand{\bfolg}{\begin{Cor}}
\newcommand{\efolg}{\end{Cor}}
\newcommand{\blem}{\begin{lemma}}
\newcommand{\elem}{\end{lemma}}
\newcommand{\bnote}{\begin{note}}
\newcommand{\enote}{\end{note}}
\newcommand{\bprf}{\begin{proof}}
\newcommand{\eprf}{\end{proof}}
\newcommand{\be}{\begin{eqnarray*}}
\newcommand{\ee}{\end{eqnarray*}}
\newcommand{\beqa}{\begin{eqnarray}}
\newcommand{\eeqa}{\end{eqnarray}}
\newcommand{\bi}{\begin{itemize}}
\newcommand{\ei}{\end{itemize}}
\newcommand{\bnum}{\begin{enumerate}}
\newcommand{\enum}{\end{enumerate}}
\newcommand{\la}{\langle}
\newcommand{\ra}{\rangle}

\newcommand{\ve}{\varepsilon}

\newcommand{\beq}{\begin{equation}}
\newcommand{\eeq}{\end{equation}}

\newcommand{\rr}{\mathbb{R}}

\newcommand{\ccc}{\mathbb{C}}

\newcommand{\vf}{\varphi}
\newcommand{\barr}{\[\begin{array}}
\newcommand{\earr}{\end{array}\]}
\newcommand{\bvec}{\left(\begin{array}{c}}
\newcommand{\evec}{\end{array}\right)}

\newcommand{\w}{\omega}

\newcommand{\s}{\sigma}

\newcommand{\ddt}{\frac{\partial}{\partial t}}

\newcommand{\bbem}{\begin{remark}}
\newcommand{\ebem}{\end{remark}}
\newcommand{\bbez}{\begin{bez}}
\newcommand{\ebez}{\end{bez}}
\newcommand{\bbsp}{\begin{bsp}}
\newcommand{\ebsp}{\end{bsp}}

\newcommand{\W}{\Omega}

\newcommand{\im}{\mathrm{im\ }}  

\renewcommand{\^}{\wedge}
\newcommand{\e}{\mathrm{e}}

\newcommand{\inter}{\hspace{-1mm}\makebox[10pt]{\rule{6pt}{.3pt}\rule{.3pt}{5pt}}}

\renewcommand{\i}{\mathrm{i}}

\begin{document}
 \title{Half-flat Structures and Special Holonomy}

 \author{V. Cort\'es} \address{Vicente Cort\'es, Department
   Mathematik, Universit\"at Hamburg, Bundesstra{\ss}e 55, D-20146
   Hamburg, Germany}
 \email{cortes@math.uni-hamburg.de}

 \author{T. Leistner} \address{Thomas Leistner, Department Mathematik,
   Universit\"at Hamburg, Bundesstra{\ss}e 55, D-20146 Hamburg,
   Germany}
 \email{leistner@math.uni-hamburg.de}

 \author{L. Sch\"afer} \address{Lars Sch\"afer, Institut
   Differentialgeometrie, Leibniz Universit\"at Hannover, Welfengarten
   1, D-30167 Hannover, Germany}
 \email{schaefer@math.uni-hannover.de}

 \author{F. Schulte-Hengesbach} \address{Fabian Schulte-Hengesbach,
   Department Mathematik, Universit\"at Hamburg, Bundesstra{\ss}e 55,
   D-20146 Hamburg, Germany}
 \email{schulte-hengesbach@math.uni-hamburg.de}

 \thanks{This work was supported by the SFB 676 of the Deutsche Forschungsgemeinschaft.}

\begin{abstract}
  It was proven by Hitchin that any solution of his evolution
  equations for a half-flat $\SU (3)$-structure on a compact
  six-manifold $M$ defines an extension of $M$ to a seven-manifold
  with holonomy in $\G_2$.  We give a new proof, which does not
  require the compactness of $M$. More generally, we prove that the
  evolution of any half-flat $G$-structure on a six-manifold $M$
  defines an extension of $M$ to a Ricci-flat seven-manifold $N$, for
  any real form $G$ of $\SL (3,\bC )$. If $G$ is noncompact, then the
  holonomy group of $N$ is a subgroup of the noncompact form $\G_2^*$
  of $\G_2^\bC$.  Similar results are obtained for the extension of
  nearly half-flat structures by nearly parallel $\G_2$- or
  $\G_2^{*}$-structures, as well as for the extension of cocalibrated
  $\G_2$- and $\G_2^*$-structures by parallel $\mathrm{Spin}(7)$- and
  $\mathrm{Spin}_0(3,4)$-structures, respectively.  As an application,
  we obtain that any six-dimensional homogeneous manifold with an
  invariant half-flat structure admits a canonical extension to a
  seven-manifold with a parallel $\G_2$- or $\G_2^*$-structure.  For
  the group $H_3\times H_3$, where $H_3$ is the three-dimensional
  Heisenberg group, we describe all left-invariant half-flat
  structures and develop a method to explicitly determine the
  resulting parallel $\G_2$- or $\G_2^*$-structure without
  integrating. In particular, we construct three eight-parameter
  families of metrics with holonomy equal to $\G_2$ and
  $\G_2^*$. Moreover, we obtain a strong rigidity result for the
  metrics induced by a half-flat structure $(\o,\rho)$ on $H_3 \times
  H_3$ satisfying $\o(\z,\z)=0$ where $\z$ denotes the centre.
  Finally, we describe the special geometry of the space of stable
  three-forms satisfying a reality condition. Considering all possible
  reality conditions, we find four different special K\"ahler
  manifolds and one special para-K\"ahler manifold.
\end{abstract}

\maketitle
\setcounter{tocdepth}{2}
\tableofcontents

\section*{Introduction}
Following Hitchin \cite{H1}, a $k$-form $\varphi$ on a differentiable
manifold $M$ is called \emph{stable} if the orbit of $\varphi (p)$
under $\GL (T_pM)$ is open in $\Lambda^k T^*_pM$ for all $p\in M$.  In
this paper we are mainly concerned with six-dimensional manifolds $M$
endowed with a stable two-form $\omega$ and a stable three-form
$\rho$.  A stable three-form defines an endomorphism field $J_\rho$ on
$M$ such that $J_\rho^2 = \ve \mathrm{id}$, see \re{Jrho}. We will assume the
following algebraic compatibility equations between $\o$ and $\rho$:
\[ \omega \wedge \rho = 0,\quad J_\rho^*\rho \wedge \rho =\frac{2}{3}
\o^3 .\] The pair $(\omega ,\rho)$ defines an $\SU (p,q)$-structure if
$\ve =-1$ and an $\SL (3,\bR )$-structure if $\ve =+1$. In the former
case, the pseudo-Riemannian metric $\o (J_\rho \cdot ,\cdot )$ has
signature $(2p,2q)$.  In the latter case it has signature $(3,3)$.
The structure is called \emph{half-flat} if the pair $(\omega ,\rho)$
satisfies the following exterior differential system:
\[ d\o^2 = 0,\quad d\rho =0.\] In \cite{H1}, Hitchin introduced the
following evolution equations for a time-dependent pair of stable
forms $(\o (t) ,\rho (t))$ evolving from a half-flat $\SU
(3)$-structure $(\o (0) ,\rho (0))$:
\[ \frac{\partial}{\partial t}\rho = d\w,\quad
\frac{\partial}{\partial t} \hat{\o} = d\hat{\rho}, \] where $\hat{\o}
= \frac{\o^2}{2}$ and $\hat{\rho}=J_\rho^*\rho$.  For compact
manifolds $M$, he showed that these equations are the flow equations
of a certain Hamiltonian system and that any solution defined on some
interval $0\in I\subset \bR$ defines a Riemannian metric on $M\times
I$ with holonomy group in $\G_2$.  We give a new proof of this
theorem, which does not use the Hamiltonian system and does not assume
that $M$ is compact. Moreover, our proof yields a similar result for
all three types of half-flat $G$-structures: $G=\SU (3), \SU(1,2)$ and
$\SL (3,\bR )$. For the noncompact groups $G$ we obtain a
pseudo-Riemannian metric of signature $(3,4)$ and holonomy group in
$\G_2^*$ on $M\times I$ (see Theorem \ref{paralleltheo}). As an
application, we prove that any six-manifold endowed with a real
analytic half-flat $G$-structure can be extended to a Ricci-flat
seven-manifold with holonomy group in $\G_2$ or $\G_2^*$, depending on
whether $G$ is compact or noncompact, see Corollary \ref{Cauchy}.

More generally, a $G$-structure $(\o ,\rho )$ is called \emph{nearly
  half-flat} if
\[ d\rho=\hat{\o}\] and a $\G_2$- or $\G_2^*$-structure defined by a
three-form $\varphi$ is called \emph{nearly parallel} if
\[ d\varphi = *_\varphi \varphi.\] We prove in Theorem
\ref{nearlyparalleltheo} that any solution $I\ni t \mapsto
\left(\w(t)= 2\widehat{d\rho} (t), \rho (t)\right)$ of the evolution
equation
\[ \dot{\rho} = d\o -\ve\hat{\rho}\] evolving from a nearly half-flat
$G$-structure $(\o (0),\rho (0) )$ on $M$ defines a nearly parallel
$\G_2$- or $\G_2^*$-structure on $M\times I$, depending on whether $G$
is compact or noncompact, see (\ref{hatsigma}) for the definition of
$\widehat{d\rho}$. For compact manifolds $M$ and $G=\SU (3)$ this
theorem was proven by Stock \cite{St}.

The above constructions are illustrated in Section \ref{NKexamples},
where we start with a nearly pseudo-K\"ahler or a nearly para-K\"ahler
six-manifold as initial structure.  These structures are both
half-flat and nearly half-flat and the resulting parallel or nearly
parallel $\G_2$- and $\G_2^*$-structures induce cone or (hyperbolic)
sine cone metrics.

In Section \ref{nilexamples}, we discuss the evolution of invariant
half-flat structures on nilmanifolds.  Lemma \ref{trick17} shows how
to simplify effectively the ansatz for a solution for a number of
nilpotent Lie algebras including the direct sum $\g=\h_3 \op \h_3$ of
two Heisenberg algebras.  Focusing on this case, we determine the
orbits of the ${\rm Aut}(\h_3 \op \h_3)$-action on non-degenerate
two-forms $\o$ on $\h_3$ which satisfy $d\o^2=0$. Based on this, we
describe all left-invariant half-flat structures $(\o,\rho)$ on $H_3
\times H_3$. A surprising phenomenon occurs in indefinite
signature. Under the assumption $\o(\z,\z)=0$, which corresponds to
the vanishing of the projection of $\o$ on a one-dimensional space,
the geometry of the metric induced by a half-flat structure
$(\o,\rho)$ is completely determined (Proposition \ref{k1=0}) and the
evolution turns out to be affine linear (Proposition \ref{afflin}).
However, this evolution produces only metrics that are decomposable
and have one-dimensional holonomy group.  On the other hand, we give
an explicit formula in Proposition \ref{k1ne0} for the parallel
three-form $\vf$ resulting from the evolution for any half-flat
structure $(\o,\rho)$ with $\o(\z,\z) \ne 0$. In fact, the formula is
completely algebraic such that the integration of the differential
equation is circumvented. In particular, we give a number of explicit
examples of half-flat structures of the second kind on $\h_3 \op \h_3$
which evolve to new metrics with holonomy group equal to $\G_2$ and
$\G_2^*$. Moreover, we construct an eight-parameter family of
half-flat deformations of the half-flat examples which lift to an
eight-parameter family of deformations of the corresponding parallel
stable three-forms in dimension seven.  Needless to say, those
examples of $\G_2^{(*)}$-metrics on $ M\times (a,b)$ for which
$(a,b)\not=\rr $ are geodesically incomplete. However, for $M$ compact
with an $\SU (3)$-structure, a conformal transformation produces
complete Riemannian metrics on $ M\times \rr$ that are conformally
parallel $\G_2$.

A $\G_2$- or $\G_2^*$-structure defined by a three-form $\varphi$ is
called \emph{cocalibrated} if
\[ d *_\varphi \varphi = 0.\] Hitchin proposed the following equation
for the evolution of a cocalibrated $\G_2$-structure $\varphi (0)$:
\[ \frac{\partial}{\partial t}( *_\varphi \varphi) = d \varphi.\] He
proved that any solution $I\ni t\mapsto \varphi (t)$ on a compact
manifold $M$ defines a Riemannian metric on $M\times I$ with holonomy
group in $\Spin (7)$. We generalise also this theorem to noncompact
manifolds and show that any solution of the evolution equation
starting from a cocalibrated $\G_2^*$-structure defines a
pseudo-Riemannian metric of signature $(4,4)$ and holonomy group in
$\Spin_0(3,4)$, see Theorem \ref{cocaThm}.

Homogeneous projective special pseudo-K\"ahler manifolds of semisimple
groups with compact stabiliser were classified in \cite{AC}.  It
follows that there is a unique homogeneous projective special
pseudo-K\"ahler manifold with compact stabiliser which admits a
transitive action of a real form of $\SL (3,\bC )$ by automorphisms of
the special K\"ahler structure, namely
\[ \frac{\SU (3,3)}{{\rm S}(\U (3)\times \U (3))}.\] Its special
K\"ahler metric is (negative) definite.  The above manifold occurred
in \cite{AC} as an open orbit of $\SU (3,3)$ on the projectivised
highest weight vector orbit of $\SL (6,\bC)$ on $ \Lambda^3(\bC^6)^*$.
The space of stable three-forms $\rho\in \Lambda^3(\bR^6)^*$, such
that $J^2_\rho = -1$, has also the structure of a special
pseudo-K\"ahler manifold \cite{H1}. The underlying projective special
pseudo-K\"ahler manifold is the manifold
\[ \frac{\SL(6,\bR)}{\U(1)\cdot \SL(3,\bC)}\] which has noncompact
stabiliser and indefinite special K\"ahler metric.  Both manifolds can
be obtained from the space of stable three-forms $\rho\in
\Lambda^3(\bC^6)^*$ by imposing two different reality conditions.  In
the last section of this paper we determine all homogeneous spaces
which can be obtained in this way and describe their special geometric
structures.  In particular, we calculate the signature of the special
K\"ahler metrics. For the projective special pseudo-K\"ahler manifold
$\SL(6,\bR)/\,(\U(1)\cdot \SL(3,\bC))$, for instance, we obtain the
signature $(6,12)$.  Apart from the two above examples, we find two
additional special pseudo-K\"ahler manifolds and also a special
para-K\"ahler manifold.  The latter is associated to the space of
stable three-forms $\rho\in \Lambda^3(\bR^6)^*$, such that $J^2_\rho =
+1$.

\section{Algebraic preliminaries}

\subsection{Stable forms}
In this section we will collect some basic facts about stable forms,
their orbits and their stabilisers.  
\bp \label{wherestable} Let $V$
be an $n$-dimensional real or complex vector space.  The general
linear group $\GL(V)$ has an open orbit in $\Lambda^kV^*$, $0\le k\le
\left[ \frac{n}{2}\right]$, if and only if $k\le 2$ or if $k=3$ and
$n=6,7$ or $8$.  \ep 
\pf The representation of $\GL(V)$ on $\Lambda^kV^*$ is
irreducible. In the complex case the result thus follows, for
instance, from the classification of irreducible complex
prehomogeneous vector spaces, \cite{KiS}. The result in the real case
follows from the complex case, since the complexification of the
$\GL(n,\bR )$-module $\Lambda^k\mathbb R^{n*}$ is an irreducible
$\GL(n,\bC )$-module.  \epf \br An open orbit is unique in the complex
case, since an orbit which is open in the usual topology is also
Zariski-open and Zariski-dense (Prop.\ 2.2, \cite{Ki}). Over the
reals, the number of open orbits is finite by a well-known theorem of
Whitney.  \er \bd A $k$-form $\rho \in \Lambda^k V^*$ is called
\emph{stable} if its orbit under $\GL(V)$ is open.  \ed
\bp \label{invmap} Let $k \in \{2,n-2\}$ and $n$ even, or $k \in \{ 3,
n-3\}$ and $n=6,7$ or $8$. There is a $\GL(V)$-equivariant mapping
\begin{eqnarray*} 
  \phi : \Lambda^k V^* \rightarrow \Lambda^n V^*, 
\end{eqnarray*} 
homogeneous of degree $\frac{n}{k}$, which
assigns a volume form to a stable $k$-form and which vanishes on
non-stable forms. Given a stable $k$-form $\rho$, the derivative of
$\phi$ in $\rho$ defines a dual $(n-k)$-form $\hat {\rho} \in
\Lambda^{n-k} V^*$ by the property
\begin{equation}
  \label{definitionhat}
  d \phi_\rho (\alpha) = \hat \rho \wedge \alpha \quad \mbox{for all $\alpha \in \Lambda^k V^*$.}
\end{equation}
The dual form $\hat \rho$ is also stable and satisfies
\[ (\Stab _{\GL(V)}(\rho))_0 = (\Stab _{\GL(V)}(\hat \rho))_0. \] A
stable form, its volume form and its dual are related by the formula
\begin{equation}
  \label{phieuler} 
  \hat \rho \wedge \rho = \frac{n}{k} \phi(\rho). 
\end{equation}
\ep \pf We consider the complex case first. As a result of the theory
of prehomogeneous vector spaces \cite{Ki}, the complement of the open
orbit is, in our situation, a hypersurface defined by a non-degenerate
homogeneous polynomial $f$ which is invariant under $\GL(V)$ up to a
non-trivial character. In other words, there is an equivariant mapping
from $\Lambda^k V^*$ to $(\Lambda^n V^*)^{\otimes s}$ for some
positive integer $s$. Taking the $s$-th root, which depends on the
choice of an orientation if $s$ is even, we obtain the equivariant map
$\phi$ with the claimed properties. The equivariance under scalar
matrices implies that the map $\phi$ is homogeneous of degree
$\frac{n}{k}$.

The derivative \[ \Lambda^k V^* \rightarrow (\Lambda^k V^*)^* \otimes
\Lambda^n V^* \stackrel{=}{\rightarrow} \Lambda^{n-k} V^* \, , \: \rho
\mapsto d_\rho \phi \mapsto \hat \rho \] inherits equivariance from
$\phi$ and is an immersion since $f$ is non-degenerate. Therefore, it
maps stable forms to stable forms such that the connected components
of the stabilisers are identical. Formula (\ref{phieuler}) is in fact
Euler's formula for the homogeneous mapping $\phi$.

Since the complexification of the $\GL(n,\bR )$-module
$\Lambda^k\mathbb R^{n*}$ is an irreducible $\GL(n,\bC )$-module, the
results in the real case are easily deduced from the complex case.
\epf

In the following, we discuss stable forms, their volume forms and
their duals in the cases which are relevant in this article. In each
case, $V$ is a real $n$-dimensional vector space.

$\mathbf{k=2, n=2m.}$ The orbit of a non-degenerate two-form is open and
there is only one open orbit in $\Lambda^2 V^*$. Thus, the stabiliser
of a stable two-form $\omega$ is isomorphic to
$\Sp(2m,\mathbb{R})$. The polynomial invariant is the Pfaffian
determinant. We normalise the associated equivariant volume form such
that is corresponds to the Liouville volume form
\begin{eqnarray*} 
  \phi(\omega)=\frac{1}{m!} \omega^m. 
\end{eqnarray*}
Differentiation of the homogeneous polynomial map $\o \mapsto
\phi(\omega)$ yields
\begin{eqnarray*}
  \hat \omega = \frac{1}{(m-1)!} \omega^{m-1}. 
\end{eqnarray*}

$\mathbf{k=(n-2), n=2m.}$ As $\L^{n-2}V^* \cong \L^2 V \otimes \L^n
V^*$, there is again only one open orbit. More precisely, an
$(n-2)$-form $\sigma$ is stable if and only if there is a stable
two-form $\omega$ with $\sigma = \hat \o$ since the mapping $\o
\mapsto \hat \o$ is an equivariant immersion. If $m$ is even, such an
$\omega$ is unique and we define the volume form $\phi(\sigma) =
\phi(\omega)$. If $m$ is odd, we need an orientation on $V$ to
uniquely define an associated volume form.  We choose the (m-1)-th
root $\omega$ with positively oriented $\omega^m$ and define again
$\phi(\sigma) = \phi(\omega)$. In both cases, we find
\begin{eqnarray}
  \label{hatsigma}
  \hat \sigma = \frac{1}{m-1} \omega 
\end{eqnarray}
with the help of (\ref{phieuler}). The stabiliser of a stable
four-form in $\GL^+(V)$ is again the real symplectic group.

$\mathbf{k=3, n=6.}$ Let $V$ be an oriented six-dimensional vector space
and let $\kappa$ denote the canonical isomorphism
\begin{eqnarray*} \kappa \, : \, \Lambda^k V^* \cong \Lambda^{6-k} V
  \otimes \Lambda^6 V^*.
\end{eqnarray*} 
Given any three-form $\rho$, we define $K \, : \, V \rightarrow V
\otimes \Lambda^6 V^*$ by
\begin{eqnarray*}
  K_\rho(v)=\kappa((v \lrcorner\, \rho) \wedge \rho)
\end{eqnarray*}
and the quartic invariant
\begin{eqnarray}
  \label{lambdarho}
  \lambda (\rho) = \frac{1}{6}\rm{tr}\, ( K_\rho^2 ) \: \in (\Lambda^6 V^*)^{\otimes 2}. 
\end{eqnarray}
Recall that, for any one-dimensional vector space $L$, an element $u
\in L^{\otimes 2}$ is defined to be positive, $u > 0$, if $u = s
\otimes s$ for some $s \in L$ and negative if $-u > 0$. Therefore, the
norm of an element $u \in L^{\otimes 2}$ is well-defined and we set
\begin{eqnarray}
  \label{phirho}
  \phi(\rho)= \sqrt{|\lambda(\rho)|} 
\end{eqnarray}
for the positively oriented square root. If $\phi(\rho) \ne 0$, we
furthermore define
\begin{eqnarray}
  \label{Jrho}
  J_\rho = \frac{1}{\phi(\rho)} K_\rho.  
\end{eqnarray}

\begin{Prop}
\label{orbitsof3forms}
A three-form $\rho$ on an oriented six-dimensional vector space $V$
with volume form $\nu$ is stable if and only $\lambda(\rho) \ne
0$. There are two open orbits.

One orbit consists of all three-forms $\rho$ satisfying one of the
following equivalent properties.
\begin{enumerate}[(a)]
\item $\lambda(\rho) > 0$
\item There are two uniquely defined real decomposable three-forms
  $\alpha$ and $\beta$ such that $\rho = \alpha + \beta$ and $\alpha
  \wedge \beta >0$.
\item The stabiliser of $\rho$ in $\GL^+(V)$ is $\SL(3,\bR) \times
  \SL(3,\bR)$.
\item It holds $\lambda(\rho) \ne 0$ and the endomorphism $J_\rho$ is
  a para-complex structure on $V$, i.e. $J^2 = \rm{id_V}$ and the
  eigenspaces for the eigenvalues $\pm 1$ are three-dimensional.
\item \label{enum1} There is a basis $\{ e_1, ..., e_6 \}$ of $V$ such
  that $\nu=e^{123456}>0$ and
  \begin{equation*}
    \rho = e^{123} + e^{456}
  \end{equation*}
  where $e^{ijk}$ is the standard abbreviation for $e^i \wedge e^j
  \wedge e^k$. In this basis, it holds $\lambda(\rho)= \nu^{\otimes
    2}$, $J_\rho e_i = e_i$ for $i \in \{ 1,2,3 \}$ and $J_\rho e_i =
  -e_i$ for $i \in \{ 4,5,6 \}$.
\end{enumerate}

The other orbit consists of all three-forms $\rho$ satisfying one of
the following equivalent properties.
\begin{enumerate}[(a)]
\item $\lambda(\rho) < 0$
\item There is a unique complex decomposable three-form $\alpha$ such
  that $\rho = \alpha + \bar \alpha$ and $\rm i (\bar \alpha \wedge
  \alpha) > 0$.
\item The stabiliser of $\rho$ in $\GL^+(V)$ is $\SL(3,\bC)$.
\item It holds $\lambda(\rho) \ne 0$ and the endomorphism $J_\rho$ is
  a complex structure on $V$.
\item There is a basis $\{ e_1, ..., e_6 \}$ of $V$ such that
  $\nu=e^{123456}>0$ and
  \begin{equation*}
    \rho = e^{135} - e^{146} - e^{236} - e^{245}.
  \end{equation*}
  In this basis, it holds $\lambda(\rho) = - 4 \nu^{\otimes 2}$,
  $J_\rho e_i = -e_{i+1}$ and $J_\rho e_{i+1}= e_{i}$ for $i \in \{
  1,3,5 \}$.
\end{enumerate}
\end{Prop}
\begin{proof}
  All properties are proved in section 2 of \cite{H2}. The only fact
  we added is the observation that $J_\rho$ is a para-complex
  structure if $\lambda(\rho)>0$ which is obvious in the standard
  basis.
\end{proof}

It is also possible to introduce a basis describing both orbits
simultaneously. Indeed, given a generic stable three-form and an
orientation, there is a basis $\{ e_1, ..., e_6 \}$ of $V$ and an $\ve
\in \{ \pm 1\}$ such that $\nu=e^{123456}>0$ and
\begin{equation}
  \label{normal3form} \rho_\ve = e^{135} + \ve (e^{146} + e^{236} + e^{245})
\end{equation}
with $\lambda(\rho) = 4 \ve \nu^{\otimes 2}$. Furthermore, it holds
$J_\rho e_i = \ve e_{i+1}$, $J_\rho e_{i+1}= e_{i}$ for $i \in \{
1,3,5 \}$ and
\begin{equation}
  \label{normalhat} J^*_{\rho_\ve} \rho_\ve = e^{246} + \ve (e^{235} + e^{145} + e^{136}).
\end{equation}
Analogies between complex and para-complex structures are elaborated
in a unified language in \cite{AC2} and \cite{SSH}. In this language,
a stable three-form always induces an $\ve$-complex structure $J_\rho$
since $J^2_\rho = \ve \rm{id}$ for the normal form $\rho_\ve$.
\bl \label{rhohatlemma} The dual of a stable three-form $\rho \in
\Lambda^3 V^*$ on an oriented six-dimensional vector space $V$ is
\begin{equation}
  \hat \rho = J_\rho^* \rho. \label{rhohat}  
\end{equation}
\el
\pf We already observed that the connected components of the
stabilisers of $\rho$ and $\hat \rho$ have to be identical. Therefore,
since the space of real three-forms invariant under
$\SL(3,\mathbb{C})$ respectively $\SL(3,\mathbb{R}) \times
\SL(3,\mathbb{R})$ is two-dimensional, we can make the ansatz \[\hat
\rho = c_1 \rho + c_2 J_\rho^* \rho\] with real constants $c_1$ and
$c_2$. Computing
\[ \frac{6}{3} \phi(\rho) \stackrel{(\ref{phieuler})}{=} \hat \rho
\wedge \rho = c_2 \, J_\rho^* \rho \wedge \rho
\stackrel{(\ref{normal3form}, \ref{normalhat})}{=} 2 c_2 \,
\phi(\rho), \] we find $c_2=1$. By
\[ d_\rho \phi (J_\rho^* \rho) = \hat \rho \wedge J_\rho^* \rho = c_1
\rho \wedge J_\rho^* \rho = - 2 c_1 \phi(\rho),\] the constant $c_1$
vanishes if the derivative of $\lambda$ (recall (\ref{phirho})) in
$\rho$ in direction of $J_\rho^* \rho$ vanishes. However, using the
normal form (\ref{normal3form}) again, we compute $\lambda(\rho + t
J_\rho^* \rho) = 4 \ve (-\ve + t^2)^2 \, (e^{123456})^{\otimes 2}$ and
the assertion follows.  \epf

A convenient way to compute the dual of $\rho$ without determining
$J_\rho$ is given by the following corollary. In fact, the corollary
explicitly shows the equivalence of the two different definitions of
$\rho \mapsto \hat \rho$ given in \cite{H1} and \cite{H2}.

\begin{Cor}
  \label{hat}
  If $\lambda(\rho)>0$ and $\rho = \alpha + \beta$ in terms of
  decomposables ordered such that $\alpha \wedge \beta >0$, the dual
  of $\rho$ satisfies $\hat \rho = \alpha - \beta$.

  If $\lambda(\rho)<0$ and $\rho$ is the real part of a complex
  decomposable three-form $\alpha$ such that $\rm i (\bar \alpha
  \wedge \alpha) > 0$, the dual of $\rho$ is the imaginary part of
  $\alpha$. In particular, the complex three-form $\a$ is a
  $(3,0)$-form w.r.t.\ $J_\rho$.
\end{Cor}
\begin{proof}
  The assertions are easily proved by comparing the claimed formulas
  for $\hat \rho$ with formula \eqref{rhohat} in the standard bases
  given in part (\ref{enum1}) of Proposition \ref{orbitsof3forms}.
\end{proof}

Finally, we note that for a fixed orientation, it holds
\begin{equation}
  \hat{\hat \rho} = - \rho \quad \mbox{and} \quad J_{\hat \rho}= - \ve J_\rho.
  \label{hathatrho}
\end{equation}

$\mathbf{k=3, n=7.}$ Given any three-form $\varphi$, we define a
symmetric bilinear form with values in $\Lambda^7 V^*$ by
\begin{equation}
  \label{bvarphi}
  b_\varphi(v,w)=(v \lrcorner\, \varphi) \wedge (w \lrcorner\, \varphi) \wedge \varphi.
\end{equation}
Since the determinant of a scalar-valued bilinear form is an element
of $(\Lambda^7 V^*)^{\otimes 2}$, we have $\mbox{det}\, b \in
(\Lambda^7 V^*)^{\otimes 9}$. If and only if $\varphi$ is stable, the seven-form
\begin{eqnarray*}
  \phi(\varphi)= (\mbox{det}\, b_\varphi)^\frac{1}{9}
\end{eqnarray*}
defines a volume form, independent of an orientation on $V$, and the
scalar-valued symmetric bilinear form
\begin{eqnarray*}
  g_\varphi = \frac{1}{\phi(\varphi)} b_\varphi
\end{eqnarray*}
is non-degenerate. Notice that $\phi(\vf) = \sqrt{\det g_\vf}$ is the
metric volume form.

It is known (\cite{B1}, \cite{Ha}) that a stable three-form defines a
multiplication ``$\cdot$'' and a vector cross product ``$\times$'' on
$V$ by the formula
\begin{equation}
  \label{stablevsO}
  \vf (x,y,z) = g_\vf(x, y \cdot z) = g_\vf (x,y \times z),
\end{equation} 
such that $(V, \times)$ is isomorphic either to the imaginary octonions
$\mathrm{Im}\, \bO$ or to the imaginary split-octonions $\mathrm{Im}\, 
\tilde \bO$. Thus, there are exactly two open orbits of stable
three-forms having isotropy groups
\begin{eqnarray}
  \label{G2stab}
  \Stab_{\GL(V)}(\varphi) \cong
  \begin{cases}
    \G_2 \subset \SO(7), & \mbox{if $g_\varphi$ is positive definite,}\\
    \G_2^* \subset \SO(3,4), & \mbox{if $g_\varphi$ is of signature (3,4)}.
  \end{cases}
\end{eqnarray}
There is always a basis $\{ e_1 ,\, ...\, , \, e_7 \}$ of $V$ such
that
\begin{equation} \label{varphistandard} \varphi = \tau e^{124} +
  \sum_{i=2}^7 e^{i \, (i+1) \, (i+3)}
\end{equation}
with $\tau \in \{ \pm 1 \}$ and indices modulo 7.  For $\tau = 1$, the
induced metric $g_\varphi$ is positive definite and the basis is
orthonormal such that this basis corresponds to the Cayley basis of
$\mathrm{Im}\, \bO$. For $\tau = -1$, the metric is of signature (3,4)
and the basis is pseudo-orthonormal with $e_1$, $e_2$ and $e_4$ being
the three spacelike basis vectors. 

The only four-forms having the same stabiliser as $\varphi$ are the
multiples of the Hodge dual $*_{g_\varphi}\varphi$, \cite[Propositions
2.1, 2.2]{B1}. Since the normal form satisfies
$g_\varphi(\varphi,\varphi) = 7$, we have by definition of the Hodge
dual $\varphi \wedge *_{g_\varphi} \varphi = 7 \, \phi(\varphi)$ and
therefore
\begin{eqnarray}
  \label{hatvsHodge1}
  \hat \varphi = \frac{1}{3} *_{g_\varphi} \varphi,
\end{eqnarray}
by comparing with (\ref{phieuler}).

\begin{lemma}
  Let $\varphi$ be a stable three-form in a seven-dimensional vector
  space $V$. Let $\beta$ be a one-form or a two-form. Then $\beta
  \wedge \varphi = 0$ if and only if $\beta=0$.
\end{lemma}
\pf For the compact case, see also \cite{Bo}. If $\beta$ is a
one-form, the proof is very easy. If $\beta$ is a two-form, we choose
a basis such that $\varphi$ is in the normal form
(\ref{varphistandard}) and $\beta = \sum_{i < j} b_{i,j} \, e^{ij}$
and compute
\begin{eqnarray*}
  \beta \wedge \varphi &=&
  (b_{2,3}-b_{1,6}) \:e^{12356}
  +(b_{2,3} - b_{4,7}) \:e^{23457}
  +(b_{1,6}+b_{4,7}) \:e^{14567}\\
  &+&(b_{5,7} \tau +b_{1,2}) \:e^{12457}
  +(b_{3,6}-b_{5,7}) \:e^{34567}
  +(b_{1,2}-b_{3,6} \tau) \:e^{12346}\\
  &-&(b_{3,7} \tau+b_{2,4}) \:e^{12347}
  +(b_{5,6} \tau+b_{2,4}) \:e^{12456}
  +(b_{3,7}+b_{5,6}) \:e^{13567}\\
  &+&(b_{2,5} - b_{4,6}) \:e^{23456}
  +(b_{4,6}-b_{1,7}) \:e^{13467}
  -(b_{2,5}+b_{1,7}) \:e^{12357}\\
  &+&(b_{4,5}+b_{2,6}) \:e^{24567}
  -(b_{1,3}+b_{2,6}) \:e^{12367}
  +(b_{4,5}+b_{1,3}) \:e^{13457}\\
  &+&(b_{3,5}+b_{6,7}) \:e^{23567}
  +(b_{1,4}-b_{3,5} \tau) \:e^{12345}
  +(b_{6,7} \tau-b_{1,4}) \:e^{12467}\\
  &+&(b_{3,4}+b_{1,5}) \:e^{13456}
  +(b_{2,7}-b_{1,5}) \:e^{12567}
  +(b_{3,4}-b_{2,7}) \:e^{23467}.
\end{eqnarray*}
The five-form is written as a linear combination of linearly independent
forms and each line contains exactly three different coefficients of
$\beta$. Inspecting the coefficient equations line by line, it is easy
to see that all coefficients of $\beta$ vanish if and only if $\beta
\wedge \varphi = 0$.  \epf

\subsection{Real forms of $\SL(3,\mathbb{C})$}
\label{SL3C}
By the following proposition, any real form of $\SL(3,\mathbb C)$ can
be written as a simultaneous stabiliser of a stable two-form and a
stable three-form.  \bp Let $V$ be a six-dimensional real vector
space. Let $\omega \in \Lambda^2 V^*$ and $\rho \in \Lambda^3 V^*$ be
stable forms which are compatible in the sense that
\begin{eqnarray}
  \label{comp01} \omega \wedge \rho &=&0.
\end{eqnarray}
Then, we have
\[ \Stab_{\GL(V)}(\rho,\omega) \cong
\begin{cases}
  \SU(p,q) \subset \SO(2p,2q) \, ,\quad p+q = 3\, , & \mbox{if $\lambda(\rho) < 0$,} \\
  \SL(3,\mathbb R) \subset \SO(3,3)\,, & \mbox{if $\lambda(\rho) > 0$,}
\end{cases} \] where $\SL(3,\mathbb R)$ is embedded in $\SO(3,3)$ such
that it acts by the standard representation and its dual,
respectively, on the maximally isotropic $\pm 1$-eigenspaces of the
para-complex structure $J_\rho$ induced by $\rho$.  \ep
\pf Let $V$ be oriented by $\phi(\omega)=\frac{1}{6}\omega^3$ and let
$J_\rho$ be the unique \mbox{(para-)} complex structure \eqref{Jrho}
associated to the three-form $\rho$ and this orientation. By $\ve \in
\{ \pm 1 \}$, we denote the sign of $\lambda(\rho)$, that is $J_\rho^2
= \ve id_V$. In the basis in which $\rho$ is in the normal form
(\ref{normal3form}), it is easy to verify that $\omega \wedge \rho =
0$ is equivalent to the skew-symmetry of $J_\rho$ with respect to
$\omega$. Equivalently, the pseudo-Euclidean metric
\begin{equation} \label{inducedmetric} g = g_{(\o,\rho)} = \ve \,
  \omega(\cdot,J_\rho\cdot),
\end{equation}
induced by $\omega$ and $\rho$, is compatible with $J_\rho$ in the
sense that $g(J_\rho \cdot, J_\rho \cdot) = -\ve g(\cdot, \cdot)$. The
stabiliser of the set of tensors $(\omega,J_\rho,g,\rho,J_\rho^*
\rho)$ satisfying this compatibility condition is well-known to be
$\SU(p,q)$ respectively $\SL(3,\bR)$.  \epf
We will call a compatible pair of stable forms $(\omega,\rho) \in
\Lambda^2 V^* \times \Lambda^3 V^*$ normalised if
\begin{equation}
  \phi(\rho) = 2\, \phi(\o) \quad \iff \quad J_\rho^* \rho \wedge \rho = \frac{2}{3}\, \omega^3 \label{normrho}.   
\end{equation}
\br By our conventions, the metric \eqref{inducedmetric} induced by a
normalised, compatible pair is of signature either $(6,0)$ or $(2,4)$
or $(3,3)$, where the first number denotes the number of spacelike
directions. We emphasise that our conventions are such that \[\omega =
g(\cdot,J_\rho \cdot). \] This sign choice turned out to be necessary in
order to achieve that $\phi(\rho)$ is indeed a positive multiple of
$\phi(\o)$ in the positive definite case. \er
Sometimes it is convenient to have a unified adapted basis. For a
compatible and normalised pair $(\omega,\rho)$, there is always a
pseudo-orthonormal basis $\{ e_1, \, ... \, ,e_6 \}$ of $V$ with dual
basis $\{ e^1, \, ... \, ,e^6 \}$ such that $\rho = \rho_\ve$ is in
the normal form (\ref{normal3form}) and
\begin{eqnarray}
  \omega &=& \tau (e^{12} + e^{34}) + e^{56}  \label{normalomega}
\end{eqnarray}
for $(\ve, \tau) \in \{ (-1,1),(-1,-1),(1,1) \}$. The signature of the
induced metric with respect to this basis is
\begin{eqnarray}
  \label{signnormal}
  (\tau , -\ve \tau, \tau, - \ve \tau, 1, -\ve) = 
  \begin{cases} 
    (+,+,+,+,+,+) & \mbox{for $\ve = -1$ and $\tau =1$,}\\
    (-,-,-,-,+,+) & \mbox{for $\ve = -1$ and $\tau =-1$,}\\
    (+,-,+,-,+,-) & \mbox{for $\ve = 1$ and $\tau =1$,}
  \end{cases}
\end{eqnarray}
and we have
\begin{eqnarray*}
  \Stab_{\GL(6,\mathbb{R})}(\omega,\rho) \cong
  \begin{cases}
    \SU(3) \subset \SO(6) & \mbox{for $\ve = -1$ and $\tau =1$,}\\
    \SU(1,2) \subset \SO(2,4) & \mbox{for $\ve = -1$ and $\tau =-1$,}\\
    \SL(3,\mathbb{R}) \subset \SO(3,3) & \mbox{for $\ve = 1$}.
  \end{cases}
\end{eqnarray*}
For instance, the following observation is easily verified using the
unified basis.
\begin{lemma}
  Let $(\o,\rho)$ be a compatible and normalised pair of stable forms
  on a six-dimensional vector space. Then, the volume form $\phi(\o)$
  is in fact a metric volume form w.r.t.\ to the induced metric
  $g=g_{(\o,\rho)}$ and the corresponding Hodge dual of $\o$ and
  $\rho$ is
  \begin{equation}
    \label{hatvsHodge2}
    *_g \o = -\ve \hat \o \, , \qquad *_g \rho = - \hat \rho
  \end{equation}
\end{lemma}

\subsection{Relation between real forms of $\SL(3,\mathbb{C})$ and $\G_2^\mathbb{C}$}
\label{6vs7}
The relation between stable forms in dimension six and seven
corresponding to the embedding $\SU(3) \subset \G_2$ is well-known. We
extend this relation by including also the embeddings $\SU(1,2)
\subset \G_2^*$ and $\SL(3,\bR) \subset \G_2^*$ as follows.  

\bp Let $V = W \oplus L$ be a seven-dimensional vector space
decomposed as a direct sum of a six-dimensional subspace $W$ and a
line $L$. Let $\alpha$ be a non-trivial one-form in the annihilator
$W^0$ of $W$ and $(\omega,\rho) \in \Lambda^2 L^0 \times \Lambda^3
L^0$ a compatible and normalised pair of stable forms inducing the
scalar product $h=h_{(\o,\rho)}$ given in \eqref{inducedmetric}.  Then, the three-form
$\varphi \in \Lambda^3 V^*$ defined by
\begin{equation}
  \label{varphi}
  \varphi = \omega \wedge \alpha + \rho 
\end{equation}
is stable and induces the scalar product
\begin{equation}
  \label{gvarphi}
  g_\varphi = h - \ve \alpha \cdot \alpha
\end{equation}
where $\ve$ denotes the sign of $\lambda(\rho)$ such that $J_\rho^2 =
\ve \rm{id}$.  The stabiliser of $\vf$ in $\GL(V)$ is
\begin{equation*}
  \Stab_{\GL(V)}(\varphi) \cong 
  \begin{cases}
    \G_2 & \mbox{for $\ve = -1$ and positive definite h,}\\
    \G_2^* & \mbox{otherwise.}\\
  \end{cases}
\end{equation*}
\ep 
\pf We choose a basis $\{ e_1 ,\, ... \,,\, e_6 \}$ of $L^0$ such
that $\omega$ and $\rho$ are in the generic normal forms
(\ref{normal3form}) and (\ref{normalomega}). With $e^7 = \alpha$, we
have
\begin{eqnarray}
  \label{g2normal}
  \varphi = \tau (e^{127} + e^{347}) + e^{567} + e^{135} + \ve (e^{146} + e^{236} + e^{245}). 
\end{eqnarray}
The induced bilinear form \eqref{bvarphi} turns out to be
\begin{eqnarray*}
  b_{\varphi}(v,w) =( - \ve \tau v^1w^1  + \tau v^2w^2  - \ve \tau v^3w^3  + \tau v^4w^4  - \ve v^5w^5  + v^6w^6  + v^7w^7) e^{1234567}
\end{eqnarray*}
for $v = \sum v^i e_i$ and $w = \sum w^i e_i$. Hence, the three-form
$\varphi$ is stable for all signs of $\ve$ and $\tau$ and its
associated volume form is
\begin{eqnarray*}
  \phi(\varphi) = (\mbox{det}\, b_{\varphi})^{\frac{1}{9}} = - \ve\, e^{1234567}.
\end{eqnarray*}
The formula \eqref{gvarphi} for the metric $g_\varphi$ induced by
$\varphi$ follows, since the basis $\{ e_1 ,\, ... \,,\, e_7 \}$ of
$V$ is pseudo-orthonormal with respect to this metric of signature
\begin{eqnarray}
  \label{gvfnormal}
  (\tau, -\ve \tau, \tau , -\ve \tau, 1, - \ve, -\ve) =
  \begin{cases} 
    (+,+,+,+,+,+,+) & \mbox{for $\ve = -1$ and $\tau =1$},\\
    (-,-,-,-,+,+,+) & \mbox{for $\ve = -1$ and $\tau =-1$},\\
    (+,-,+,-,+,-,-) & \mbox{for $\ve = 1$ and $\tau =1$.}
  \end{cases}
\end{eqnarray}
The assertion on the stabilisers now follows from \eqref{G2stab}.
\epf \bl Under the assumptions of the previous proposition, the dual
four-form of the stable three-form $\varphi$ is
\begin{equation}
  \label{Hodgevf}
  3 \hat \varphi = *_\varphi \varphi = - \,\ve\, (\alpha \wedge \hat \rho + \hat \o )= \ve \alpha \wedge *_h \rho + *_h \o , 
\end{equation}
where $*_\varphi$ denotes the Hodge dual with respect to the metric
$g_\varphi$ and the orientation induced by $\phi(\varphi)$.  \el \pf
In the basis of the previous proof, the Hodge dual of $\varphi$ is
\[ *_\varphi \varphi = -\ve \tau ( e^{3456} + e^{1256} ) - \ve
e^{1234} + \ve\, e^{2467} + e^{2357} + e^{1457} + e^{1367}. \] The
second equality follows when comparing this expression with $\ve ( e^7
\wedge \hat \rho + \frac{1}{2} \omega^2)$ in this basis using
\eqref{normalhat} and (\ref{normalomega}). The first and the third
equality are just the formulas (\ref{hatvsHodge1}) and
(\ref{hatvsHodge2}), respectively.  \epf
The inverse process is given by the following construction.
\begin{Prop}
  \label{7to6_prop}
  Let $V$ be a seven-dimensional real vector space and $\vf \in \L^3
  V^*$ a stable three-form which induces the metric $g_\vf$ on
  $V$. Moreover, let $n \in V$ be a unit vector with $g_\vf(n,n)=-\ve
  \in \{ \pm 1\}$ and let $W= n^\perp$ denote the orthogonal
  complement of $\bR \!\cdot\! n$.  Then, the pair $(\o,\rho) \in
  \Lambda^2 W^* \times \Lambda^3 W^*$ defined by
  \begin{eqnarray}
    \label{7to6}
    \o = \, n \,\inter\, \vf \, , \qquad \rho = \vf_{|W} \,,
  \end{eqnarray}
  is a pair of compatible normalised stable forms.  The metric $h =
  h_{(\o,\rho)}$ induced by this pair on $W$ satisfies
  $h=(g_\vf)_{|W}$ and the stabiliser is
  \begin{eqnarray*}
    \Stab_{\GL(W)}(\omega,\rho) \cong
    \begin{cases}
      \SU(3), & \mbox{if $g_\vf$ is positive definite,}\\
      \SU(1,2), & \mbox{if $g_\vf$ is indefinite and $\ve = -1$,}\\
      \SL(3,\mathbb{R}), & \mbox{if $\ve = 1$}.
    \end{cases}
  \end{eqnarray*}  
  When $(V,\vf)$ is identified with the imaginary octonions,
  respectively, the imaginary split-octonions, by \eqref{stablevsO}, the
  $\ve$-complex structure induced by $\rho$ is given by
  \begin{equation}
    \label{Jcross}
    J_\rho v = - n\cdot v =  - n \times v \qquad \qquad \mbox{for $v \in V$}.  
  \end{equation}
\end{Prop}
\begin{proof}
  Due to the stability of $\vf$, we can always choose a basis $\{
  e_1,\dots,e_7 \}$ of $V$ with $n=e_7$ such that $\vf$ is given by
  \eqref{g2normal} where $\ve = - g_\vf(n,n)$ and $\tau \in \{ \pm 1
  \}$ depends on the signature of $g_\vf$.  As this basis is
  pseudo-orthonormal with signature given by \eqref{gvfnormal}, the
  vector $n$ has indeed the right scalar square and $\{ e_1,\dots,e_6
  \}$ is a pseudo-orthonormal basis of the complement $W=n^\perp$.
  Since the pair $(\o,\rho)$ defined by \eqref{7to6} is now exactly in
  the generic normal form given by (\ref{normal3form}) and
  (\ref{normalomega}), it is stable, compatible and normalised and the
  induced endomorphism $J_\rho$ is an $\ve$-complex structure.  The
  identity $h=(g_\vf)_{|W}$ for the induced metric $h_{(\o,\rho)}$
  follows from comparing the signatures \eqref{gvfnormal} and
  \eqref{signnormal} and the assertion for the stabilisers is an
  immediate consequence. Finally, the formula for the induced
  $\ve$-complex structure $J_\rho$ is another consequence of
  $g=(g_\vf)_{|W}$ since we have
  \[ g_\vf (x,n \times y) \stackrel{\eqref{stablevsO}}{=} \vf (x,n,y)
  = - \o(x,y) = - h(x,J_\rho y) \] for all $x,y \in W$.
\end{proof}
Notice that, for a fixed metric $h$ of signature $(2,4)$ or $(3,3)$,
the compatible and normalised pairs $(\o,\rho)$ of stable forms
inducing this metric are parametrised by the homogeneous spaces
$\SO(2,4) /\,\SU(1,2)$ and $\SO(3,3) /\,\SL(3,\bR)$, respectively.
Thus, the mapping $(\o,\rho) \mapsto \vf$ defined by formula
(\ref{varphi}) yields isomorphisms
\begin{equation*}
  \frac{\SO(2,4)}{\SU(1,2)} \cong \frac{\SO(3,4)}{\G_2^*}\,, \qquad 
  \frac{\SO(3,3)}{\SL(3,\bR)} \cong \frac{\SO(3,4)}{\G_2^*}\,,
\end{equation*}
since the metric $h$ completely determines the metric
$g_\vf$ by the formula (\ref{gvarphi}).

\subsection{Relation between real forms of  $\G_2^\mathbb{C}$ and $\mathrm{Spin}(7,\mathbb{C})$}
\label{7vs8}
It is possible to extend this construction to dimension eight as
follows. Starting with a stable three-form $\vf$ on a
seven-dimensional space $V$, we can consider the four-form
\begin{equation}
  \label{spin4form}
  \Phi = e^8 \wedge \varphi + *_\varphi \varphi.
\end{equation}
on the eight-dimensional space $V \op \bR e_8$. Although the four-form
$\Phi$ is not stable, it is shown in \cite{B1} that it induces the
metric 
\begin{equation}
  \label{gPhi}
  g_\Phi = g_\vf + (e^8)^2  
\end{equation}
 on $V \op \bR e_8$ and that its stabiliser is
\begin{eqnarray*}
  \Stab_{\GL(V \op \bR e_8)}(\Phi) \cong
  \begin{cases}
    \Spin(7) \subset \SO(8), & \mbox{if $g_\vf$ is positive definite,}\\
    \Spin_0(3,4) \subset \SO(4,4), & \mbox{if $g_\vf$ is indefinite.}
  \end{cases}
\end{eqnarray*}  
The index ``$0$'' denotes, as usual, the connected
component. Starting conversely with a four-form $\Phi$ on $V \op \bR
e_8$ such that its stabiliser in ${\GL(V \op \bR e_8)}$ is
isomorphic to $\Spin(7)$ or $\Spin_0(3,4)$, the process can be
reversed by setting $\vf = e_8 \, \inter \Phi$. As before, the metric induced
by $\Phi$ on $V \op \bR e_8$ is determined by the metric $g_\vf$
induced by $\vf$ on $V$. Thus, the indefinite
analogue of the well-known isomorphisms
\begin{equation*}
  \bR \mathbb{P}^7 \cong \frac{\SO(6)}{\SU(3)} 
  \cong \frac{\SO(7)}{\G_2} \cong \frac{\SO(8)}{\Spin(7)}
\end{equation*}
is given by
\begin{equation}
  \label{vierhomogene}
  \frac{\SO(2,4)}{\SU(1,2)} 
  \cong \frac{\SO(3,3)}{\SL(3,\bR)}
  \cong \frac{\SO(3,4)}{\G_2^*}
  \cong \frac{\SO(4,4)}{\Spin_0(3,4)}.
\end{equation}

\section{Hitchin's flow equations}

\subsection{Half-flat structures and parallel $\G_2^{(*)}$-structures}
Now we want to put the algebraic structures considered in the previous
section onto smooth manifolds. This is best done in terms of
reductions of the bundle of frames of the manifold. This bundle has
$\GL(n,\rr)$ as its structure group if $n$ is the dimension of the
manifold. A subbundle whose structure group is a subgroup $G$ of
$\GL(n,\rr)$ is a called a reduction of the frame bundle, or a
$G$-{\em structure}.  For example, if $G\subset \Og(p,q)$, for
$p+q=n$, the reduction determines a pseudo-Riemannian metric of
signature $(p,q)$ and the distinguished frames are orthonormal with
respect to this metric.  If $G\subset \Og(p,q)$, again with $p+q=n$,
then a $G$-structure is called {\em parallel} if the $G$-subbundle is
invariant under the parallel transport defined by the Levi-Civita
connection of the corresponding metric. This is equivalent to the
property that the holonomy group of the Levi-Civita is contained in
$G$.

In the following we will consider $G$-structures that are given by the
groups described in the previous sections.  According to the notations
given there, we denote by $H^{\ve,\tau}$ a real form of $\SL(3,\ccc)$
and $G^{\ve, \tau}$ the corresponding real form of $\G_2^\ccc$ in
which $H^{\ve,\tau}$ is embedded, i.e. $H^{-1,1}=\SU(3)\subset
\SO(6)$, $H^{-1,-1}=\SU(1,2)\subset \SO(2,4)$,
$H^{1,1}=\SL(3,\rr)\subset \SO(3,3)$, $G^{-1,1}=\G_2\subset \SO (7)$,
and $G^{-1,-1}=G^{1,1}=\G_2^*\subset \SO(3,4)$. We will also use the
notation $\G_2^{(*)}$ as a shorthand for ``$\G_2$ respectively
$\G_2^*$''.

An $H^{\ve,\tau}$-structure is equivalent to a pair of everywhere
stable forms $\w\in \W^2M$ and $\rho\in\W^3M$ on $M$, considered up to
rescaling of $\rho$ by a nonzero constant, that satisfy the
compatibility condition
\begin{equation}
  \rho\wedge\w=0\label{comp1}
\end{equation}
corresponding to (\ref{comp01}) and in addition
\begin{align}
  \label{comp2} 
  \phi(\rho) = c\, \phi(\omega)\,,  \qquad \mbox{i.e.\ 
  $J_\rho^* \rho \wedge \rho = \frac{1}{3}\, c\, \omega^3$,}
\end{align}
for a positive real constant $c$.  Indeed, if an
$H^{\ve,\tau}$-structure is given, these forms are obtained by
applying the formulae (\ref{normal3form}) and (\ref{normalomega}) to
one of the frames of the $H^{\ve,\tau}$-structure. By construction,
the stable forms then satisfy (\ref{comp2}) with $c=2$.

On the other hand, if $\omega \in \W^2 M$ and $\rho \in \W^3 M$ are
everywhere stable and satisfy (\ref{comp1}) and (\ref{comp2}), we can
find a local frame, in which they are in normal form after rescaling
$\rho$ by a constant. This frame then determines the
$H^{\ve,\tau}$-structure.

Note that stable forms define an $H^{\ve,\tau}$-structure, even if they
only satisfy (\ref{comp1}) but not the second compatibility condition
(\ref{comp2}). In this case $\rho$ can always be rescaled by a smooth
function such that (\ref{comp2}) holds. When we say that the pair of
stable forms defines an $H^{\ve,\tau}$-structure,
we will always assume that {\em both} compatibility conditions are
satisfied.  We will call the $H^{\ve,\tau}$-structure {\em normalised}
if $c=2$. This seems to be a common normalisation for
$\SU(3)$-structures in the literature.

Furthermore, one can show that the $H^{\ve,\tau}$-structure is
parallel if and only if $\rho$, $\hat\rho$, and $\w$ are closed. The
proof of this fact given in \cite[p. 567]{H2} generalises to $\SU
(1,2)$-structures and also to $\SL (3,\rr)$-structures, in the latter
case using Frobenius' Theorem instead of the Newlander-Nirenberg
Theorem. In all cases the parallel $H^{\ve,\tau}$-structure is
equivalent to $M$ being a Ricci-flat (para-)K\"ahler manifold.

Now we consider a weaker condition, that will turn out to be related
to parallel $\G_2^{(*)}$-structures.  \bde An $H^{\ve,\tau}$-structure
$(\rho, \w)$ is called {\em half-flat} if
\begin{eqnarray}\label{drho}
  d\rho & = &0\\
  \label{dsigma}
  d\sigma & = &0,
\end{eqnarray}
where $2\s=\w^2$.  \ede

Similarly, a smooth seven-manifold admits a $\G_2^{(*)}$-structure if
and only if there is a stable three-form $\vf$. Again, this structure is
parallel if and only if $\vf$ is closed and co-closed, i.e. $d\vf =
d\ast\vf=0$, where $\ast$ denotes the Hodge operator with respect to
the metric induced by the $\G_2^{(*)}$-structure. For a proof in both
cases see \cite[Theorem 4.1]{gray69}.

Note that any orientable hypersurface in a manifold with $\G_2$- or
$\G_2^{*}$-structure admits an $H^{\ve,\tau}$-structure by the
algebraic construction described in Proposition \ref{7to6_prop}.  If
the $\G_2^{(*)}$-structure $\vf$ is parallel,
the induced $H^{\ve,\tau}$-structure is half-flat due to equations
\eqref{varphi} and \eqref{Hodgevf}.  For the various results on the
$\SU(3)$-structures on hypersurfaces in $\G_2$-structures, we refer to
\cite{calabi58}, \cite{cabrera06} and references
therein.

On the other hand, certain one-parameter families of half-flat
structures define parallel $\G_2^{(*)}$-structures.
\bs \label{parallelprop} Let $H^{\ve,\tau}$ be a real form of
$\SL(3,\ccc)$, $G^{\ve, \tau}$ the corresponding real form of
$\G_2^\ccc$ and $(\rho,\w)$ a one-parameter family of
$H^{\ve,\tau}$-structures on a six-manifold $M$ with a parameter $t$
from an interval $I$. 
Then, the three-form
\[ \vf = \w\wedge dt +\rho \] defines a parallel
$G^{\ve,\tau}$-structure on $M\times I$ if and only if the
$H^{\ve,\tau}$-structure $(\rho,\w)$ is half-flat for all $t$ and
satisfies the following evolution equations
\begin{eqnarray}
  \dot{\rho}&=&d\w\label{ev1} 
  \\
  \dot{\s}&=&d\hat{\rho}\label{ev2} 
\end{eqnarray}
with $\s=\frac{1}{2}\w^2$.  \es \bprf Let $(\rho,\w)$ be an
$H^{\ve,\tau}$-structure and $\vf = \w\wedge dt +\rho$ a stable
three-form on $\check{M}:=M\times I$.  By \eqref{Hodgevf}, the
Hodge-dual of $\vf$ is given by
\[\ast \vf= \ve\left(\hat\rho\wedge dt-\s\right).\]
Denoting by $\check d$ the differential on $\check M$ and by $d$ the
differential on $M$ we calculate
\begin{eqnarray}
  \label{dfi}
  \check{d}\vf
  &=& d\w \wedge dt +dt \wedge\dot\rho + d\rho\ = \ (d\w - \dot\rho )\wedge dt +d\rho\\
  \label{dsternfi}\check{d}\ast \vf
  &=& \ve\left( d\hat\rho\wedge dt - dt\wedge \dot \s - d\s \right) 
  \ =\ \ve (d\hat\rho-\dot\s)\wedge dt - \ve d\s
\end{eqnarray}
Thus, $\vf$ defines a parallel $G^{\ve,\tau}$-structure if and only if
the evolution equations (\ref{ev1}) and (\ref{ev2}) and the half-flat
equations are satisfied. \eprf 
The evolution equations \eqref{ev1} and \eqref{ev2} are the {\em
  Hitchin flow equations}, as found in \cite{H1} for
$\SU(3)$-structures, applied to $H^{\ve,\tau}$-structures. Their
solutions $(\rho, \w)$, called {\em Hitchin flow}, have to satisfy
possibly dependent conditions in order to yield a parallel
$\G_2^{(*)}$-structure: the evolution equations and the compatibility
equations for the family of half-flat structures.  The following
theorem shows that the evolution equations together with an initial
condition already ensure that the family consists of half-flat
structures. A special version of this theorem was proved in \cite{H1}
under the assumption that $M$ is compact and that $H=\SU(3)$.

\btheo\label{paralleltheo} Let $(\rho_{0}, \w_{0})$ be a half-flat
$H^{\ve,\tau}$-structure on a six-manifold $M$. Furthermore, let
$(\rho, \w) \in \W^3M \times \W^2M$ be a one-parameter family of
stable forms with parameters from an interval $I$ satisfying the
evolution equations \eqref {ev1} and \eqref{ev2}. If
$(\rho(t_0),\w(t_0))=(\rho_0,\w_0)$ for a $t_0 \in I$, then $(\rho,
\w)$ is a family of half-flat $H^{\ve,\tau}$-structures. In
particular, the three-form
\begin{equation}
  \label{varphidt}
  \vf = \w\wedge dt +\rho   
\end{equation}
defines a parallel $G^{\ve,\tau}$-structure on $M\times I$ and the
induced metric
\begin{equation}
  \label{gvarphidt}
  g_\vf = g(t) - \ve dt^2,  
\end{equation}
has holonomy contained in $G^{\ve,\tau}$, where $g = g(t)$ is the
family of metrics on $M$ associated to $(\rho, \w)$.  \etheo

\bprf Differentiating the evolution equations (\ref{ev1}) and
(\ref{ev2}) gives $d\dot{\rho}=d\dot{\s}=0$. The initial condition for
$t_0$ was that $(\rho_0, \w_0)$ is half-flat. This implies \be
d\rho&=&0
\\
d\s&=&0 \ee for all $t\in I$. Hence, in order to obtain a family of
half-flat structures we have to verify that the compatibility
condition (\ref{comp1}) holds for all $t\in I$. 
\blem \label{lielemma} Let $M$ be a six-manifold with
$H^{\ve,\tau}$-structure $(\rho, \w)$, $\phi: \W^3M\rightarrow \W^6 M$
defined pointwise by the map $\phi:\Lambda^3 T^*_pM\rightarrow
\Lambda^6 T^*_p M$ given in Proposition \ref{invmap} and $\hat\rho$
defined by $d\phi_\rho (\xi)=\hat\rho\wedge \xi$ for all $\xi\in
\W^3M$. If $\cal L_X$ denotes the Lie derivative, then
\[ \cal L_X (\phi(\rho))= \hat{\rho}\wedge \cal L_X\rho.\] \elem 
\bprf First note that the $\GL(n,\rr)$-equivariance of the map
$\phi:\Lambda^3T^*_pM\rightarrow \Lambda^6T^*_pM$ implies that the
corresponding map $\phi:\W^3M\rightarrow \W^6M$ is equivariant under
diffeomorphisms. Indeed, if $\psi$ is a (local) diffeomorphism of $M$
we get that \[\psi^*(\phi(\rho))=\phi(\psi^*\rho).\] Let $\psi_t$ be
the flow of the vector field $X$. Then the Lie derivative is given by
\[
\cal L_X (\phi(\rho))= \frac{d}{dt}\left(
  \psi_t^*\phi(\rho)\right)|_{t=0} = \frac{d}{dt}\phi(
\psi_t^*\rho)|_{t=0} = d\phi_\rho( \cal L_X\rho),
\]
implying the statement.  \eprf
\begin{lemma}
  \label{formel1lemma}
  A stable three-form $\rho \in \Omega^3M$ on a six-manifold satisfies
  for any $X \in \mathfrak{X}(M)$
  \begin{eqnarray}
    \label{formel1}
    \hat \rho_X \wedge \rho &=& - \hat \rho \wedge \rho_X\,, \\
    \label{formel2}
    (d \hat \rho)_X \wedge \rho &=& \hat \rho \wedge
    (d \rho)_X\,,
  \end{eqnarray}
  where $\rho_X$ denotes the interior product of $X$ with the form $\rho$.
\end{lemma}
\begin{proof}
  In order to verify the first identity, we can assume that
  $\rho=\rho_p$ is a stable three-form on $V=T_pM$ and $X \in V$ for a
  $p\in M$.  If $\lambda(\rho)<0$, the stabiliser $\SL (3,\ccc)$ of
  $\rho$ in $\GL^+(V)$ acts transitively on $V \setminus \{ 0\}$.  If
  $\lambda(\rho)>0$, we can decompose $V$ in the $\pm 1$-eigenspaces
  $V^\pm$ of $J_\rho$. The stabiliser $\SL(3,\rr)\times \SL(3,\rr)$ of
  $\rho$ in $\GL^+(V)$ acts transitively on the dense open subset $V^+
  \setminus \{ 0\} \times V^- \setminus \{ 0\} \subset V$ and there is
  an automorphism exchanging $V^+$ and $V^-$ which stabilises $\rho$.
  Thus, it suffices to verify the first identity for the normal form
  \eqref{normal3form}, \eqref{normalhat} and $X=e_1$, which is easy.

  For the second identity, using Lemma \ref{lielemma} in the second
  step, we compute 
  \be (d \hat \rho)_X \wedge \rho - \hat{\rho}\wedge(d\rho)_X
  &=& - d \hat \rho \wedge \rho_X + \hat \rho \wedge d (\rho_X) - \hat{\rho}\wedge \mathcal L_X\rho \\
  &
  =& - d(\hat{\rho}\wedge\rho_X) - \mathcal L_X (\phi(\rho)) \\
  &=& - d(\hat{\rho}\wedge\rho_X + \phi(\rho)_X) \\
  &\stackrel{(\ref{phieuler})}{=}& - \frac{1}{2}
  d(\hat{\rho}\wedge\rho_X + \hat \rho_X \wedge \rho).  \ee Hence, the
  first identity \eqref{formel1} implies \eqref{formel2}.
\end{proof}
Using this lemma, we calculate the $t$-derivative of the six-form
$\w_X\wedge\w\wedge\rho=\s_X\wedge\rho$ for any vector field $X$:
\be \ddt (\s_X \wedge
\rho)&=&\dot{\s}_X\wedge \rho + \s_X\wedge \dot{\rho}
\\
&\stackrel{(\ref{ev1} ),(\ref{ev2} )}{=}& (d\hat{\rho})_X\wedge \rho +
\s_X\wedge d\w
\\
&\stackrel{(\ref{formel2})}{=}& \hat \rho \wedge (d \rho)_X
+\w_X\wedge\w\wedge d\w\\
&\stackrel{\eqref{drho},\eqref{dsigma}}{=}& 0.  \ee
Together with the initial condition $\o_0\wedge\rho_0=0$ this implies
that $\s_X\wedge \rho=0$ for all $t\in I$ and for all vector fields
$X$. Since $\w$ is non degenerate, the product of any one-form with
$\w\wedge \rho$ vanishes and thus, the compatibility condition
$\w\wedge \rho=0$ holds for all $t$.

The preservation of the normalisation \eqref{comp2} in time is shown
in \cite{H1}, in the final part of the proof of Theorem 8. The idea is
to compute the second derivative of the volume form assigned to a
stable three-form. In fact, the proof holds literally for all
signatures since all it uses is the first compatibility condition we
have just proved.  \eprf
\bfolg
\label{Cauchy}
Let $M$ be a real analytic six-manifold with a half-flat
$H^{\ve,\tau}$-structure that is given by a pair of analytic stable
forms $(\omega_0,\rho_0)$.
\begin{enumerate}[(i)]
\item Then, there exists a unique maximal solution $(\omega,\rho)$ of
  the evolution equations \eqref{ev1}, \eqref{ev2} with initial value
  $(\omega_0,\rho_0)$, which is defined on an open neighbourhood
  $\Omega \subset \rr\times M$ of $\{0\}\times M$.  In particular,
  there is a parallel $G^{\ve,\tau}$-structure on $\Omega$.
\item Moreover, the evolution is natural in the sense that, given a
  diffeomorphism $f$ of $M$, the pullback $(f^*\omega,f^*\rho)$ of the
  solution with initial value $(\omega_0,\rho_0)$ is the solution of
  the evolution equations for the initial value
  $(f^*\omega_0,f^*\rho_0)$.

  In particular, if $f$ is an automorphism of the initial structure
  $(\omega_0,\rho_0)$, then, for all $t \in \rr$, $f$ is an
  automorphism of the solution $(\omega (t),\rho (t))$ defined on the
  (possibly empty) open set $ U_t=\{p\in M\mid (t,p)\in \W\text{ and }
  (t,f(p))\in \W\}$.
\item Furthermore, assume that $M$ is compact or a homogeneous space
  $M=G/K$ such that the $H^{\ve,\tau}$-structure is
  $G$-invariant. Then there is a unique maximal interval $I\ni 0$ and
  a unique solution $(\omega,\rho)$ of the evolution equations
  \eqref{ev1}, \eqref{ev2} with initial value $(\omega_0,\rho_0)$ on
  $I\times M$. In particular, there is a parallel
  $G^{\ve,\tau}$-structure on $I\times M$.
\end{enumerate}
\efolg
\begin{proof}
  If the manifold and the initial structure $(\omega_0,\rho_0)$ are
  analytic, there exists a unique maximal solution of the evolution
  equations on a neighbourhood $\Omega$ of $M \times \{ 0 \}$ in $M
  \times \bR$ by the Cauchy-Kovalevskaya theorem.
  The naturality of the solution is an immediate consequence of the
  uniqueness due to the naturality of the exterior derivative.  If $M$
  is compact, there is a maximal interval $I$ such that the solution
  is defined on $M \times I$. The same is true for a homogeneous
  half-flat structure $(\omega_0,\rho_0)$ as it is determined by
  $(\omega_0,\rho_0)_{|p}$ for any $p \in M$.
\end{proof}
We remark that, for a homogeneous half-flat structure
$(\omega_0,\rho_0)$, the evolution equations reduce to a system of
ordinary differential equations due to the naturality assertion of the
corollary. This simplification will be used in Section
\ref{examplesH3H3} to construct metrics with holonomy equal to $\G_2$
and $\G_2^*$.

\subsection{Remark on completeness: geodesically complete conformal $\G_2$-metrics}
The $\G_2^{(*)}$-metrics arising from the Hitchin flow on a
six-manifold $N$ are of the form $(I\times N, dt^2+ g_t)$ with an open
interval $I=(a,b)$ and a family of Riemannian
metrics $g_t$ depending on $t\in I $ (formula \eqref{gvarphidt} in Theorem
\ref{paralleltheo}). As curves of the form $t\mapsto
(t, x)$ are geodesics for this metric, they are obviously geodesically
incomplete if $a$ or $b \in \rr$.

For the {\em Riemannian} case and {\em compact} manifolds $N$, we
shall explain how one easily obtains {\em complete} metrics by a
conformal change of the $\G_2$-metric.  
\blem Let $N$ be a compact manifold with a family $g_r$ of Riemannian
metrics. Then the Riemannian metric on $\rr\times N$ defined by
$h=dr^2+g_r$ is geodesically complete.  \elem \bprf Denote by $d$ the
distance on $\rr\times N$ induced by the Riemannian metric
$h=dr^2+g_r$ and by $d_r$ the distance on $N$ induced by $g_r$.  For a
curve $\gamma $ in $M=\rr\times N$ we have that the length of
$\gamma(t)=(r(t),x(t))$ satisfies
\[ \ell(\gamma) =\int_0^1 \sqrt{ \dot{r}(t)^2
  +g_{r(t)}(\dot{x}(t),\dot{x}(t)) }dt \ \ge \ \int_0^1 | \dot{r}(t) |
dt\ \ge \ |r(1)-r(0)|. \]
As the distance of two points $p=(r,x)$ and $q=(s,y)$ is defined as
the infimum of the lengths of all curves joining them,
this inequality implies that
\begin{equation}
  \label{dist1}
  d(p,q)\ge |r-s|.
\end{equation}
Note also that a curve $\gamma (t) = \left((s-r)t+r, x\right)$ joining
$p=(r,x)$ and $q=(s,x)$ in $\rr\times \{x\}$ has length
$\ell(\gamma)=|r-s| $ and thus, for such $p$, $q$ we get that
$d(p,q)=|r-s|$.  On the other hand, for $p=(r,x)$ and $q=(r,y)$ with
the same $\rr$-projection $r$ we only get that $d(p,q)\le d_r(x,y)$.

Since $h$ has Riemannian signature we can use the Hopf-Rinow Theorem
and consider a Cauchy sequence $p_n=(r_n,x_n)\in \rr\times N$
w.r.t.\ the distance $d$.  Equation \eqref{dist1} then implies that the
sequence $r_n$ is a Cauchy sequence in $\rr$. Hence, $r_n$ converges
to $r\in \rr$.  Since $N$ is compact, the sequence $x_n$ has a
subsequence $x_{n_k}$ converging to $x\in N$. For $p=(r,x)$ and
$q_{n_k}:=(r,x_{n_k})$ the triangle inequality implies that
\[ d(p,p_{n_k})\ \le\ d\left(p,q_{n_k}\right)+ d\left( q_{n_k},
  p_{n_k}\right) \ \le\ d_r(x,x_{n_k}) + d\left( q_{n_k},
  p_{n_k}\right)\ =\ d_r(x,x_{n_k})+ |r-r_{n_k}|.\] Hence, $p_{n_k}$
converges to $p$. As $p_n$ was a Cauchy sequence, we have found $p$ as
a limit for $p_n$. By the Theorem of Hopf and Rinow, $M$ is
geodesically complete.  \eprf
The consequence of the lemma is
\bs \label{complete} Let $(M=I\times N,h= dt^2+ g_t)$ be a Riemannian
metric on a product of an open interval $I$ and a compact manifold
$N$. Then $(M,h)$ is globally conformally equivalent to a metric on
$\rr\times N$ that is geodesically complete. The scaling factor
depends only on $t\in I$ and is determined by a diffeomorphism $\vf :
\rr\to I$.  \es
\bprf Let $\vf :\rr\to I$ be a diffeomorphism with inverse
$r=\vf^{-1}$. Changing the coordinate $t$ to $r$, the metric $h$ on
$I\times N$ can be written as
\[ h\ =\ \left(\vf'(r) dr\right)^2 + g_{\vf(r)}\ =\ \vf'(r)^2\left(
  dr^2 + \frac{1}{ \vf'(r)^2} g_{\vf(r)}\right).\] Hence, $h$ is
globally conformally equivalent to the metric $ dr^2 +
\frac{1}{\vf'(r)^2} g_{\vf(r)}$ on $\rr\times N$. By the lemma, this
metric is geodesically complete.  \eprf
Regarding the solution of the Hitchin flow equations, using Theorem
\ref{paralleltheo}, Corollary \ref{Cauchy}, and Proposition
\ref{complete} we obtain the following consequence.

\bfolg\label{completecoroll} Let $M$ be a compact analytic
six-manifold with half-flat $\SU(3)$-structure given by analytic
stable forms $(\rho_0,\w_0)$. Then there is a complete metric on
$\rr\times M$ that is globally conformal to the parallel $\G_2$-metric
obtained by the Hitchin flow.  \efolg 
In Example \ref{completeexample} of Section \ref{examplesH3H3} we will
construct explicit examples of this type.  Finally, note that due to
the Cheeger-Gromoll splitting Theorem, see for example \cite[Theorem
6.79]{Besse87}, one cannot expect to obtain by the Hitchin flow
irreducible $\G_2$-metrics that are complete without allowing
degenerations of $g_t$.

\subsection{Nearly half-flat structures and nearly parallel $\G_2^{(*)}$ -structures}
A $\G_2^{(*)}$-structure $\varphi$ on a seven-manifold $N$ is called
\emph{nearly parallel} if
\begin{eqnarray}
  d \varphi = \mu *_\varphi \varphi
\end{eqnarray}
for a constant $\mu \in \mathbb R^*$. Nearly parallel $\G_2$- and
$\G_2^*$-structures are also characterised by the existence of a
Killing spinor, refer \cite{FKMS} respectively \cite{Ka}. 

By Proposition \ref{7to6_prop}, a $\G_2^{(*)}$-structure on a
seven-manifold $(N,\vf)$ induces an $H^{\ve,\tau}$-structure
$(\o,\rho)$ on an oriented hypersurface in $(N,\vf)$. If the
$\G_2^{(*)}$-structure is nearly parallel, the
$H^{\ve,\tau}$-structure satisfies the equation $d \rho = - \ve \mu
\hat \omega$ due to the formulas (\ref{varphi}) and
(\ref{Hodgevf}). This observation motivates the following definition.
\begin{definition}
  An $H^{\ve,\tau}$-structure $(\omega,\rho)$ on a six-manifold $M$ is
  called \emph{nearly half-flat} if
  \begin{eqnarray}
    \label{nhf}
    d \rho = \frac{\lambda}{2} \omega^2 = \lambda \sigma
  \end{eqnarray}
  for some constant $\lambda \in \mathbb R^*$.
\end{definition}
The notion of a nearly half-flat $\SU(3)$-structure was introduced in
\cite{FIMU}, where also evolution equations on six-manifolds leading
to nearly parallel $\G_2$-structures are considered. For compact
manifolds $M$, it is shown in \cite{St} that a solution which is a 
nearly half-flat $\SU(3)$-structure for a time $t=t_0$ already defines 
a nearly parallel $\G_2$-structure. In the following, we extend these
evolution equations to all possible signatures and give a simplified
proof for the properties of the solutions which also holds for
non-compact manifolds.
\bs \label{nearlyprop} Let $H^{\ve,\tau}$ be a real form of
$\SL(3,\ccc)$, $G^{\ve, \tau}$ the corresponding real form of
$\G_2^\ccc$ and $(\rho,\w)$ a one-parameter family of
$H^{\ve,\tau}$-structures on a six-manifold $M$ with a parameter $t$
from an interval $I$. Then, the three-form
\[ \vf = \w\wedge dt +\rho \] defines a nearly parallel
$G^{\ve,\tau}$-structure for the constant $\mu \ne 0$ on $M\times I$
if and only if the $H^{\ve,\tau}$-structure $(\rho,\w)$ is nearly
half-flat for the constant $-\ve \mu$ for all $t \in I$ and satisfies
the evolution equation
\begin{eqnarray}
  \dot{\rho}&=&d\w - \ve \mu \hat \rho \label{nev1}. 
\end{eqnarray}
\es \bprf The assertion follows directly from the following
computation, analogously to the proof of Proposition
\ref{parallelprop}:
\begin{eqnarray*}
  \check{d}\vf
  &=& d\w \wedge dt +dt \wedge\dot\rho + d\rho\ = \ (d\w - \dot\rho )\wedge dt +d\rho,\\
  \mu \ast \vf &=& \ve \mu \left(\hat\rho\wedge dt-\s\right).
\end{eqnarray*}
\eprf

The main theorem for the parallel case generalises as follows. Recall
\eqref{hatsigma} that for a stable four-form $\sigma = \frac{1}{2}
\omega^2 = \hat \omega$, the application of the operator $\sigma
\mapsto \hat \sigma$ yields the stable two-form \[ \hat {\hat \omega}
= \hat \sigma = \frac{1}{2} \omega.\]

\btheo\label{nearlyparalleltheo} Let $(\rho_{0}, \w_{0})$ be a nearly
half-flat $H^{\ve,\tau}$-structure for the constant $\lambda \ne 0$ on
a six-manifold $M$. Let $M$ be oriented such that
$\omega_0^3>0$. Furthermore, let $\rho \in \W^3M$ be a one-parameter
family of stable forms with parameters coming from an interval $I$
such that $\rho(t_0) = \rho_0$ and such that the evolution equation
\begin{eqnarray}
  \dot{\rho}&=& \frac{2}{\lambda} d (\widehat{d \rho}) + \lambda \, \hat \rho \label{nev_rho}
\end{eqnarray}
is satisfied for all $t \in I$.  Then $(\rho, \w = \frac{2}{\lambda}
\widehat{d\rho})$ is a family of nearly half-flat
$H^{\ve,\tau}$-structures for the constant $\lambda$. In particular,
the three-form
\[\vf = \w\wedge dt +\rho\]
defines a nearly parallel $G^{\ve,\tau}$-structure for the constant $-
\ve \lambda$ on $M\times I$.  \etheo

\begin{proof}
  First of all, we observe that $d\rho$ is stable in a neighbourhood
  of the stable form $d\rho_0 = \lambda \sigma_0$, since stability is
  an open condition. Furthermore, the operator $d \rho \mapsto
  \widehat{d\rho}$ is uniquely defined by the orientation induced from
  $\omega_0$. Therefore, the evolution equation is locally
  well-defined and we assume that $\rho$ is a solution on an interval
  $I$. The only possible candidate for a nearly half-flat structure
  for the constant $\lambda$ is $(\rho, \w = \frac{2}{\lambda}
  \widehat{d\rho})$ since only this two-form $\omega$ satisfies the
  nearly half-flat equation $\sigma= \hat \omega = \frac{1}{\lambda}
  d\rho$. Obviously, it holds
  \begin{equation}
    \label{formel3000}
    d\sigma = 0 = d\omega \wedge \omega.
  \end{equation}
  By Proposition \ref{nearlyprop}, it only remains to show that this
  pair of stable forms defines an $H^{\ve,\tau}$-structure, or
  equivalently, that the compatibility conditions \eqref{comp1} and
  \eqref{comp2} are preserved in time. By taking the exterior
  derivative of the evolution equation, we find
  \begin{equation}
    \label{nev2}
    \dot \sigma = \frac{1}{\lambda} d \dot \rho = d \hat \rho   
  \end{equation}
  which is in fact the second evolution equation of the parallel
  case. Completely analogous to the parallel case, the following
  computation implies the first compatibility condition: \be \ddt
  (\s_X \wedge \rho)&=&\dot{\s}_X\wedge \rho + \s_X\wedge \dot{\rho}
  \\
  &\stackrel{(\ref{nev_rho} ),(\ref{nev2} )}{=}& (d\hat{\rho})_X
  \wedge \rho + \s_X \wedge d\w + \lambda \, \sigma_X \wedge \hat\rho
  \\
  &\stackrel{(\ref{formel2}),(\ref{nhf})}{=}& \hat \rho \wedge (d
  \rho)_X
  +\w_X\wedge\w\wedge d\w + (d\rho)_X \wedge \hat\rho \\
  &\stackrel{\eqref{formel3000}}{=}& 0.  \ee

  The proof of the second compatibility condition in \cite{H1} again
  holds literally since the term $\hat \rho \wedge \dot \rho = \hat
  \rho \wedge d \omega$ is the same as in the case of the parallel
  evolution.
\end{proof}

The system (\ref{nev_rho}) of second order in $\rho$ can easily be
reformulated into a system of first order in $(\omega,\rho)$ to which
we can apply the Cauchy-Kovalevskaya theorem. Indeed, a solution
$(\o,\rho)$ of the system
\begin{eqnarray}
  \label{sys}
  \dot{\rho}=d\w + \lambda \hat \rho \;, \qquad \dot \sigma = d \hat \rho ,
\end{eqnarray}
with nearly half-flat initial value $(\o(t_0),\rho(t_0))$ is nearly
half-flat for all $t$ and also satisfies the system
\eqref{nev_rho}. Conversely, (\ref{nev_rho}) implies (\ref{sys}) with
$\s=\hat \o=\frac{1}{\lambda}d \rho$.

Therefore, for an initial nearly half-flat structure which satisfies
assumptions analogous to those of Corollary \ref{Cauchy}, we obtain
existence, uniqueness and naturality of a solution of the system
\eqref{sys}, or, equivalently, of (\ref{nev_rho}).

\subsection{Cocalibrated $\G_2^{(*)}$-structures and parallel $\Spin(7)$- and $\Spin_0(3,4)$-structures}
In \cite{H1}, another evolution equation is introduced which relates
cocalibrated $\G_2$-structures on compact seven-manifolds $M$ to
parallel $\Spin(7)$-structures. As before, we generalise the evolution
equation to non-compact manifolds and indefinite metrics.

As we have already seen in Section \ref{7vs8}, the stabiliser in
$\GL(V)$ of a four-form $\Phi_0$ on an eight-dimensional vector space
$V$ is $\Spin(7)$ or $\Spin_0(3,4)$ if and only if it can be written
as in (\ref{spin4form}) for a stable three-form $\vf$ on a
seven-dimensional subspace with stabiliser $\G_2$- or $\G_2^*$,
respectively. Thus, a $\Spin(7)$- or $\Spin_0(3,4)$-structure on an
eight-manifold $M$ is defined by a four-form $\Phi \in \Omega^4 M$
such that $\Phi_p \in \L^4T^*_pM$ has this property for all $p$. By
formula \eqref{gPhi} for the metric $g_\Phi$ induced by $\Phi$, an
oriented hypersurface in $(M,\Phi)$ with spacelike unit normal vector
field $n$ with respect to $g_\Phi$ carries a natural $\G_2$- or
$\G_2^*$-structure, respectively, defined by $\vf = n \, \inter \Phi$.

A $\Spin(7)$- or $\Spin_0(3,4)$-structure $\Phi$ is \emph{parallel}
if and only if $d \Phi = 0$. We remark that the proof for the
Riemannian case given in \cite[Lemma 12.4]{Sa} is not hard to transfer
to the indefinite case when considering \cite[Proposition 2.5]{B1} and
using  the complexification of the two spin groups. 

Due to this fact, the induced $\G_2^{(*)}$-structure $\vf$ on an
oriented hypersurface in an eight-manifold $M$ with parallel
$\Spin(7)$- or $\Spin_0(3,4)$-structure $\Phi$ is \emph{cocalibrated},
i.e.\ it satisfies
\begin{eqnarray}
  d *_\varphi \varphi = 0.
\end{eqnarray}
Conversely, a cocalibrated $\G_2^{(*)}$-structure can be embedded in
an eight-manifold with parallel $\Spin(7)$- or
$\Spin_0(3,4)$-structure as follows.
\begin{theorem}\label{cocaThm}
  Let $M$ be a seven-manifold and $\varphi \in \Omega^3 M$ be a
  one-parameter family of stable three-forms with a parameter $t$ in
  an interval $I$ satisfying the evolution equation
  \begin{eqnarray}
    \label{spinvolution}
    \frac{\partial}{\partial t}( *_\varphi \varphi) &=& d \varphi.
  \end{eqnarray}
  If $\varphi$ is cocalibrated at $t=t_0 \in I$, then $\varphi$
  defines a family of cocalibrated $\G_2$- or $\G_2^{(*)}$-structures
  for all $t\in I$. Moreover, the four-form
  \begin{equation}
    \label{spinform}
    \Phi = dt \wedge \varphi + *_\varphi \varphi
  \end{equation}
  defines a parallel $\Spin(7)$- or $\Spin_0(3,4)$-structure
  on $M \times I$, respectively, which induces the metric 
  \begin{equation}
    \label{gtPhi}
    g_\Phi = g_\vf + dt^2.
  \end{equation}
\end{theorem}
\begin{proof}
  Since the time derivative of $d * \varphi$ vanishes when inserting
  the evolution equation, the family stays cocalibrated if it is
  cocalibrated at an initial value. As before, we denote by $\check d$
  the exterior differential on $\check M := M \times I$ and
  differentiate the four-form \eqref{spinform}:
  \[ \check d \Phi = - dt \wedge d \varphi + d( * \varphi) + dt \wedge
  \frac{\partial}{\partial t}(* \varphi).\] Obviously, this four-form
  is closed if and only the evolution equation is satisfied and the
  family is cocalibrated. The formula for the induced metric
  corresponds to formula (\ref{gPhi}).
\end{proof}
As before, the Cauchy-Kovalevskaya theorem guarantees existence and
uniqueness of solutions if assumptions analogous to those of Corollary
\ref{Cauchy} are satisfied.
\begin{remark}
  We observe that nearly parallel $\G_2$- and $\G_2^*$-structures are
  in particular cocalibrated such that analytic nearly half-flat
  structures in dimension six can be embedded in parallel $\Spin(7)$-
  or $\Spin_0(3,4)$-structures in dimension eight by evolving them
  twice with the help of the Theorems \ref{nearlyparalleltheo} and
  \ref{cocaThm}.
\end{remark}

\section{Evolution of nearly $\ve$-K\"ahler manifolds}
\label{NKexamples}
In this section, we consider the evolution of nearly pseudo-K\"ahler
and nearly para-K\"ahler six-manifolds which can be unified by the
notion of a nearly $\ve$-K\"ahler manifold. The explicit solution of
the Hitchin flow yields a simple and unified proof for the
correspondence of nearly $\ve$-K\"ahler manifolds and parallel
$\G_2^{(*)}$-structures on cones. We complete the picture by
considering similarly the evolution of nearly K\"ahler structures to
nearly parallel $\G_2^{(*)}$-structures on (hyperbolic) sine cones and
the evolution of nearly parallel $\G_2^{(*)}$-structures to parallel
$\Spin(7)$- and $\Spin_0(3,4)$-structures on cones. Our presentation
in terms of differential forms unifies various results in the
literature, which were originally obtained using spinorial methods,
and applies to all possible real forms of the relevant groups.

\subsection{Cones over nearly $\ve$-K\"ahler manifolds}
In the language of \cite{AC2} and \cite{SSH}, an \emph{almost
  $\ve$-Hermitian manifold} $(M^{2m},g,J)$ is defined by an almost
$\ve$-complex structure $J$ which squares to $\ve \mathrm{id}$ and a
pseudo-Riemannian metric $g$ which is $\ve$-Hermitian in the sense
that $g(J\cdot,J\cdot) =-\ve g(\cdot,\cdot)$. Consequently, a
\emph{nearly $\ve$-K\"ahler manifold} is defined as an almost
$\ve$-Hermitian manifold such that $\nabla J$ is skew-symmetric.  On a
six-manifold $M$, a nearly $\ve$-K\"ahler structure $(g,J,\o)$ with $|
\n J |^2 = 4$ (i.e.\ of constant type $1$ in the terminology of
\cite{gray76}) is equivalent to a normalised $H^{\ve,\tau}$-structure
$(\omega,\rho)$ which satisfies 
\begin{eqnarray}
  d \omega &=& 3 \rho \,, \label{NK1}\\  
  d \hat \rho &=& 4 \hat \omega. \label{NK2}
\end{eqnarray}
This result is well-known for Riemannian signature \cite{RC} and is
generalised to arbitrary signature in \cite[Theorem 3.14]{SSH}. In
particular, nearly $\ve$-K\"ahler structures $(\omega,\rho)$ in
dimension six are half-flat and the structure $(\o,\hat \rho)$ is
nearly half-flat (for the constant $\lambda=4$).

\begin{Prop}
  \label{prop_NKvscone}
  Let $(M,h_0)$ be a pseudo-Riemannian six-manifold of signature
  $(6,0)$, $(4,2)$ or $(3,3)$ and let $(\bar M= M \times \bR^+,\bar
  g_\ve = h_0 - \ve dt^2)$ be the timelike cone for $\ve=1$ and the
  spacelike cone for $\ve=-1$. There is a one-to-one correspondence
  between nearly $\ve$-K\"ahler structures $(h_0,J)$ with $|\n J|^2
  =4$ on $(M,h_0)$ and parallel $\G_2$- and
  $\G_2^*$-structures $\vf$ on $\bar M$ which induce the cone metric
  $\bar g_\ve$.
\end{Prop}
\begin{proof}
  This well-known fact is usually proved using Killing spinors, see
  \cite{B}, \cite{Gru} and \cite{Ka2}. We give a proof relying
  exclusively on the framework of stable forms and the Hitchin
  flow. For Riemannian signature, this point of view is also
  adopted in \cite{ChSa} and \cite{Bu}.

  The $H^{\ve,\tau}$-structures inducing the given metric $h_0$ are
  the reductions of the bundle of orthonormal frames of $(M,h_0)$
  to the respective group $H^{\ve,\tau}$. Given any
  $H^{\ve,\tau}$-reduction $(\omega_0, \rho_0)$ of $h_0$, we consider
  for $t\in \bR^+$ the one-parameter family
  \begin{equation}
    \label{NKfamily}
    \omega = t^2 \omega_0\, , \: \rho = t^3 \rho_0,   
  \end{equation}
  which induces the family of metrics $h=t^2 h_0$.  By formula
  \eqref{gvarphidt}, the metric $g_\varphi$ on $\bar M$ induced by the
  stable three-form $\varphi = \omega \wedge dt + \rho$ is exactly the
  cone metric $\bar g_\ve$.  

  It is easily verified that the family \eqref{NKfamily} consists of
  half-flat structures satisfying the evolution equations if and only
  if the initial value $(\omega(1),\rho(1))=(\omega_0, \rho_0)$
  satisfies the exterior system \eqref{NK1}, \eqref{NK2}. Therefore,
  the stable three-form $\vf$ on the cone $(\bar M,\bar g_\ve)$ is
  parallel if and only if the $H^{\ve,\tau}$-reduction $(\omega_0,
  \rho_0)$ of $h_0$ is a nearly $\ve$-K\"ahler structure with $| \n J
  |^2 =4$.

  Conversely, let $\vf$ be a stable three-form on $\bar M$ which
  induces the cone metric $\bar g_\ve$. Since $\partial_t$ is a normal
  vector field for the hypersurface $M = M \times \{ 1 \}$ satisfying
  $\bar g(\partial_t,\partial_t) = -\ve$, we obtain an
  $H^{\ve,\tau}$-reduction $(\o_0,\rho_0)$ of $h_0$ defined by 
  \begin{equation}
    \label{inducedNK}
    \o_0 = \partial_t \, \inter\, \vf \, , \qquad \rho_0 = \vf_{|TM}
  \end{equation}
  with the help of Proposition \ref{7to6_prop}. Since the two
  constructions are inverse to each other, the proposition follows.
\end{proof}
\begin{example}
  \label{pseudospheres}
  Consider the flat $(\bR^{(3,4)} \setminus \{ 0 \},\la .,.\ra)$ which
  is isometric to the cone $(M^\ve \times \R^+, t^2 h_\ve - \ve dt^2)$
  over the pseudo-spheres $M^\ve := \{ p \in \bR^{(3,4)} \, | \,
  \langle p,p \rangle = -\ve \}$, $\ve=\pm 1$, with the standard
  metrics $h_\ve$ of constant sectional curvature $-\ve$ and signature
  $(2,4)$ for $\ve=-1$ and $(3,3)$ for $\ve=1$. Obviously, a stable
  three-form $\vf$ inducing the flat metric $\la .,.\ra$ is parallel
  if and only if it is constant. Thus, the previous discussion and
  Proposition \ref{7to6_prop}, in particular formula \eqref{Jcross},
  yield a bijection
  \begin{eqnarray*}
    \SO(3,4)/\,\G_2^* &\rightarrow & \{ \mbox{$\ve$-complex structures $J$ on $M^\ve$ such that $(h_\ve,J)$ is nearly $\ve$-K\"ahler} \} \\
    \vf &\mapsto& J \qquad \mbox{with} \quad J_p(v) = - p \times v,\qquad \forall \, p \in M^\ve
  \end{eqnarray*}
  where the cross-product $\times$ induced by $\vf$ is defined by
  formula \eqref{stablevsO}. In other words, the pseudo-spheres
  ($M^\ve,h_\ve$) admit a nearly $\ve$-K\"ahler structure which is
  unique up to conjugation by the isometry group $\Og(3,4)$ of
  $h_\ve$. In fact, these $\ve$-complex structures on the
  pseudo-spheres are already considered in \cite{L} and the nearly
  para-K\"ahler property for $\ve=1$ is for instance shown in
  \cite{Be}.
\end{example}

\subsection{Sine cones over nearly $\ve$-K\"ahler manifolds} 
For Riemannian signature, it has been shown in \cite{FIMU} that the
evolution of a nearly K\"ahler $\SU(3)$-structure to a nearly parallel
$\G_2$-structure induces the Einstein sine cone metric. This result
can be extended as follows. We prefer to consider (hyperbolic) cosine
cones since they are defined on all of $\bR$ in the hyperbolic case.
\begin{Prop}
  \label{prop_cosinecone}
  Let $(M,h_0)$ be a pseudo-Riemannian six-manifold.
  \begin{enumerate}[(i)]
  \item If $h_0$ is Riemannian, or has signature $(2,4)$,
    respectively, there is a one-to-one correspondence between nearly
    (pseudo-)K\"ahler structures $(h_0,J)$ on $M$ with $| \n J |^2 =4$
    and nearly parallel $\G_2$-structures, or $\G_2^*$-structures,
    respectively, for the constant $\mu=-4$ on the spacelike cosine
    cone
    \[(M \times (-\frac{\pi}{2},\frac{\pi}{2}), \, \cos^2(t) h_0 +
    dt^2).\]
  \item If $h_0$ has signature $(3,3)$, there is a one-to-one
    correspondence between nearly para-K\"ahler structures $(h_0,J)$
    on $M$ with $| \n J |^2 =4$ and nearly parallel $\G_2^*$-structures
    for the constant $\mu=4$ on the timelike hyperbolic cosine cone
    \[ (M \times \bR, \, - \cosh^2(t) h_0 - dt^2).\]
  \end{enumerate}
\end{Prop}
\begin{proof}
  \begin{enumerate}[(i)]
  \item Starting with any $\SU(3)$- or $\SU(1,2)$-reduction
    $(\omega_0,\rho_0)$ of $h_0$, the one-parameter family
    \[ \omega = \cos^2(t) \omega_0 \, , \: \rho = - \cos^3(t) (\sin(t)
    \rho_0 + \cos(t) \hat \rho_0) \] with $(\omega(0),\rho(0)) =
    (\omega_0,-\hat \rho_0)$ defines a stable three-form $\vf = \o
    \wedge dt + \rho$ on $M \times
    (-\frac{\pi}{2},\frac{\pi}{2})$. Since $z \Psi_0= z(\rho_0 +\i
    \hat \rho_0)$ is a $(3,0)$-form w.r.t. the induced almost complex
    structures $J_{\mathrm{Re}\,(z \Psi_0)}$ for all $z \in \bC^*$, the
    structure $J_\rho=J_{\rho_0}$ is constant in $t$. Thus, the metric
    $g_\vf$ induced by $\vf$ is the cosine cone metric. Moreover, it
    holds $\hat \rho = - \cos^3(t) (\sin(t) \hat \rho_0 - \cos(t)
    \rho_0)$ due to Corollary \ref{hat}.

    It takes a short calculation to verify that the one-parameter
    family is nearly half-flat (for the constant $\lambda = - 4$) and
    satisfies the evolution equation \eqref{nev1} if and only
    $(\omega_0,\rho_0)$ satisfies the exterior system \eqref{NK1},
    \eqref{NK2}. Thus, applying Proposition \ref{nearlyprop}, the
    three-form $\vf = \w\wedge dt +\rho$ defines a nearly parallel
    $G^{\ve,\tau}$-structure on $M \times
    (-\frac{\pi}{2},\frac{\pi}{2})$ (for the constant $\mu =-4$) if
    and only if $(h_0,J_{\rho_0})$ is nearly $\ve$-K\"ahler with $| \n
    J |^2 = 4$.

    The inverse construction is given by \eqref{inducedNK} in analogy
    to the case of the ordinary cone.
  \item The proof in the para-complex case is completely analogous if
    we consider the one-parameter family
    \[ \omega = \cosh^2(t) \omega_0 \, , \: \rho = - \cosh^3(t)
    (\sinh(t) \rho_0 + \cosh(t) \hat \rho_0) \] which is defined for
    all $t \in \bR$. We note the following subtleties regarding
    signs. By Proposition \ref{invmap}, we know that the mapping $\rho
    \mapsto \hat \rho$ is homogeneous of degree $1$, but not
    linear. Indeed, by applying Corollary \ref{hat}, we find
    \[\widehat{\sinh(t) \rho_0 + \cosh(t) \hat \rho_0} = - \sinh(t)
    \hat \rho_0 - \cosh(t) \rho_0.\] Using this formula, one can check
    that $J_\rho = J_{\hat \rho_0}=-J_{\rho_0}$ is constant in $t$ such
    that the metric induced by $(\o,\rho)$ is in fact $h=-\cosh^2(t)h_0$.
  \end{enumerate}
\end{proof}
The fact that the (hyperbolic) cosine cone over a six-manifold
carrying a Killing spinor carries again a Killing spinor was proven in
\cite{Ka}. By relating spinors to differential forms, these results
also imply the existence of a nearly parallel $\G_2^{(*)}$-structures
on the (hyperbolic) cosine cone over a nearly $\ve$-K\"ahler manifold.
\begin{example}
  \label{ex_cosinecone}
  The (hyperbolic) cosine cone of the pseudo-spheres $(M^\ve,h_\ve)$
  of Example \ref{pseudospheres} has constant sectional curvature $1$,
  for instance due to \cite[Corollary 2.3]{ACGL}, and is thus
  (locally) isometric to the pseudo-sphere $S^{3,4} = \{ p \in
  \bR^{(4,4)} \, | \, \langle p,p \rangle = 1 \} = \Spin_0(3,4)/\,\G_2^*$.
\end{example}

\subsection{Cones over nearly parallel $\G_2^{(*)}$-structures}
By Lemma 9 in \cite{B}, there is a one-to-one correspondence on a
Riemannian seven-manifold $(M,g_0)$ between nearly parallel
$\G_2$-structures
and parallel $\Spin(7)$-structures on the Riemannian cone. In order to
illustrate the evolution equations for nearly parallel
$G_2^*$-structures, we extend this result to the indefinite case by
applying Theorem \ref{cocaThm}. This is possible since nearly parallel
$\G_2^*$-structures are in particular cocalibrated. Again, the fact
that the cone over a nearly parallel $\G_2^*$-manifold admits a
parallel spinor can be derived from the connection to Killing spinors as 
observed in \cite{Ka}.
\begin{Prop}
  \label{prop_npvsparallelspin}
  Let $(M,g_0)$ be a pseudo-Riemannian seven-manifold of signature
  $(3,4)$. There is a one-to-one correspondence between nearly
  parallel $\G_2^*$-structures for the constant $4$ which induce the
  given metric $g_0$ and parallel $\Spin_0(3,4)$-structures on $M
  \times \bR^+$ inducing the cone metric $\bar g = t^2 g_0 + dt^2$.
\end{Prop}
\begin{proof}
  Let $\vf_0$ be any cocalibrated $G_2^*$-structure on $M$ inducing
  the metric $g_0$. The one-parameter family of three-forms defined by
  $\vf = t^3 \vf_0$ for $t \in \bR^+$ induces the family of metrics
  $g=t^2 g_0$ such that the Hodge duals are $*_\vf \vf = t^4 *_{\vf_0}
  \vf_0$. By \eqref{gtPhi}, the $\Spin_0(3,4)$-structure $\Psi = dt
  \wedge \vf + *_\vf \vf$ on $M \times \bR^+$ induces the cone metric
  $\bar g$. Conversely, given a $\Spin_0(3,4)$-structure $\Psi$ on the
  cone $(M \times \bR^+,\bar g)$, we have the cocalibrated
  $G_2^*$-structure $\vf_0 = \partial_t \, \inter \Psi$ on $M$, which
  also induces the given metric $g_0$. Since the evolution equation
  (\ref{spinvolution}) is satisfied if and only if the initial value
  $\vf_0$ is nearly parallel for the constant $4$ and since the two
  constructions are inverse to each other, the assertion follows from
  Theorem \ref{cocaThm}.
\end{proof}
\begin{example}
  We consider again the easiest example, i.e.\ the flat $\bR^{(4,4)}
  \setminus \{0\}$ which is isometric to the cone over the
  pseudo-sphere $S^{3,4}$. Analogous to Example \ref{pseudospheres},
  the proposition just proved yields a proof of the fact that the
  nearly parallel $G_2^*$-structures for the constant $4$ on $S^{3,4}$
  are parametrised by $\SO(4,4)/\,\Spin_0(3,4)$, i.e.\ by the four
  homogeneous spaces (\ref{vierhomogene}). In particular, these
  structures are conjugated by the isometry group $\Og(4,4)$ of
  $S^{3,4}$. 

  Summarising the application of the three Propositions
  \ref{prop_NKvscone}, \ref{prop_cosinecone} and
  \ref{prop_npvsparallelspin} to pseudo-spheres, we find a mutual
  one-to-one correspondence between 
  \begin{enumerate}[(1)]
    \item nearly pseudo-K\"ahler structures with $|\n J|^2 \ne 0$ on $(S^{2,4},g_{can})$,
    \item nearly para-K\"ahler structures with $|\n J|^2 \ne 0$ on $(S^{3,3},g_{can})$,
    \item parallel $\G_2^*$-structures on $(\bR^{(3,4)},g_{can})$,
    \item nearly parallel $\G_2^*$-structures on the spacelike cosine
      cone over $(S^{2,4},g_{can})$, 
    \item nearly parallel $\G_2^*$-structures on the
      timelike hyperbolic cosine cone over $(S^{3,3},g_{can})$, 
    \item nearly parallel $\G_2^*$-structures on $(S^{3,4},g_{can})$ and 
    \item parallel $\Spin_0(3,4)$-structures on $(\bR^{(4,4)},g_{can})$.
  \end{enumerate}
  This geometric correspondence is reflected in the
  algebraic fact that the four homogeneous spaces (\ref{vierhomogene})
  are isomorphic.
\end{example}

\section{The evolution equations on nilmanifolds $\Gamma \, \backslash
  \, H_3 \!\times\! H_3 $}
\label{nilexamples}
Let $H_3$ be the three-dimensional real Heisenberg group with Lie
algebra $\h_3$.  In this section, we will develop a method to
explicitly determine the parallel $\G_2^{(*)}$-structure induced by an
arbitrary invariant half-flat structure on a nilmanifold $\Gamma \,
\backslash \, H_3 \!\times\! H_3 $ without integrating. In particular,
this method is applied to construct three explicit large families of
metrics with holonomy equal to $\G_2$ or $\G_2^*$, respectively.

\subsection{Evolution of invariant half-flat structures on
  nilmanifolds}
\label{nivolution}
Left-invariant half-flat structures $(\omega_0,\rho_0)$ on a Lie group
$G$ are in one-to-one correspondence with normalised pairs $(\o,\rho)$
of compatible stable forms on the Lie algebra $\mathfrak{g}$ of $G$
which satisfy $d\rho=0$ and $d\o^2=0$. To shorten the notation, we
will speak of a \emph{half-flat structure on a Lie algebra}.

Given as initial value a half-flat structure on a Lie algebra, the
evolution equations
\begin{eqnarray}
  \label{evnochmal}
  \dot \rho = d \omega \; , \qquad \dot \sigma = d \hat \rho \; ,
\end{eqnarray}
reduce to a system of ordinary differential equations and a unique
solution exists on a maximal interval $I$. Due to the structure of
the equation, the solution differs from the initial values by adding
exact forms to $\sigma_0$ and $\rho_0$. In other words, an initial
value $(\sigma_0,\rho_0)$ evolves within the product $[\sigma_0]
\times [\rho_0]$ of their respective Lie algebra cohomology classes.

Every nilpotent Lie group $N$ with rational structure constants admits
a cocompact lattice $\Gamma$ and the resulting compact quotients
$\Gamma \backslash N$ are called nilmanifolds. Recall that a geometric
structure on a nilmanifold $\Gamma \backslash N$ is called
\emph{invariant} if is induced by a left-invariant geometric structure
on $N$.

Explicit solutions of the Hitchin flow equations on several nilpotent
Lie algebras can be found for instance in \cite{CF} and \cite{AS}. In
both cases, a metric with holonomy contained in $\G_2$ has been
constructed before by a different method and this information is used
to obtain the solution. For a symplectic half-flat initial value,
another explicit solution on one of these Lie algebras is given in
\cite{CT}. In all cases, the solution depends only on one variable.

At least for four nilpotent Lie algebras including $\h_3 \op \h_3 $, a
reason for the simple structure of the solutions has been observed in
\cite{AS}. Indeed, the following lemma shows that the evolution of
$\sigma$ takes place in a one-dimensional space. As usual, we define a
nilpotent Lie algebra by giving the image of a basis of one-forms
under the exterior derivative, see for instance \cite{Sa}. The same
reference also contains a list of all six-dimensional nilpotent Lie
algebras.
\begin{lemma}
  \label{trick17}
  Let $\rho$ be a closed stable three-form with dual three-form $\hat
  \rho$ on a six-dimensional nilpotent Lie algebra $\g$.
  \begin{enumerate}[(i)]
  \item If $\g$ is one of the three Lie algebras
    \[ (0,0,0,0,e^{12},e^{34})\; , \quad
    (0,0,0,0,e^{13}+e^{42},e^{14}+e^{23})\, , \quad
    (0,0,0,0,e^{12},e^{14}+e^{23}),\] then $d\hat \rho \in \L^4 U$
    for the four-dimensional kernel $U$ of $d : \Lambda^1 \g^*
    \rightarrow \Lambda^2 \g^*$.
  \item If $\g$ is the Lie algebra
    \[ (0,0,0,0,0,e^{12}+e^{34}),\] then $d\hat \rho \in \L^4 U$
    for the four-dimensional subspace $U = {\rm span}\{ e^1, e^2, e^3,
    e^4 \}$ of $\ker d$.
  \end{enumerate}
\end{lemma}
\begin{remark}
  The assertion of the lemma is not true for the remaining
  six-dimensional nilpotent Lie algebras with $b_1=\mathrm{dim}
  (\mathrm{ker}\,d) = 4$ or $b_1=5$. In each case, we have constructed
  a closed stable $\rho$ such that $d \hat \rho$ is \textrm{not}
  contained in $\L^4 (\mathrm{ker}\,d)$.
\end{remark}
In fact, this lemma can also be viewed as a corollary of the following
lemma which we will prove first.
\begin{lemma}
  \label{Jinvar}
  Let $\rho$ be a closed stable three-form on one of the four Lie
  algebras of Lemma \ref{trick17} and let $U$ be the four-dimensional
  subspace of $\ker d$ defined there. In all four cases, the space $U$
  is $J_\rho$-invariant where $J_\rho$ denotes the almost
  (para-)complex structure induced by $\rho$.
\end{lemma}
\begin{proof}
  For $\lambda (\rho) < 0$, the assertion is similar to that of
  \cite[Lemma 2]{AS}. However, since the only proof seems to be given
  for the Iwasawa algebra for integrable $J$ in \cite[Theorem
  1.1]{KeS}, we give a complete proof.

  Let $\g$ be one of the three Lie algebras given in part (i) of Lemma
  \ref{trick17} and $U=\ker d$. Obviously, the two-dimensional image
  of $d$ lies within $\Lambda^2 U$ in all three cases. By $J=J_\rho$
  we denote the almost (para-)complex structure associated to the
  closed stable three-form $\rho$.  As before, we denote by $\ve \in
  \{ \pm 1 \}$ the sign of $\lambda(\rho)$ such that $J_\rho = \ve
  \rm{id}$. Let the symbol $\i_\ve$ be defined by the property
  $\i_\ve^2=\ve$ such that the para-complex numbers and the complex
  numbers can be unified by $\bC_\ve = \bR[\rm i_\ve]$. Thus, a
  $(1,0)$-form can be defined for both values of $\ve$ as an eigenform
  of $J_\rho$ in $\Lambda^1 \g^* \otimes \bC_\ve$ for the eigenvalue
  $\i_\ve$.

  We define the $J$-invariant subspace $W:= U \cap J^* U$ of $\g$ such
  that $2 \le \dim W \le 4$. In fact, $\dim W = 4$ is equivalent to
  the assertion.  The other two cases are not possible, which can be
  seen as follows. To begin with, assume that $W$ is
  two-dimensional. When choosing a complement $W'$ of $W$ in $U$, we
  have by definition of $W$ that
  \[V = W \oplus W' \oplus J^* W'.\] We observe that, for $\ve=1$, the
  $\pm 1$-eigenspaces of $J$ restricted to $W' \oplus J^* W'$ are both
  two-dimensional. Therefore, we can choose for both values of $\ve$ a
  basis $\{ e^1 , e^2, e^3, e^4=J^*e^1, e^5=J^*e^2, e^6=J^*e^3 \}$ of
  $V$ such that $e^1, e^2, e^3$ and $e^4$ are closed and $de^5, de^6
  \in \Lambda^2 U $.  Since $\rho + \i_\ve J^*_\rho \rho$ is a
  (3,0)-form in both cases, it is possible to change the basis vectors
  $e^1, e^4$ within $W \subset \ker d$ such that
  \[ \rho + \i_\ve J^*_\rho \rho = (e^1 + \i_\ve e^4)\^(e^2 + \i_\ve
  e^5)\^(e^3 + \i_\ve e^6) \]
  and thus
  \[ \rho = e^{123} + \ve e^{156} - \ve e^{246} + \ve e^{345}. \] By
  construction of the basis, we have that \[ 0 = d \rho = - \ve e^1
  \wedge d e^5 \wedge e^6 + \ve e^1 \wedge e^5 \wedge d e^6 +
  \alpha \] with $\alpha \in \Lambda^4 U$. As the first two summands
  are linearly independent and not in $\Lambda^4 U$, we conclude that
  both $e^1 \wedge de^5$ and $e^1 \wedge de^6$ vanish. Thus, the
  closed one-form $e^1$ has the property that the wedge product of
  $e^1$ with any exact two-form vanishes. However, an inspection of
  the standard basis of each of the three Lie algebras in question
  reveals that such a one-form does not exist on these Lie algebras
  and we have a contradiction to $\dim W=2$.

  Since a $J$-invariant space cannot be three-dimensional for
  $\ve=-1$, the proof is finished for this case. However, if $\ve=1$,
  the case $\dim W =3$ cannot be excluded that easy. Assuming that it
  is in fact $\dim W =3$, we choose again a complement $W'$ of $W$ in
  $U$ and find a decomposition
  \[V = W \oplus W' \oplus J^* W' \oplus W'' \] with $J^*W'' =
  W''$. Without restricting generality, we can assume that $J$ acts
  trivially on $W''$. Then, we find a basis for $V$ such that the
  $+1$-eigenspace of $J$ is spanned by $\{ e^1, e^4 + e^5, e^6 \}$ and
  the $-1$-eigenspace by $\{ e^2,e^3, e^4- e^5 \}$, where $e^1, e^2,
  e^3$ and $e^4$ are closed and $e^5 = J^*e^4$. Since the given closed
  three-form $\rho$ generates this $J$, it has to be of the form
  \[ \rho = a e^1 \^ (e^4 + e^5) \^ e^6 + b e^{23} \wedge (e^4- e^5)\]
  for two real constants $a,b$.  The vanishing exterior derivative
  \[ d \rho = a e^1 \wedge d( e^{56} ) \qquad \mbox{mod} \quad
  \Lambda^4 U \] leads to the same contradiction as in the first case
  and part (i) is shown.

  In fact, the same arguments apply to the Lie algebra of part
  (ii). The four-dimensional space $U \subset \ker d$ spanned by $\{
  e^1, ..., e^4 \}$ also satisfies $\im d \subset \Lambda^2 U$. Going
  through the above arguments, the only difference is that $e^5$ or
  $e^6$ may be closed. However, at least one of them is not closed and
  its image under $d$ generates the exact two-forms. Again, there is
  no one-form $\beta \in U$ such that $\beta \wedge \gamma=0$ for all
  exact two-forms $\gamma$ and the arguments given in part (i) lead to
  contradictions for both $\dim W =2$ and $\dim W =3$.
\end{proof}

\begin{proof}[Proof of Lemma \ref{trick17}]
  Let $\rho$ be a closed stable three-form on one of the four
  nilpotent Lie algebras and $U \subset \ker d$ as defined in the
  lemma. For both values of $\ve$, we can apply Lemma \ref{Jinvar} and
  choose two linearly independent closed $(1,0)$-forms $E^1$ and $E^2$
  within the $J_\rho$-invariant space $U \ot \bC_\ve$. Considering
  that $\rho + \rm i_\ve \hat \rho$ is a $(3,0)$-form for both values
  of $\ve$, there is a third $(1,0)$-form $E^3$ such that $\rho + \rm
  i_\ve \hat \rho = E^{123}$. Since $d\rho=0$ and $\im d \subset
  \Lambda^2 U$, it follows that the exterior derivative
  \[ d \hat \rho = \ve \i_\ve d (E^{123}) = \ve \i_\ve E^{12}
  \wedge dE^3\] is an element of $\Lambda^4 U$.
\end{proof}

\subsection{Left-invariant half-flat structures on $H_3 \times H_3 $}
\label{halfflatonH3H3}
From now on, we focus on the Lie algebra $\g=\h_3 \op \h_3$. Apart
from describing all half-flat structures on this Lie algebra, i.e.\
all initial values for the evolution equations, we give various
explicit examples and prove a strong rigidity result concerning the
induced metric.

Obviously, pairs of compatible stable forms on a Lie algebra which are
isomorphic by a Lie algebra automorphism induce equivalent
$H^{\ve,\tau}$-structures on the corresponding simply connected Lie
group. Thus, we derive, to begin with, a normal form modulo Lie
algebra automorphisms for stable two-forms $\o \in \Lambda^2 \g^*$
which satisfy $d\o^2=0$.

A basis $\{ e_1, e_2, e_3, f_1, f_2, f_3 \}$ for $\h_3 \oplus \h_3$
such that the only non-vanishing Lie brackets are given by
\[de^3 = e^{12}, \qquad df^3 = f^{12},\] will be called a {\it
  standard basis}. The connected component of the automorphism group
of the Lie algebra $\h_3 \oplus \h_3$ in the standard basis is
\begin{equation}
  \label{aut}
  \rm{Aut}_0(\h_3 \oplus h_3) = 
  \left\{
    \begin{pmatrix}
      A & 0 & 0 & 0 \\
      a^t & \det(A) & c^t & 0\\
      0 & 0 & B & 0\\
      d^t & 0 & b^t & \det(B)
    \end{pmatrix} , \: A,B \in \GL(2,\bR), \; a,b,c,d \in \bR^2
  \right\}.
\end{equation}
We denote by $\g_i$, $i=1,2$, the two summands, by $\z_i$ their
centres and by $\z$ the centre of $\g$. The annihilator of the centre
is $\z^0 = \rm{ker}\,d$ and similarly for the summands by restricting
$d$. We have the decompositions
\begin{eqnarray*}
  \g^* &\cong& \z_1^0 \op \z_2^0 \op \frac{ \g_1^*}{\z_1^0} \op \frac{\g_2^*}{\z_2^0}, \\
  \L^2 \g^* &\cong& \Lambda^2(\z^0) \op 
  \underbrace{(\frac{\g_1^*}{\z_1^0} \wedge \frac{\g_2^*}{\z_2^0})}_{\k_1} \op
  \underbrace{(\z_1^0 \wedge \frac{\g_2^*}{\z_2^0})}_{\k_2}  \op 
  \underbrace{(\z_2^0 \wedge \frac{\g_1^*}{\z_1^0})}_{\k_3} \op 
  \underbrace{(\z_1^0 \wedge \frac{\g_1^*}{\z_1^0})  \op 
    (\z_2^0 \wedge \frac{\g_2^*}{\z_2^0}) }_{\k_4}.
\end{eqnarray*}
By $\o^{\k_i}$ we denote the projection of a two-form $\o$ onto one of
the spaces $\k_i$, $i=1,2,3,4$, defined as indicated in the
decomposition. We observe that $\k_1 = \L^2 (\frac{\g^*}{\z^0})$ and
$\o^{\k_1}=0$ if and only if $\o(\z,\z)=0$.
\begin{lemma}
  Consider the action of $\rm{Aut}(\h_3 \oplus \h_3)$ on the set of
  non-degenerate two-forms $\o$ on $\g$ with $d\omega^2 = 0$. The
  orbits modulo rescaling are represented in a standard basis by the
  following two-forms:
  \begin{eqnarray*}
    \label{onormal}
    \omega_1 = & e^1f^1 +  e^2f^2 + e^3f^3, & \mbox{if $\o^{\k_1} \ne 0$,} \\
    \omega_2 = & e^2f^2 + e^{13} + f^{13}, & \mbox{if $d\omega = 0  \iff \o^{\k_1} = 0$, $\o^{\k_2} = 0$, $\o^{\k_3} = 0$,} \\
    \omega_3 = & e^1f^3 + e^2f^2 + e^3f^1, & \mbox{if $\o^{\k_1} = 0$, $\o^{\k_2} \ne 0$, $\o^{\k_3} \ne 0$, $\o^{\k_4} = 0$,}\\
    \omega_4 = & e^1f^3 + e^2f^2 + e^3f^1 + e^{13} + \beta f^{13}, & \mbox{if $\o^{\k_1} = 0$, $\o^{\k_2} \ne 0$, $\o^{\k_3} \ne 0$, $\o^{\k_4} \ne 0$,}\\
    \omega_5 = & e^1f^3 + e^2f^2 + e^{13} + f^{13} & \mbox{otherwise,}
  \end{eqnarray*}
  where $\beta \in \bR$ and $\beta \ne -1$.
\end{lemma}
\begin{proof}
  Let \[\omega = \sum \a_{i} e^{(i+1)(i+2)} + \sum \beta_{i}
  f^{(i+1)(i+2)} + \sum \gamma_{i,j} e^{i}f^{j}\] be an arbitrary
  non-degenerate two-form expressed in a standard basis.  We will give
  in each case explicitly a change of standard basis by an
  automorphism of the form \eqref{aut} with the notation
  \[ A = \begin{pmatrix} a_1 & a_2 \\ a_3 & a_4 \end{pmatrix}, \:
  a^t=(a_5,a_6), \: B = \begin{pmatrix} b_1 & b_2 \\ b_3 &
    b_4 \end{pmatrix}, \: b^t=(b_5,b_6), \: c^t = (c_1,c_2), \: d^t =
  (d_1,d_2). \]

  First of all, if $\o^{\k_1} \ne 0$, the term $\ga_{3,3} e^3f^3$ is
  different from zero and we rescale such that $\ga_{3,3}=1$. Then,
  the application of the change of basis
  \begin{eqnarray*}
    a_1 &=& 1, \:\: a_2 = 0, \:\: a_{3} = 0, \:\: a_4 = 1, \:\: a_{5} = -\gamma_{1, 3}, \:\: a_{6} = -\gamma_{2, 3},  \\ 
    b_{1} &=& \gamma_{2, 2}-\gamma_{2, 3} \gamma_{3, 2}-\alpha_{1} \beta_{1}, \:\:
    b_{2} = -\gamma_{1, 2}-\beta_{1} \alpha_{2}+\gamma_{3, 2} \gamma_{1, 3}, \:\:
    b_{3} = -\gamma_{2, 1}+\gamma_{3, 1} \gamma_{2, 3}-\beta_{2} \alpha_{1},   \\ 
    b_{4} &=& \gamma_{1, 1}-\alpha_{2} \beta_{2}-\gamma_{1, 3} \gamma_{3, 1}, \:\:
    b_{5} = -\gamma_{3, 1} \gamma_{2, 2}+\gamma_{3, 1} \alpha_{1} \beta_{1}+\gamma_{3, 2} \gamma_{2, 1}+\gamma_{3, 2} \beta_{2} \alpha_{1},  \\ 
    b_{6} &=& \gamma_{3, 1} \gamma_{1, 2}+\gamma_{3, 1} \beta_{1} \alpha_{2}-\gamma_{3, 2} \gamma_{1, 1}+\gamma_{3, 2} \alpha_{2} \beta_{2}, \:\: \\
    c_{1} &=& \beta_{2} \gamma_{2, 2}-\beta_{2} \gamma_{2, 3} \gamma_{3, 2}+\beta_{1} \gamma_{2, 1}-\beta_{1} \gamma_{3, 1} \gamma_{2, 3}, \:\: d_{1} = -\alpha_{2}, \\
    c_{2} &=& -\beta_{2} \gamma_{1, 2}+\beta_{2} \gamma_{3, 2} \gamma_{1, 3}-\beta_{1} \gamma_{1, 1}+\beta_{1} \gamma_{1, 3} \gamma_{3, 1},  \:\: d_{2} = \alpha_{1},  
  \end{eqnarray*}
  transforms $\omega$ into $\tilde \omega = \tilde \gamma_{1,1}
  (e^1f^1 + e^2f^2 + e^3f^3) + \tilde \alpha_3 e^{12} + \tilde \beta_3
  f^{12}$, $\tilde \ga_{1,1} \ne 0$. This two-form satisfies $d\tilde
  \omega^2 = 0$ if and only if $\tilde \alpha_3= 0, \tilde \beta_3 =0$
  and the normal form $\o_1$ is achieved by rescaling.

  Secondly, the vanishing of $d\o$ corresponds to $\o^{\k_1} = 0$,
  $\o^{\k_2} = 0$, $\o^{\k_3} = 0$ or $\ga_{3,3}=\ga_{1,3} = \ga_{2,3}
  =\ga_{3,1} = \ga_{3,2} =0$ in a standard basis. By non-degeneracy,
  at least one of $\a_1$ and $\a_2$ is not zero and we can always
  achieve $\a_1=0,\a_2 \ne 0$. Indeed, if $\a_1 \ne 0$, we apply the
  transformation \eqref{aut} with $a_1=1$, $a_2 = 1$, $a_4=
  \frac{\a_2}{\a_1}$, $B=\mathbbm{1}$ and all remaining entries zero. With an
  analogous argument, we can assume that $\b_1=0$, $\b_2 \ne 0$. Since
  $\ga_{2,2} \ne 0$ by non-degeneracy, we can rescale $\o$ such that
  $\ga_{2,2} =1$. Now, the transformation of the form \eqref{aut}
  given by
  \begin{eqnarray*}
    && a_{1} = 1, \:\: a_{2} = 0, \:\: a_{3} = 0, \:\: a_{4} = -\beta_{2}, \:\: 
    b_{1} = 1, \:\: b_{2} = 0, \:\: b_{3} = 0, \:\: b_{4} = -\alpha_{2}, \:\: a_{5} = 0, \:\: \\
    && a_{6} = -\frac{\alpha_{3}\beta_{2}}{\alpha_{2}}, \:\:
    b_{5} = 0, \:\: b_{6} = -\frac{\alpha_{2}\beta_{3}}{\beta_{2}}, \:\: 
    c_{1} = \frac{\gamma_{1, 1}}{\alpha_{2}}, \:\: c_{2} = -\gamma_{1, 2}, \:\: d_{1} = 0, \:\: d_{2} = \gamma_{2, 1}, 
  \end{eqnarray*}
  maps $\omega$ to a multiple of the normal form $\omega_2$.

  Thirdly, we assume that $\o$ is non-degenerate with $\o^{\k_1}=0$,
  i.e.\ $\ga_{3,3}=0$ and both $\o^{\k_2} \ne 0$, i.e.\ $\ga_{1,3}$ or
  $\ga_{2,3}\ne 0$, and $\o^{\k_3} \ne 0$, i.e.\ $\ga_{3,1}$ or
  $\ga_{3,2}\ne 0$. Similar as before, we can achieve
  $\ga_{2,3}=0$, $\ga_{1,3} \ne 0$ by applying, if $\ga_{2,3} \ne 0$, the
  transformation \eqref{aut} with $a_1=1$, $a_2 = 1$, $a_4= -
  \frac{\ga_{1, 3}}{\ga_{2, 3}}$, $B=\mathbbm{1}$ and all remaining entries
  zero. Analogously, we can assume $\ga_{3,2} =0$, $\ga_{3,1} \ne 0$ and
  rescaling yields $\ga_{2,2}=1$, which is non-zero by
  non-degeneracy. After this simplification, the condition $d \o ^2=0$
  implies that $\a_1=\b_1=0$ and the transformation
  \begin{eqnarray*}
    a_{1} &=& 1, \:\: a_{2} = 0, \:\: a_{3} = \frac{\alpha_{2}\beta_{3}-\gamma_{3, 1}\gamma_{1, 2}}{\gamma_{3, 1}}, \:\: a_{4} = \gamma_{1, 3}, \:\: a_{5} = 0, \:\: a_{6} = 0, \:\: b_{1} = 1, \:\: b_{2} = 0, \\
    b_{3} &=& \frac{\beta_{2}\alpha_{3}-\gamma_{1, 3}\gamma_{2, 1}}{\gamma_{1, 3}}, \:\:  b_{4} = \gamma_{3, 1}, \:\: b_{5} = \frac{\gamma_{1, 2}\gamma_{1, 3}\gamma_{2, 1}\gamma_{3, 1}-\gamma_{1, 1}\gamma_{1, 3}\gamma_{3, 1}-\alpha_{2}\alpha_{3}\beta_{2}\beta_{3}}{\gamma_{1, 3}^2\gamma_{3, 1}}, \:\: b_{6} = 0, \:\: \\
    c_{1} &=& 0, \:\: c_{2} = \beta_{3}, \:\: d_{1} = 0, \:\: d_{2} = -\alpha_{3}, 
  \end{eqnarray*}
  maps $\o$ to $\tilde \o = e^1f^3 + e^2f^2 + e^3f^1 + \tilde \alpha_2
  e^{31} + \tilde \beta_2 f^{31}$. The condition $\omega^{\k_4}=0$
  corresponds to $\tilde \alpha_2=0$, $\tilde \b_2=0$, i.e.\ normal
  form $\omega_3$. If $\omega^{\k_4} \ne 0$, we can achieve $\tilde
  \a_2 \ne 0$ by possibly changing the summands. Now, the
  transformation
  \begin{eqnarray*}
    a_1 = 1,\:\: a_2 = 0,\:\: a_3 = 0 ,\:\: a_4 = -\frac{1}{\tilde \a_2},\:\:  a_5 = 0,\:\: a_6 = 0,\:\: c_1 = 0,\:\:  c_2 = 0, \:\: \\
    b_1 = -\frac{1}{\tilde \a_2},\:\: b_2 = 0,\:\: b_3 = 0 ,\:\: b_4 = -\frac{1}{\tilde \a_2},\:\: b_5 = 0,\:\:  b_6 = 0,\:\:  d_1 = 0,\:\: d_2 = 0,\:\: 
  \end{eqnarray*}
  maps $\tilde \o$ to the fourth normal form $\o_4$.

  The cases that remain are $\o^{\k_1} = 0$ and either $\o^{\k_2} \ne
  0, \o^{\k_3} = 0$ or $\o^{\k_3} = 0,\o^{\k_2} \ne 0$. After changing
  the summands if necessary, we can assume $\o^{\k_3} = 0$ and
  $\o^{\k_2} \ne 0$, i.e. $\ga_{3,1}=\ga_{3,2}=\ga_{3,3} =0$ and at
  least one of $\ga_{1,3}$ or $\ga_{2,3}$ non-zero. As before, we can
  achieve $\ga_{2,3}=0$ by the transformation $a_1=1, a_2=1, a_4 =
  -\frac{\ga_{1,3}}{\ga_{2,3}}$. Evaluating $d \o^2 =0$ yields $\a_1
  =0$. Now, non-degeneracy enforces that $\beta_1 \ne 0$ or $\beta_2
  \ne 0$, and after another similar transformation
  $\beta_1=0$. Finally, the simplified $\o$ is non-degenerate if and
  only if $\ga_{2,2}\a_2\b_2 \ne 0$ and, after rescaling such that
  $\ga_{2,2}=1$, the transformation
  \begin{eqnarray*}
    a_{1} &=& 1,\:\:  a_{2} = 0,\:\:  a_{3} = 0,\:\:   
    a_{4} = - \frac{\gamma_{1, 3}^2}{\beta_{2}},\:\:  
    a_{5} = 0,\:\:  
    a_{6} = \frac{\gamma_{1, 3}^2(\gamma_{1, 3}\gamma_{2, 1}-\alpha_{3}\beta_{2})}{\alpha_{2}\beta_{2}^2},\:\:  \\
    b_{1} &=& -\frac{\gamma_{1, 3}}{\beta_{2}},\:\:  b_{2} = 0,\:\:   
    b_{3} = 0,\:\:  b_{4} = -\alpha_{2},\:\:  b_{5} = 0,\:\: 
    b_{6} = -\frac{\beta_{3}\alpha_{2}}{\beta_{2}},\:\:  \\
    c_{1} &=& 0,\:\: 
    c_{2} = -\frac{\gamma_{1, 2}\beta_{2}+\gamma_{1, 3}\beta_{3}}{\beta_{2}} ,\:\:  
    d_{1} = -\frac{\gamma_{1, 1}}{\beta_{2}} ,\:\: 
    d_{2} = \frac{\gamma_{1, 3}^2 \gamma_{2, 1}}{\beta_{2}^2}, 
  \end{eqnarray*}
  maps $\o$ to a multiple of the fifth normal form $\o_5$.
\end{proof}

Using this lemma, it is possible to describe all half-flat structures
$(\o,\rho)$ on $\h_3 \op \h_3$ as follows. In a fixed standard basis
such that $\o$ is in one of the normal forms, the equations $d\rho=0$
and $\o \wedge \rho =0$ are linear in the coefficients of an arbitrary
three-form $\rho$. Thus, it is straightforward to write down all
compatible closed three-forms for each normal form which depend on
nine parameters in each case. The stable forms in this
nine-dimensional space are parametrised by the complement of the
zero-set of the polynomial $\lambda(\rho)$ of order four. One
parameter is eliminated when we require a stable $\rho$ to be
normalised in the sense of \eqref{normrho}. We remark that the
computation of the induced tensors $J_\rho$, $\hat \rho$ and
$g_{(\o,\rho)}$ may require computer support, in particular, the
signature of the metric is not obvious. However, stability is an open
condition: If a single half-flat structure $(\o_0,\rho_0)$ is
explicitly given such that $\o_0$ is one of the normal forms, then the
eight-parameter family of normalised compatible closed forms defines a
deformation of the given half-flat structure $(\o_0,\rho_0)$ in some
neighbourhood of $(\o_0,\rho_0)$.

For instance, the closed three-forms which are compatible with the
first normal form
\begin{equation}
  \label{o1}
  \omega=e^{1}f^1 + e^2f^2 + e^3f^3
\end{equation}
in a standard basis can be parametrised as follows:
\begin{eqnarray}
  \label{rho1}
  \rho = \rho (a_1, \dots , a_9) &=& a_{1} \:e^{123}+ a_{2} \:f^{123} + a_3 \:e^{1}f^{23} +a_4 \:e^{2}f^{13} + a_5 \:e^{23}f^{1} + a_6 \:e^{13}f^{2}\\
  &+& a_7 \:(e^{2}f^{23}-\:e^{1}f^{13}) + a_8 \:(e^{12}f^{3}-e^{3}f^{12}) + a_9 \:(e^{23}f^{2}- \:e^{13}f^{1}). \nonumber
\end{eqnarray}
The quartic invariant $\lambda(\rho)$ depending on the nine parameters
is
\begin{eqnarray*}
  \lambda(\rho)&=&(2 a_6 a_4 a_8^2+2 a_1 a_2 a_8^2+2 a_8^2 a_3 a_5-4 a_5 a_7^2 a_6-4 a_9^2 a_4 a_3-4 a_9^2 a_2 a_8+4 a_7^2 a_8 a_1\\ 
  &+&4 a_7 a_8^2 a_9 +a_1^2 a_2^2+a_6^2 a_4^2+a_3^2 a_5^2+a_8^4-2 a_6 a_4 a_3 a_5+4 a_5 a_7 a_9 a_3+4 a_9 a_4 a_6 a_7\\
  &-&4 a_5 a_2 a_6 a_8 + 4 a_4 a_8 a_1 a_3-4 a_9 a_2 a_1 a_7-2 a_1 a_2 a_6 a_4-2 a_1 a_2 a_3 a_5) \, (\e^{123}f^{123})^{\otimes 2}.
\end{eqnarray*} 
\begin{example}
  \label{3ex}
  For each possible signature, we give an explicit normalised
  half-flat structure with fundamental two-form \eqref{o1}. The first
  and the third example appear in \cite{SH}. To begin with, the closed
  three-form
  \begin{eqnarray}
    \label{ex1}
    \rho  &=&  \frac{1}{\sqrt 2} (e^{123} - f^{123} - e^{1}f^{23} + e^{23}f^{1} - e^{2}f^{31} + e^{31}f^{2} - e^{3}f^{12} + e^{12}f^{3})
  \end{eqnarray}
  induces a half-flat $\SU(3)$-structure $(\omega,\rho)$ such that the
  standard basis is orthonormal. Similarly, the closed three-form
  \begin{eqnarray}
    \label{ex2}
    \rho  &=&  \frac{1}{\sqrt 2} (e^{123} - f^{123} - e^{1}f^{23} + e^{23}f^{1} + e^{2}f^{31} - e^{31}f^{2} + e^{3}f^{12} - e^{12}f^{3})
  \end{eqnarray}
  induces a half-flat $\SU(1,2)$-structure $(\omega,\rho)$ such that
  the standard basis is pseudo-orthonormal with $e_1$ and $e_4$ being
  spacelike. Finally, the closed three-form
  \begin{eqnarray}
    \label{ex3}
    \rho & = & \sqrt 2\, (e^{123} + f^{123}),
  \end{eqnarray}
  induces a half-flat $\SL(3,\bR)$-structure $(\omega,\rho)$ such that
  the two $\h_3$-summands are the eigenspaces of the para-complex
  structure $J_\rho$, which is integrable since also $d \hat \rho =
  0$. The induced metric is
  \begin{eqnarray*}
    g &=& 2 \left(\, e^1 \cdot e^4 + e^2 \cdot e^5 + e^3 \cdot e^6  \right).
  \end{eqnarray*}
\end{example}

In fact, half-flat structures with Riemannian metrics are only
possible if $\o$ belongs to the orbit of the first normal form.
\begin{lemma}
  Let $(\o,\rho)$ be a half-flat $\SU(3)$-structure on $\h_3 \op
  \h_3$. Then it holds $\o^{\k_1} \ne 0$. In particular, there is a
  standard basis such that $\o= \o_1 = e^{1}f^1 + e^2f^2 + e^3f^3$.
\end{lemma}
\begin{proof}
  Suppose that $(\o,\rho)$ is a half-flat $\SU(3)$-structure on $\h_3
  \op \h_3$ with $\o^{\k_1} = 0$. Thus, we can choose a standard basis
  such that $\o$ is in one of the normal forms $\o_2, \dots, \o_5$ of
  Lemma \ref{onormal} and $\rho$ belongs to the corresponding
  nine-parameter family of compatible closed three-forms. We claim
  that the basis one-form $e^1$ is isotropic in all four cases which
  yields a contradiction since the metric of an $\SU(3)$-structure is
  positive definite. The quickest way to verify the claim is the
  direct computation of the induced metric, which depends on nine
  parameters, with the help of a computer. In order to verify the
  assertion by hand, the following formulas shorten the calculation
  considerably. For all one-forms $\a$, $\b$ and all vectors $v$, the
  $\ve$-complex structure $J_\rho$ and the metric $g$ induced by a
  compatible pair $(\o,\rho)$ of stable forms satisfy
  \begin{eqnarray*}
    \a \wedge J_\rho^* \beta \wedge \o^2 &=& g(\a,\b) \frac{1}{3} \o^3, \\
    J_\rho^* \a (v) \phi (\rho) &=& \a \wedge \rho \wedge (v \inter \rho),
  \end{eqnarray*}
  which is straightforward to verify in the standard basis
  \eqref{normal3form}, \eqref{normalomega}, cf. also \cite[Lemmas 2.1,
  2.2]{SH}. For instance, for the second normal form $\o_2$, it holds
  $e^1 \wedge \o_2^2 = - 2 e^{12}f^{123}$. Thus, by the first formula,
  it suffices to show that $ J_\rho^* e^1 (e_3) = e^1(J_\rho e_3) = 0$
  which is in turn satisfied if $\e^1 \wedge \rho \wedge (e_3 \inter
  \rho)=0$ due to the second formula. A similar simplification applies
  to the other normal forms and we omit the straightforward
  calculations.
\end{proof}
Moreover, the geometry turns out to be very rigid if $\o^{\k_1}=0$. We
recall that simply connected para-hyper-K\"ahler symmetric spaces with
abelian holonomy are classified in \cite{ABCV}, \cite{C}. In
particular, there exists a unique simply connected four-dimensional
para-hyper-K\"ahler symmetric space with one-dimensional holonomy
group, which is defined in \cite{ABCV}, Section 4. We denote the
underlying pseudo-Riemannian manifold as $(N^4,g_{PHK})$.
\begin{Prop}
  \label{k1=0}
  Let $(\o,\rho)$ be a left-invariant half-flat structure with
  $\o^{\k_1} =0$ on $H_3 \times H_3$ and let $g$ be the
  pseudo-Riemannian metric induced by $(\o,\rho)$. Then, the
  pseudo-Riemannian manifold $(H_3 \times H_3,g)$ is either flat or
  isometric to the product of $(N^4,g_{PHK})$ and a two-dimensional
  flat factor. In particular, the metric $g$ is Ricci-flat.
\end{Prop}
\begin{proof}
  Due to the assumption $\o^{\k_1} =0$, we can choose a standard basis
  such that $\o$ is in one of the normal forms $\o_2, \dots, \o_5$. In
  each case separately, we do the following. We write down all
  compatible closed three-forms $\rho$ depending on nine
  parameters. With computer support, we calculate the induced metric
  $g$. For the curvature considerations, it suffices to work up to a
  constant such that we can ignore the rescaling by $\lambda(\rho)$
  which is different from zero by assumption. Now, we transform the
  left-invariant co-frame $\{ e^1, \dots, f^3 \}$ to a coordinate
  co-frame $\{ dx_1, \dots, dy_3 \}$ by applying the transformation
  defined by
  \begin{equation}
    \label{trafo}
    e^1 = dx_1 ,\: e^2 = dx_2,\: e^3 =  dx_3 + x_1 dx_2, \: f^1 = dy_1 ,\: f^2 = dy_2,\: f^3 =  dy_3 + y_1 dy_2,
  \end{equation}
  such that the metric is accessible for any of the numerous packages
  computing curvature. The resulting curvature tensor $R \in \Gamma(\mathrm{End}\, \L^2 TM)$, $M = H^3 \times H^3$, has in each case only one
  non-trivial component
  \begin{equation}
    \label{Rpara}
    R( \partial_{x_1} \wedge \partial_{y_1}) = c \, \partial_{x_3} \wedge \partial_{y_3}  
  \end{equation}
  for a constant $c \in \bR$ and $R$ is always parallel. Thus, the
  metric is flat if $c=0$ and symmetric with one-dimensional holonomy
  group if $c\ne 0$, for $H_3 \times H_3$ is simply connected and a
  naturally reductive homogeneous metric is complete.

  Furthermore, it turns out that the metric restricted to $TN:=span
  \{ \partial_{x_1}, \partial_{x_3}, \partial_{y_1}, \partial_{y_3}\}$
  is non-degenerate and of signature $(2,2)$ for all parameter
  values. Thus, the manifold splits in a four-dimensional symmetric
  factor with neutral metric and curvature tensor \eqref{Rpara} and
  the two-dimensional orthogonal complement which is flat. Since a
  simply connected symmetric space is completely determined by its
  curvature tensor and the four-dimensional para-hyper-K\"ahler
  symmetric space $(N^4,g_{PHK})$ has the same signature and curvature
  tensor, the four-dimensional factor is isometric to
  $(N^4,g_{PHK})$. Finally, the metric $g$ is Ricci-flat since
  $g_{PHK}$ is Ricci-flat.
\end{proof}
\begin{example}
  The following examples define half-flat normalised
  $\SU(1,2)$-structures with $\o^{\k_1} =0$ in a standard basis. None
  of the examples is flat. Thus, the four structures are equivalent as
  $\SO(2,4)$-structures due to Proposition \ref{k1=0}, but the
  examples show that the geometry of the reduction to $\SU(1,2)$ is
  not as rigid.
  \begin{eqnarray*}
    \omega = \omega_2 \, , \quad \rho &=& e^{12}f^{3}+\sqrt 2e^{13}f^{2}+e^{1}f^{23}+e^{23}f^{1}-e^{3}f^{12}+\sqrt 2f^{123} , \\
    g &=& - \, ( e^2 )^2 - \, ( f^2 )^2 +2 \, e^1 \ccdot e^3 -2 \sqrt 2 \, e^1 \ccdot f^3 + 2 \sqrt 2 \, e^3 \ccdot f^1 -2 \, f^1 \ccdot f^3 , \\
    && \mbox{(Ricci-flat pseudo-K\"ahler since $d\omega = 0$, $d \hat \rho = 0$);}\\
    \omega = \omega_3 \, , \quad \rho &=& e^{123}+e^{12}f^{3}+e^{13}f^{2}+e^{1}f^{12}-2 e^{1}f^{23}+e^{2}f^{13}-e^{3}f^{12}, \\
    g&=& - \, ( e^2 )^2 -2 \, ( f^2 )^2 +2 \, e^1 \ccdot f^1 +2 \, e^1 \ccdot f^3 +2 \, e^2 \ccdot f^2 -2 \, e^3 \ccdot f^1 -2 \, f^1 \ccdot f^3 , \\
    && \mbox{($d\o \ne 0$, $J_\rho$ integrable since $d \hat \rho = 0$);}\\
    \omega = \omega_4 \, , \quad 
    \rho &=& \b\, e^{12}f^{3} - \b\, e^{13}f^{2} + \b\, e^{1}f^{23}
    + \frac{\b+1}{\b^3} \, e^{23}f^{1} + \frac{\b^4-\b-1}{\b^3}\, e^{2}f^{13} \\
    && \quad - \b \,e^{3}f^{12} - (\b^2+2 \b) \,f^{123} , 
    \qquad \qquad \mbox{($d\o \ne 0$, $d \hat \rho \ne 0$),}\\
    g&=&  - \frac{1}{\b^2} \, ( e^2 )^2 - \b^2 \, ( f^2 )^2 
    + 2 \b^2 \, e^1 \ccdot f^3 - \frac{2}{\b^2(\b+1)} \, e^3 \ccdot f^1
    - \frac{2(\b^4+\b+1)}{\b^2} \, f^1 \ccdot f^3 ; \\
    \omega = \omega_5 \, , \quad  
    \rho &=& e^{12}f^{3}+e^{13}f^{2}-e^{1}f^{23}+e^{23}f^{1}-e^{3}f^{12}+f^{123} , 
    \qquad \qquad \mbox{($d\o \ne 0$, $d \hat \rho \ne 0$),}\\
    g&=& - \, ( e^2 )^2 -2 \, ( f^2 )^2 +2 \, e^1 \ccdot e^3 +2 \, e^2 \ccdot f^2 +2 \, f^1 \ccdot f^3 .
  \end{eqnarray*}
\end{example}
\begin{example}
  Moreover, we give examples of half-flat normalised
  $\SL(3,\bR)$-structures with $\o^{\k_1} =0$. Again, none of the
  structures is flat.
  \begin{eqnarray*}
    \omega = \omega_2 \, , \quad 
    \rho &=& \sqrt 2 \, ( \,e^{1}f^{23} + e^{23}f^{1} \,) , 
    \qquad \qquad \mbox{($d\o = 0$, $d \hat \rho = 0$),} \\
    g&=& 2 \, e^1 \ccdot e^3 -2 \, e^2 \ccdot f^2 -2 \, f^1 \ccdot f^3 ; \\
    \omega = \omega_3 \, , \quad  
    \rho &=& \sqrt 2 \, ( \, e^{12}f^{3} + e^{13}f^{2} + e^{1}f^{12} - e^{3}f^{12} \, ) , 
    \qquad \qquad \mbox{($d\o \ne 0$, $d \hat \rho \ne 0$),}\\
    g&=& -2 \, ( e^1 )^2 +2 \, e^1 \ccdot e^3 -2 \, e^1 \ccdot f^3 +2 \, e^2 \ccdot f^2 -2 \, f^1 \ccdot f^3 ; \\
    \omega = \omega_4 \, , \quad 
    \rho &=& -\sqrt{2\b +2} \, ( \:e^{12}f^{3}-e^{1}f^{23}+e^{2}f^{13}-e^{3}f^{12} ), 
    \qquad \qquad \mbox{($d\o \ne 0$, $d \hat \rho \ne 0$),} \\
    g&=& -2 \, ( f^2 )^2 +2 \, e^1 \ccdot e^3 +2 \, e^1 \ccdot f^3 +2 \, e^2 \ccdot f^2 -2 \, e^3 \ccdot f^1 - (2\b+4) \, f^1 \ccdot f^3 ; \\
    \omega = \omega_5 \, , \quad 
    \rho &=& \sqrt 2\,( \, e^{123} + f^{123}\, ), 
    \qquad \qquad \mbox{($d\o \ne 0$, $d \hat \rho = 0$),} \\
    g&=& 2 \, e^1 \ccdot f^3 +2 \, e^2 \ccdot f^2 +2 \, e^3 \ccdot f^1 . \\
  \end{eqnarray*}
\end{example}
\subsection{Solving the evolution equations on $H_3 \times H_3$}
\label{examplesH3H3}
Due to the preparatory work of the Lemmas \ref{trick17} and
\ref{onormal}, it turns out to be possible to explicitly evolve every
half-flat structure on $\h_3 \op \h_3$ without integrating.
\begin{Prop}
  \label{afflin}
  Let $(\omega_0,\rho_0)$ be any half-flat $H^{\ve,\tau}$-structure on
  $\h_3 \oplus \h_3$ with $\omega_0^{\k_1} = 0$. Then, the solution of
  the evolution equations \eqref{evnochmal} is affine linear in the
  sense that
  \begin{eqnarray}
    \label{linearsolution}
    \sigma(t) = \sigma_0 + t \, d\hat \rho_0\, , \qquad \rho(t) = \rho_0 + t \, d\o_0
  \end{eqnarray}
  and is well-defined for all $t \in \bR$.
\end{Prop}
\begin{proof}
  Let $\{ e_1, \dots , f_3 \}$ be a standard basis such that
  $\omega_0$ is in one of the normal forms $\omega_2, \dots, \omega_5$
  of Lemma \ref{onormal} which satisfy $\omega_0^{\k_1} = 0$. By Lemma
  \ref{trick17} and the second evolution equation, we know that there
  is a function $y(t)$ with $y(0)=0$ such that
  \[\sigma(t) = \sigma_0 + y(t) e^{12}f^{12} = \frac{1}{2} \o_0^2 +
  y(t) e^{12}f^{12}.\] For each of the four normal forms, the unique
  two-form $\o(t)$ with $\frac{1}{2} \o(t)^2 = \sigma(t)$ and
  $\o(0)=\o_0$ is
  \[ \o(t) = \o_0 - y(t) e^1 f^1 .\] However, the two-form $e^1f^1$ is
  closed such that the exterior derivative $d\o(t) = d\o_0$ is
  constant. Therefore, we have $\rho(t) = \rho_0 + t \, d\o_0$ by the
  first evolution equation. Moreover, the two-form $\o(t)$ is stable
  for all $t \in \bR$ since it holds $\phi(\o(t))=\phi(\o_0)$ for each
  of the normal forms and for all $t \in \bR$. It remains to show that
  $d\hat \rho (t)$ is constant in all four cases which implies that
  the function $y(t)$ is linear by the second evolution equation.

  As explained in section \ref{halfflatonH3H3}, it is easy to write
  down, for each normal form $\o_0$ separately, all compatible, closed
  three-forms $\rho_0$, which depend on nine parameters. For $\rho(t)
  = \rho_0 + t \, d\o_0$, we verify with the help of a computer that
  $\lambda(\rho(t)) = \lambda(\rho_0)$ is constant such that $\rho(t)$
  is stable for all $t \in \bR$ since $\rho_0$ is stable. When we also
  calculate $J_{\rho(t)}$ and $\hat \rho (t) = J^*_{\rho(t)} \rho(t)$,
  it turns out in all four cases that $d\hat \rho (t)$ is
  constant. This finishes the proof.
\end{proof}
We cannot expect that this affine linear evolution of spaces which
have one-dimensional holonomy, due to Proposition \ref{k1=0}, yields
metrics with full holonomy $\G_2^{*}$. Indeed, due to the following
result the geometry does not change significantly compared to the
six-manifold.
\begin{Cor}
  Let $(\omega_0,\rho_0)$ be a half-flat $H^{\ve,\tau}$-structure on
  $\h_3 \oplus \h_3$ with $\omega_0^{\k_1} = 0$ and let $g_\vf$ be the
  Ricci-flat metric induced by the parallel stable three-form $\vf$ on
  $M\times \bR$ defined by the solution \eqref{linearsolution} of the
  evolution equations with initial value $(\omega_0,\rho_0)$.  Then,
  the pseudo-Riemannian manifold $(M \times \bR,g_\vf)$ is either flat
  or isometric to the
  product of the four-dimensional para-hyper-K\"ahler symmetric space
  $(N^4,g_{PHK})$ and a three-dimensional flat factor.
\end{Cor}
\begin{proof}
  By formula \eqref{gvarphi}, the metric $g_\vf$ is determined by the
  time-dependent metric $g(t)$ induced by $(\o(t),\rho(t))$. All
  assertions follow from the analysis of the curvature of $g_\vf$
  completely analogous to the proof of Proposition \ref{k1=0}.
\end{proof}
The situation changes completely when we consider the first normal
form $\o_1$ of Lemma \ref{onormal}.
\begin{Prop}
  \label{k1ne0}
  Let $(\omega_0,\rho_0)$ be any normalised half-flat
  $H^{\ve,\tau}$-structure on $\h_3 \oplus \h_3$ with $\omega_0^{\k_1}
  \ne 0$. There is always a standard basis $\{ e_1, \dots , f_3 \}$
  such that $\omega_0=e^{1}f^1 + e^2f^2 + e^3f^3$.  In such a basis,
  we define $(\omega(x),\rho(x))$ by
  \begin{eqnarray*}
    \rho(x) &=& \rho_0 + x ( e^{12}f^{3} - e^{3}f^{12} ), \\
    \omega (x) &=& 2 \, (\ve \kappa(x) )^{-\frac{1}{2}} \; \left( \; 
      \frac{1}{4} \ve \kappa(x) \, e^1 f^1 
      + \frac{1}{4} \ve \kappa(x) \, e^2 f^2 + e^3 f^3 \; \right),
  \end{eqnarray*} 
  where $\kappa(x) \, ({e^{123}f^{123}})^{\otimes 2} = \lambda
  (\rho(x))$. Furthermore, let $I$ be the maximal interval containing
  zero such that the polynomial $\kappa(x)$ of order four does not
  vanish for any $x \in I$. The parallel stable three-form
  \eqref{varphidt} on $M \times I$ obtained by evolving
  $(\omega_0,\rho_0)$ along the Hitchin flow \eqref{evnochmal} is
  \[ \varphi = \frac{1}{2} \sqrt{\ve \kappa(x)} \, \omega (x) \wedge
  dx + \rho (x). \] The metric induced by $\varphi$, which has
  holonomy contained in $G^{\ve,\tau}$, is by \eqref{gvarphidt} given
  as
  \begin{equation}
    \label{G2metric}
    g_\varphi = g(x) - \frac{1}{4} \kappa(x) dx^2,
  \end{equation}
  where $g(x)$ denotes the metric associated to $(\omega(x),\rho(x))$
  via \eqref{Jrho} and \eqref{inducedmetric}. The variable $x$ is
  related to the parameter $t$ of the Hitchin flow by the ordinary
  differential equation \eqref{ODE}.
\end{Prop}
\begin{proof}
  Since $\omega_0^{\k_1} \ne 0$, we can always choose a standard basis
  such that $\omega_0=e^{1}f^1 + e^2f^2 + e^3f^3$ is in the first
  normal form of Lemma \ref{onormal}. Then $\rho_0$ is of the form
  (\ref{rho1}).

  Moreover, by Lemma \ref{trick17}, there is a function $y(t)$ which
  is defined on an interval containing zero and satisfies $y(0)=0$
  such that the solution of the second evolution equation can be
  written
  \[\sigma(t) = \sigma_0 + y(t) e^{12}f^{12}. \]
  The unique $\omega(t)$ that satisfies $\omega(0) = \omega_0 $ and
  $\frac{1}{2}\omega(t)^2 = \sigma(t)$ for all $t$ is
  \[ \omega (t) = \sqrt{1-y(t)} \; e^1 f^1 + \sqrt{1-y(t)}\; e^2 f^2 +
  \frac{1}{ \sqrt{1-y(t)}}\; e^3 f^3. \] Since
  \begin{equation}
    \label{domegat}
    d \omega (t) = \frac{1}{ \sqrt{1-y(t)}} (e^{12} f^3 - e^3 f^{12}),  
  \end{equation}
  there is another function $x(t)$ with $x(0)=0$ such that the
  solution of the first evolution equation can be written
  \begin{equation}
    \label{rhot}
    \rho(t) = \rho_0 + x(t) (e^{12} f^3 - e^3 f^{12}). 
  \end{equation}
  This three-form is compatible with $\omega(t)$ for all $t$, as one
  can easily see from (\ref{rho1}). Furthermore, the solution is
  normalised by Theorem \ref{paralleltheo}, which implies
  \[ \sqrt{\ve \lambda(\rho(t)) } = \phi(\rho(t)) = 2 \phi (\omega(t))
  = - 2 \sqrt{1-y(t)} \, e^{123}f^{123}.\] Hence, we can eliminate
  $y(t)$ by
  \[ y(t) = 1 - \frac{1}{4} \ve \kappa(x(t)). \] We remark that the
  normalisation of $\rho_0 = \rho(0)$ corresponds to $\kappa(0)
  =4\ve$.  Comparing \eqref{domegat} and \eqref{rhot}, the evolution
  equations are equivalent to the single ordinary differential
  equation
  \begin{eqnarray}
    \label{ODE}
    \dot x = \frac{2}{\sqrt{\ve \kappa(x(t))}}
  \end{eqnarray}
  for the only remaining parameter $x(t)$. In fact, we do not need to
  solve this equation in order to compute the parallel $\G_2^ {(*)}$-form
  when we substitute the coordinate $t$ by $x$ via the local
  diffeomorphism $x(t)$ satisfying $dt = \frac{1}{2} \sqrt {\ve
    \kappa(x(t))} \; dx$. Inserting all substitutions into the
  formulas \eqref{varphidt} and \eqref{gvarphidt} for the stable
  three-form $\vf$ on $M \times I$ and the induced metric $g_\vf$, all
  assertions of the proposition follow immediately from Theorem
  \ref{paralleltheo}.
\end{proof}

\begin{example}
  \label{exfullG2}
  The invariant $\kappa(x)$ and the induced metric $g(x)$ for the three explicit half-flat structures of Example \ref{3ex} are the following. \\
  If $(\o_0,\rho_0)$ is the $\SU(3)$-structure \eqref{ex1}, it holds 
  \begin{eqnarray}
    \nonumber \kappa(x) &=& (x- \sqrt 2)^3(x+\sqrt 2), \qquad \qquad I = (-\sqrt 2,\sqrt 2), \\
    \nonumber
    g(x) &=& (1 - \frac{1}{2}\sqrt 2\, x)\left((e^1)^2 + (e^2)^2 \!-\! 4\kappa(x)^{-1} (e^3)^2 + (e^4)^2 + (e^5)^2 \!-\! 4\kappa(x)^{-1}(e^6)^2  \right) \\ 
    \nonumber &+& \sqrt 2\, x (1 - \frac{1}{2}\sqrt 2\, x)\left( e^1 \ccdot e^4 + e^2 \ccdot e^5 + 4\kappa(x)^{-1} e^3 \ccdot e^6\right).
  \end{eqnarray}
  If $(\o_0,\rho_0)$ is the $\SU(1,2)$-structure \eqref{ex2}, we have
  \begin{eqnarray*}
    \kappa(x) &=& (x- \sqrt 2)(x+\sqrt 2)^3, \qquad \qquad I = (-\sqrt 2,\sqrt 2), \\
    g(x) &=& (1 + \frac{1}{2}\sqrt 2\, x) \left( (e^1)^2 - (e^2)^2 + 4\kappa(x)^{-1}  (e^3)^2 + (e^4)^2 - (e^5)^2 + 4\kappa(x)^{-1} (e^6)^2  \right) \\ 
    &-& \sqrt 2\, x (1 + \frac{1}{2}\sqrt 2\, x) \left( e^1 \ccdot e^4 + e^2 \ccdot e^5 + 4\kappa(x)^{-1}  e^3 \ccdot e^6\right).
  \end{eqnarray*}
  And for the $\SL(3,\bR)$-structure \eqref{ex3}, it holds 
  \begin{eqnarray*}
    \kappa(x) &=& (2+x^2)^2 , \qquad \qquad I = \bR,\\
    g(x) &=& (2 + x^2) \left( e^1 \ccdot e^4 + e^2 \ccdot e^5 \right) + 4 (2-x^2) \kappa(x)^{-1}  e^3 \ccdot e^6 + 4 \sqrt 2 \,x \kappa(x)^{-1}   \left( (e^3)^2 -(e^6)^2  \right).
  \end{eqnarray*}
\end{example}
\smallskip
\begin{theorem}
  Let $(\o(x),\rho(x))$ be the solution of the Hitchin flow with one
  of the three half-flat structures $(\o_0,\rho_0)$ of Example
  \ref{3ex} as initial value (see Proposition \ref{k1ne0} for the
  explicit solution and Example \ref{exfullG2} for the corresponding
  metric $g(x)$, defined for $x \in I$).

  Then, the holonomy of the metric $g_\vf$ on $M \times I$ defined by
  formula \eqref{G2metric} equals $\G_2$ for the $\SU(3)$-structure
  $(\o_0,\rho_0)$ and $\G_2^*$ for the other two structures.

  Moreover, restricting the eight-parameter family of half-flat
  structures given by \eqref{rho1} to a small neighbourhood of the
  initial value $(\rho_0,\o_0)$ yields in each case an eight-parameter
  family of metrics of holonomy equal to $\G_2$ or $\G_2^*$.
\end{theorem}
\begin{proof}
  For all three cases, we can apply the transformation \eqref{trafo}
  and calculate the curvature $R$ of the metric $g_\varphi$ defined by
  \eqref{G2metric}. Carrying this out with the package ``tensor''
  contained in Maple 10, we obtained that the rank of the curvature
  viewed as endomorphism on two-vectors is $14$. This implies that the
  holonomy of $g_\varphi$ in fact equals $\G_2$ or $\G_2^*$.

  The assertion for the eight-parameter family is an immediate
  consequence. Indeed, by construction, the rank of the curvature
  endomorphism is bounded from above by $14$ and being of maximal rank
  is an open condition.
\end{proof}
To conclude this section we address the issue of completeness and use
the Riemannian family in Example \ref{exfullG2} and Corollary
\ref{completecoroll} to construct a complete conformally parallel
$\G_2$-metric on $\rr\times (\Gamma\backslash H_3\times H_3)$.
\begin{example}\label{completeexample}
  Let $H_3$ be the Heisenberg group and $N=\Gamma \backslash H_3\times
  H_3$ be a compact nilmanifold given by a lattice $\Gamma$. Let us
  denote by $x:I\to (-\sqrt{2},\sqrt{2})$ the maximal solution to the
  equation
  \[ \dot{x}(t) = \frac{2}{\sqrt{ (\sqrt{2}-
      x(t))^3(x(t)+\sqrt{2})}},\] with initial condition $x(0)=0$,
  defining the $t$-dependent family of Riemannian metrics \be g_t&=&
  \frac{\sqrt{2}-x(t)}{\sqrt{2} } \left( \, (e^1)^2 + (e^2)^2 +
    (e^4)^2 + (e^5)^2 \, \right) + x(t)\left( \sqrt{2}-x(t) \right)
  \left( \, e^1\cdot e^4 + e^2\cdot e^5 \, \right)
  \\
  && + \frac{2\sqrt{2} }{(\sqrt{2}-x(t))^2(x(t)+\sqrt{2})} \left( \,
    (e^3)^2+ (e^6)^2 \, \right)
  -\frac{4x(t)}{(\sqrt{2}-x(t))^2(x(t)+\sqrt{2})}e^3\cdot e^6.  \ee If
  $\vf : \rr\to I$ is a diffeomorphism, then the metric
  \[dr^2 + \frac{1}{\vf' (r)^2} g_{\vf(r)}\] is globally conformally
  parallel $\G_2$ and geodesically complete.
\end{example}

\section{Special geometry of real forms of 
the symplectic $\SL (6,\bC)$-module
$\wedge^3\bC^6$} 
Homogeneous projective special K\"ahler manifolds of semisimple groups
with possibly indefinite metric and compact stabiliser were classified
in \cite{AC}. This includes the case of manifolds with (positive or
negative) definite metrics, for which the stabiliser is automatically
compact.  Projective special K\"ahler manifolds with negative definite
metric play an important role in supergravity and string theory. The
space of local deformations of the complex structure of a Calabi Yau
three-fold, for instance, is an example of a projective special
K\"ahler manifold with negative definite metric. As a particular
result of the classification \cite{AC}, there is an interesting
one-to-one correspondence between complex simple Lie algebras $\gl$ of
type A, B, D, E, F and G and homogeneous projective special K\"ahler
manifolds of semisimple groups with negative definite metric. The
resulting spaces are certain Hermitian symmetric spaces of non-compact
type.  The homogeneous projective special K\"ahler manifold associated
to the complex simple Lie algebra of type $E_6$, for instance, is
precisely the Hermitian symmetric space $\SU(3,3)/\,{\rm
  S}(\U(3)\times \U(3))$.

Under the above assumptions, the homogeneous projective special
K\"ahler manifold $G/K$ is realised as an open orbit of a real
semisimple group $G$ acting on a smooth projective algebraic variety
$X\subset P(V)$, where $V$ is the complexification of a real
symplectic module $V_0$ of $G$ and the cone $C(X):=\{ v\in V \,|\, \pi (v)
\in X\} \subset V$ over $X$ is Lagrangian.  Here $\pi : V\setminus \{
0\}\rightarrow P(V)$ denotes the canonical projection.  In fact,
$C(X)$ is the orbit of the highest weight vector of the $G^\bC$-module
$V$ under the complexified group $G^\bC$.  In the case $G/K =
\SU(3,3)/\,{\rm S}(\U(3)\times \U(3)) $ the symplectic module is
given by $V=\wedge^3\bC^6$. It was shown in \cite{BC} that the real
symplectic $G$-module $V_0$ always admits a homogeneous quartic
invariant $\lambda$, which is related to the hyper-K\"ahler part of
the curvature tensor of a symmetric quaternionic K\"ahler manifold
associated to the given complex simple Lie algebra $\gl$. Moreover,
the level sets $\{ \lambda = c\}$ are proper affine hyperspheres for
$c\neq 0$ and the affine special K\"ahler manifold $M$ underlying the
homogeneous projective special K\"ahler manifold $\bar{M}=G/K$ can be
realised as one of the open orbits of $\bR^* \cdot G$ on $V_0$
\cite{BC}.

In the following we shall describe all real forms $(G,V_0)$ of the
$\SL(6,\bC )$-module $V=\wedge^3\bC^6$ and study the affine special
geometry of the corresponding open orbits of $\bR^* \cdot G$.  As a
consequence, we obtain a list of projective special K\"ahler
manifolds, which admit a transitive action of a real form of
$\SL(6,\bC )$ by automorphisms of the special K\"ahler structure.
Besides the unique stationary compact example
\[ \SU(3,3)/\,{\rm S} (\U(3)\times \U(3)),\] we obtain the
homogeneous projective special K\"ahler manifolds
\[ \SU(3,3)/\,{\rm S}(\U(2,1)\times \U(1,2)),\quad \SU(5,1)/\,{\rm
  S}(\U(3)\times \U(2,1)) \quad\mbox{and}\quad \SL(6,\bR )/\,(\U(1)\cdot
\SL(3,\bC)),\] which are symmetric spaces with indefinite metrics and
non-compact stabiliser. The Hermitian signature of the metric is
$(4,5)$, $(6,3)$ and $(3,6)$, respectively. The latter result $(3,6)$
corrects Proposition 7 in \cite{H1}, according to which the Hermitian
signature of the underlying affine special K\"ahler manifold
$\GL^+(6,\bR )/\,\SL(3,\bC)$ is $(1,9)$. The correct Hermitian signature
of the affine special K\"ahler manifold is $(4,6)$.

Finally, we find that one of the two open orbits of $\SL(6,\bR )$ on
$\wedge^3\bR^6$ carries affine special para-K\"ahler geometry, the
geometry of $N=2$ vector multiplets on Euclidian rather than
Minkowskian space-time \cite{CMMS}. The corresponding homogeneous
projective special para-K\"ahler manifold is the symmetric space
\[ \SL(6,\bR )/\,{\rm S} (\GL(3,\bR)\times \GL(3,\bR)).\]
\subsection{The symplectic $\SL(6,\bC)$-module $V=\wedge^3\bC^6$ and
  its Lagrangian cone $C(X)$ of highest weight vectors}
We consider the 20-dimensional irreducible $\SL(6,\bC)$-module
$V=\wedge^3\bC^6$ equipped with a generator $\nu$ of the line
$\wedge^6\bC^6$.  The choice of $\nu$ determines an
$\SL(6,\bC)$-invariant symplectic form $\Omega$, which given by
\begin{equation}\label{OEqu} \Omega (v,w)\nu =v\wedge w,\quad v,w\in
  V.
\end{equation}

The highest weight vectors in $V$ are precisely the non-zero
decomposable three-vectors. They form a cone $C(X)\subset V$ over a
smooth projective variety $X\subset P(V)$, namely the Grassmannian
$Gr_3(\bC^6)$ of complex three-planes in $\bC^6$.  The group
$\SL(6,\bC)$ acts transitively on the cone $C(X)$ and, hence, on the
compact variety
\[ X\cong \SL(6,\bC)/P \cong \SU(6)/\,{\rm S} (\U(3)\times \U(3)),\]
where $P=\SL(6,\bC)_x\subset \SL(6,\bC)$ is the stabiliser of a point
$x\in X$ (a parabolic subgroup).  \bp The cone $C(X)=\{ v\in
\wedge^3\bC^6\setminus \{ 0\}\,|\, v\; \mbox{is decomposable}\} \subset V$
is Lagrangian.  \ep

\pf Let $(e_1, \ldots ,e_6)$ be a basis of $\bC^6$ and put
$p=e_{123}$. Then
\[ T_pC(X)= {\rm span}\{ e_{ijk}\,|\, \# \{i,j,k\}\cap \{1,2,3\} \ge 2\}\]
is ten-dimensional and is clearly totally isotropic with respect to
$\Omega$.
\end{proof}
\subsection{Real forms $(G,V_0)$ of the complex module $(\SL(6,\bC),V)$} 
Let $G$ be a real form of the complex Lie group
$\SL(6,\bC)$. There exists a $G$-invariant real structure $\tau$ on 
$V=\wedge^3\bC^6$ if and only if $G = \SL(6,\bR)$, $\SU(3,3)$, $\SU(5,1)$. 
In the first case $\tau$ is simply complex conjugation with respect to
$V_0=\wedge^3\bR^6$. In order to describe the real structure in the other 
two cases, we first 
endow $\bC^6$ with the standard $\SU(p,q)$-invariant pseudo-Hermitian
form $\langle \cdot ,\cdot \rangle$. The pseudo-Hermitian
form $\langle \cdot ,\cdot \rangle$ on $\bC^6$ induces an $\SU(p,q)$-invariant  
pseudo-Hermitian form $\gamma$ on $V$ such that 
\begin{equation} \label{gammadetEqu} 
\gamma (v_1\wedge v_2\wedge v_3,w_1\wedge w_2\wedge w_3)= 
\det ( \langle v_i,w_j\rangle ),
\end{equation}
for all $v_1, \ldots , w_3\in \bC^6$.  
Then we define an $\SU(p,q)$-invariant anti-linear map $\tau : V \rightarrow V$
by the equation 
\[ \tau := \sqrt{-1}\gamma^{-1}\circ \Omega .\]
Notice that $\Omega :V\rightarrow V^*, v\mapsto \Omega (\cdot ,v)$ is linear, whereas 
$\gamma :V\rightarrow V^*, v\mapsto \gamma (\cdot ,v)$ and $\gamma^{-1} : V^* \rightarrow V$ are 
anti-linear.  

\bp The anti-linear map $\t$ is an $\SU(p,q)$-invariant real structure
on $V=\wedge^3\bC^6$ if and only if $p-q\equiv 0\pmod{4}$.  In that
case, the $\SU(p,q)$-invariant pseudo-Hermitian form $\gamma =
\sqrt{-1}\Omega (\cdot ,\tau \cdot )$ on $V$ has signature $(10,10)$.
Otherwise, $\tau$ is an $\SU(p,q)$-invariant quaternionic structure on
$V$.  \ep

\pf  We present the calculations in the relevant  
cases $(p,q)=(3,3)$ and $(p,q)=(5,1)$.
The calculations in the other cases are similar.

Case $(p,q)=(3,3)$. Let  $(e_1,e_2,e_3,f_1,f_2,f_3)$ be a
unitary basis of $(\bC^6,\langle \cdot ,\cdot \rangle ) =\bC^{3,3}$, 
such that $\langle e_i,e_i \rangle = -\langle f_i,f_i \rangle =1$.   
We consider the following basis of $V$:
\begin{eqnarray*} &&(e_{123}, e_1\wedge f_{12},e_1\wedge f_{13},e_1\wedge f_{23},
e_2\wedge f_{12},e_2\wedge f_{13},e_2\wedge f_{23},
e_3\wedge f_{12},e_3\wedge f_{13},e_3\wedge f_{23},\\
&& 
f_{123},e_{23}\wedge f_3,-e_{23}\wedge f_2,e_{23}\wedge f_1,-e_{13}\wedge f_3,
e_{13}\wedge f_2,-e_{13}\wedge f_1,
e_{12}\wedge f_3,-e_{12}\wedge f_2,e_{12}\wedge f_1).
\end{eqnarray*}  
With respect to that basis and $\nu=e_{123} \wedge f_{123}$ we have
\begin{equation} \label{gammaOmegaEqu} \gamma = \left( \begin{array}{cc} \mathbbm{1}_{10}&0\\
0&-\mathbbm{1}_{10}\end{array}\right),\quad \Omega = \left( \begin{array}{cc} 0&\mathbbm{1}_{10}\\
-\mathbbm{1}_{10}&0\end{array}\right),\quad 
 \t = \sqrt{-1}\gamma^{-1}\circ \Omega = \sqrt{-1} \left( \begin{array}{cc} 0&\mathbbm{1}_{10}\\
\mathbbm{1}_{10}&0\end{array}\right). \end{equation}
This implies $\tau^2= {\rm Id}$, since $\t$ is anti-linear. 

Case $(p,q)=(5,1)$. Let  $(e_1,\ldots ,e_5,f)$ be a
unitary basis of $(\bC^6,\langle \cdot ,\cdot \rangle ) =\bC^{5,1}$, 
such that $\langle e_i,e_i \rangle = -\langle f,f \rangle =1$. 
With respect to the basis 
\begin{eqnarray} \label{51Equ} &&(e_{123}, e_{124},e_{125},e_{134},
e_{135},e_{145},e_{234},
e_{235},e_{245},e_{345},\\
&& 
e_{45}\wedge f,-e_{35}\wedge f,e_{34}\wedge f,e_{25}\wedge f,-e_{24}\wedge f,
e_{23}\wedge f,-e_{15}\wedge f,
e_{14}\wedge f,-e_{13}\wedge f,e_{12}\wedge f)\nonumber
\end{eqnarray}    
of $V$ and $\nu=e_{12345} \wedge f$ we have again the formulas
\re{gammaOmegaEqu} and $\t^2={\rm Id}$.
\end{proof} 

\subsection{Classification of open $G$-orbits on the Grassmannian $X$ and 
corresponding special K\"ahler manifolds} \label{mainSec} 
For each of the real forms $(G,V_0)$ obtained in the previous section,
we will now describe all open orbits of the real simple Lie group 
$G$ on the Grassmannian 
$X=Gr_3(\bC^6)=\{E \subset \bC^6$ a three-dimensional subspace$\} 
\hookrightarrow P(V)$, 
$E\mapsto \wedge^3 E$. We will also describe the projective 
special K\"ahler structure of these orbits $\bar{M}\subset P(V)$ 
and the (affine)  
special K\"ahler structure of the corresponding 
cones $M=C(\bar{M})\subset V$. The resulting homogeneous projective
special K\"ahler manifolds are listed in Table \ref{Tab}.  
\begin{table}\begin{center}
\begin{tabular}
{|c   |c       | } \hline
$G/H$ & Hermitian signature\\\hline
$\SU(3,3)/\, {\rm S} (\U(3)\times \U(3))$ & (0,9)\\\hline
$\SU(3,3)/\, {\rm S} (\U(2,1)\times \U(1,2))$& (4,5)\\\hline
$\SU(5,1)/\, {\rm S} (\U(3)\times \U(2,1))$& (6,3)\\\hline
$\SL(6,\bR)/\,(\U(1)\cdot \SL(3,\bC))$&(3,6)\\
\hline
\end{tabular}\\[2ex]
\caption{\it Homogeneous projective special K\"ahler manifolds $\bar{M}=G/H$ of
real simple groups $G$ of type $A_5$. Notice that $\dim_{\bC} \bar{M}=9$.} 
\label{Tab}
\end{center}\end{table}
Let us first recall
some definitions and constructions from special K\"ahler geometry. 
\smallskip
\subsubsection{Basic facts about special K\"ahler manifolds} 
\bd A {\em (pseudo-)K\"ahler manifold} $(M,J,g)$ is a
pseudo-Riemannian manifold $(M,g)$ endowed with a parallel
skew-symmetric complex structure $J\in \Gamma ({\rm End}\,TM)$. The
symplectic form $\o = g(\cdot,J \cdot)$ is called {\em K\"ahler
  form}. A {\em special K\"ahler manifold} $(M,J,g,\n)$ is a
(pseudo-)K\"ahler manifold $(M,J,g)$ endowed with a flat torsion-free
connection $\n$ such that $\n \o =0$ and $d^{\n} J=0$, where $d^{\n}
J$ is the exterior covariant derivative of the vector valued one-form
$J$.

A {\em conical special K\"ahler manifold} $(M,J,g,\n, \xi )$ is a 
special K\"ahler manifold $(M,J,g,\n)$ endowed with a timelike or a 
spacelike vector field $\xi$ such that $\n \xi =D \xi = {\rm Id}$,
where $D$ is the Levi-Civita connection. The vector field $\xi$ is called
{\em Euler vector field}. 
\ed     
The vector fields $\xi$ and $J\xi$ generate a free holomorphic action of a 
two-dimensional Abelian Lie algebra. If the action can be 
integrated to a free holomorphic $\bC^*$-action such that the 
quotient map $M \rightarrow \bar{M}:=M/\bC^*$ is a holomorphic submersion, then 
$\bar{M}$ is called a  {\em projective special K\"ahler manifold}. 
We will see now that $\bar{M}$ carries a  
canonical (pseudo-)K\"ahler metric $\bar{g}$ compatible with the 
induced complex structure $J$ on $\bar{M}=M/\bC^*$.  
Multiplying the metric $g$ with $-1$ if necessary, we can assume that 
$g(\xi ,\xi ) >0$. Then $S=\{ p\in M\,|\, g(\xi (p),\xi (p) )=1\}$ is a
smooth hypersurface invariant under the isometric $S^1$-action generated by 
the Killing vector field $J\xi$ and we can recover $\bar{M}$ as
the base of the circle bundle $S \rightarrow \bar{M}=S/S^1$. Then $\bar{M}$
carries a unique pseudo-Riemannian metric $\bar{g}$ 
such that $S \rightarrow \bar{M}$ is a Riemannian submersion. 
$(\bar{M},J, \bar{g})$ is in fact the 
K\"ahler quotient of $(M,J,g)$ by the $S^1$-action generated by 
the Hamiltonian Killing vector field $J\xi$. 

Next we explain the extrinsic construction of special K\"ahler manifolds from
\cite{ACD}. 
Let $(V,\Omega )$ be a complex symplectic vector space of dimension $2n$ 
endowed with a real structure $\tau$ such that 
\begin{equation}\label{tauOEqu}\ol{\Omega (v,w)}=\Omega (\tau v, \tau w),\quad 
\mbox{for all}\quad v,w\in V. 
\end{equation}
Then the  pseudo-Hermitian form
\begin{equation} \label{gammaEq} \gamma := \sqrt{-1}\Omega (\cdot ,\tau \cdot )
\end{equation}
has signature $(n,n)$. 
\bd  A holomorphic immersion $\phi : M \rightarrow V$ from an n-dimensional complex 
manifold $(M,J)$ into $V$ is called 
\begin{itemize}
\item[(i)] {\em nondegenerate} if $\phi^*\gamma$ is nondegenerate,
\item[(ii)] {\em Lagrangian} if $\phi^*\Omega =0$ and 
\item[(iii)] {\em conical} if $\phi (p) \in d\phi(T_pM)$ 
and $\gamma (\phi (p),\phi (p))\neq 0$  for all $p\in M$. 
\end{itemize} 
\ed  

\begin{samepage}
\bt \cite{ACD} \label{ACDThm} \begin{itemize}
\item[(i)] Any nondegenerate Lagrangian immersion $\phi : M\rightarrow
  V$ induces on the complex manifold $(M,J)$ the structure of a
  special K\"ahler manifold $(M,J,g,\n )$, where $g={\rm Re}\,
  \phi^*\gamma$ and $\n$ is determined by the condition $\n \phi^* \a
  =0$ for all $\a \in V^*$ which are real valued on $V^\tau$.
\item[(ii)] Any conical nondegenerate Lagrangian immersion $\phi :
  M\rightarrow V$ induces on $(M,J)$ the structure of a conical
  special K\"ahler manifold $(M,J,g,\n ,\xi )$. The vector field $\xi$
  is determined by the condition $d\phi\,\xi(p) = \phi (p)$.
\end{itemize} 
\et 
\end{samepage}

\subsubsection{ The case $G=\SL(6,\bR)$} 
Using the complex conjugation $\tau : v\mapsto \bar{v}$ 
on $\bC^6$ we can decompose $X$ into 
$G$-invariant real algebraic subvarieties  
$X_{(k)} = X_{(k)}(\tau ) := \{ E \in X\,|\, \dim (E\cap \bar{E}) = k\}\subset X$,
where $E\subset \bC^6$ runs through all three-dimensional subspaces and
$k\in \{0,1,2,3\}$. Notice that only $X_{(0)} \subset X$ is open. 
\bp The group $\SL(6,\bR)$ acts transitively on the open real subvariety  
$X_{(0)}=\{ E\in X \,|\,E\cap \bar{E} =0\} \subset X$. 
\ep 
\pf Given bases $(e_1,e_2,e_3)$, $(e_1',e_2',e_3')$ of $E, E'\in X_{(0)}$,
respectively, let $\varphi$ be the linear transformation, which maps
the basis $(e_1,e_2,e_3,\bar{e}_1,\bar{e}_2,\bar{e}_3)$ of $\bC^6$
to the basis $(e_1',e_2',e_3',\bar{e}_1',\bar{e}_2',\bar{e}'_3)$ of $\bC^6$. 
Then $\varphi\in \SL(6,\bR)$ and $\varphi E =E'$. 
\end{proof} 
\bt \label{SL6Thm} The group $\SL(6,\bR)$ has a unique open orbit
$X_{(0)}\cong \SL(6,\bR)/\,(\U(1)\cdot \SL(3,\bC))$ on the highest
weight variety $P(V) \supset X\cong Gr_3(\bC^6)$ of $V=\wedge^3\bC^6$.
The cone $V \supset M=C(X_{(0)})\cong \GL^+(6,\bR)/\,\SL(3,\bC)$ carries
an $\SL(6,\bR)$-invariant special K\"ahler structure of Hermitian
signature $(4,6)$, which induces on $\bar{M}=X_{(0)}$ the structure of
a homogeneous projective special K\"ahler manifold $\bar{M}$ of
Hermitian signature $(3,6)$.  \et

\pf Let $e_1,e_2,e_3\in \bC^6$ be three vectors which span a
three-dimensional subspace $E\subset \bC^6$ such that $E\in
X_{(0)}$. Then $\bC^6=E\oplus \bar{E}$ and the tangent space of
$M=C(X_{(0)})$ at $p=e_{123}$ is given by $T_pM=\wedge^3E \oplus
\wedge^2E\wedge \bar{E}$. We choose the real generator $\nu =
\sqrt{-1}e_{123}\wedge \bar{e}_{123}\in \wedge^3\bR^6$ and compute
$\gamma=\sqrt{-1}\Omega (\cdot, \t \cdot )$ on $T_pM$ using the
formula \re{OEqu}. The matrix of $\gamma|_{T_pM}$ with respect to the
basis
$(e_{123},e_{12}\wedge\bar{e}_3,e_{13}\wedge\bar{e}_2,e_{23}\wedge\bar{e}_1,
e_{12}\wedge \bar{e}_1,e_{13}\wedge \bar{e}_1,e_{12}\wedge \bar{e}_2,
-e_{23}\wedge \bar{e}_3,e_{23}\wedge \bar{e}_2,e_{13}\wedge
\bar{e}_3)$ is given by
\begin{equation} \label{MatEqu} \left( \begin{array}{cccc} 1&0&0&0\\
0&-\mathbbm{1}_3&0&0\\
0&0&0&\mathbbm{1}_3\\
0&0&\mathbbm{1}_3&0
\end{array}
\right) .
\end{equation} 
This shows that $\gamma$ has signature $(4,6)$ on $T_pM$.
Since $\GL^+(6,\bR)$ acts transitively
on $M$ and preserves the pseudo-Hermitian
form $\gamma$ up to a positive factor ($\SL(6,\bR )$ acts isometrically), 
the signature of the restriction of $\gamma$
to $M$  
does not depend on the 
base point. This shows that the inclusion $M\subset V$ is a 
holomorphic conical \emph{nondegenerate} 
Lagrangian immersion.  By Theorem \ref{ACDThm}, it induces 
a conical special \mbox{(pseudo-)}K\"ahler structure $(J,g,\n ,\xi)$ on $M$. 
It follows that the image $\bar{M} = \pi (M) = X_{(0)}
\cong \SL(6,\bR)/\,(\U(1)\cdot \SL(3,\bC))$ of $M$ under the 
projection $\pi : V\setminus \{ 0\}\rightarrow P(V)$ is a homogeneous 
projective special K\"ahler manifold.  
The induced pseudo-Hermitian form $\bar{\gamma}$  
on $T_{\pi (p)}\bar{M}$
has signature $(3,6)$. 
The latter statement follows from formula $\bar{\gamma}(d\pi_p X,d\pi_p Y)=
\frac{\gamma (X,Y)}{\gamma (p,p)}$ for $X,Y\in T_pM \cap p^\perp\subset V$ 
(see \cite{AC}), since 
$\gamma (p,p)=1$ for $p=e_{123}$.  
\end{proof} 
\subsubsection{ The case $G=\SU(3,3)$} 
Using the pseudo-Hermitian form 
$h=\langle \cdot ,\cdot \rangle$ on $\bC^6$ invariant under $G=\SU(3,3)$ 
we can decompose
the Grassmannian $X=Gr_3(\bC^6)$ into 
the $G$-invariant real algebraic subvarieties  
$X_{(k)} = X_{(k)}(h) :=\{ E \in X\,|\, {\rm rk}(h|_E)=k\}$, $k\in \{ 0,1,2,3\}$. 
Notice that only $X_{(3)} \subset X$ is open and that it can be decomposed 
further according to the possible signatures of $h|_E$:
\[ X_{(s,t)} := \{ E \in X\,|\, E\; \mbox{has signature}\; (s,t)\},\]
where $(s,t)\in \{ (3,0), (2,1), (1,2), (0,3)\}$. 

\bt \label{su33Thm} The group $\SU(3,3)$ has precisely four open
orbits on the highest weight variety $P(V) \supset X\cong Gr_3(\bC^6)$
of $V=\wedge^3\bC^6$, namely $X_{(3,0)}$, $X_ {(2,1)}$, $X_{(1,2)}$
and $X_{(0,3)}$.  In all four cases the cone $M=C(X_{(s,t)})\subset V$
carries an $\SU(3,3)$-invariant special K\"ahler structure.
\[ C(X_{(3,0)})
\cong \bR^*\cdot \SU(3,3)/\,
\SU(3)\times \SU(3)\] 
has Hermitian signature $(1,9)$. $C(X_{(0,3)})
\cong \bR^*\cdot \SU(3,3)/\,
\SU(3)\times \SU(3)$ 
has Hermitian signature $(9,1)$. 
For $\{s,t\}=\{2,1\}$, 
\[ C(X_{(s,t)})\cong \bR^*\cdot \SU(3,3)/\,
\SU(2,1)\times \SU(1,2)\]  
has Hermitian signature $(5,5)$. 
In all cases, the conical special K\"ahler manifold $M=C(X_{(s,t)})$ 
induces on $\bar{M}=X_{(s,t)}$ the structure of a homogeneous projective
special K\"ahler manifold $\bar{M}$. For $\{s,t\}=\{3,0\}$, 
\[ \bar{M}=X_{(s,t)}
\cong X_{(3,0)} \cong 
\SU(3,3)/\, {\rm S} (\U(3)\times \U(3)) \]  
has Hermitian signature $(0,9)$, for 
$\{s,t\}=\{2,1\}$, 
\[ \bar{M}=X_{(s,t)} \cong X_ {(2,1)}
\cong \SU(3,3)/\, {\rm S} (\U(2,1)\times \U(1,2))\]  
has Hermitian signature $(4,5)$. 
The special K\"ahler manifolds $C(X_{(s,t)})$ and $C(X_{(t,s)})$ are
equivalent. In fact, they are 
related by a holomorphic $\n$-affine anti-isometry, which induces 
a holomorphic isometry between the corresponding
projective special K\"ahler manifolds. 
\et 

\pf $X_{(3)}\subset X$ is Zariski open and is decomposed into the four 
open (in the standard topology) orbits $X_{(s,t)}$ of $G=\SU(3,3)$. 
Let $(e_1,e_2,e_3,f_1,f_2,f_3)$ be a
unitary basis of $(\bC^6,h)$, 
such that $\langle e_i,e_i \rangle = -\langle f_i,f_i \rangle =1$.   
Then $X_{(3,0)}$, $X_{(2,1)}$, $X_{(1,2)}$, $X_{(0,3)}$ 
are the $G$-orbits of the lines generated by the elements 
$e_{123}, e_1\wedge f_{12}, e_{12}\wedge f_1$, $f_{123}\in V$,  respectively. 

For $p=e_{123}$ and $M=C(X_{(3,0)})$, we calculate 
\[ T_pM= {\rm span}\{ e_{123},f_i\wedge e_{jk}\}.\] {}From
\re{gammadetEqu} we see that the matrix of $\gamma|_{T_pM}$ with
respect to the basis $(e_{123},f_i\wedge e_{jk})$ is ${\rm
  diag}(1,-\mathbbm{1}_9)$.  Therefore $\gamma$ has signature $(1,9)$ on
$T_pM$.  Since $\gamma (e_{123},e_{123})=1$, the signature of the
induced pseudo-Hermitian form $\bar{\gamma}$ on $T_p\bar{M}$ is
$(0,9)$.

For $p=e_{12}\wedge f_1$ and $M=C(X_{(2,1)})$, we obtain 
\[ T_pM= {\rm span}\{ e_1\wedge f_{12},e_2\wedge f_{12}, e_2\wedge
f_{13}, e_1\wedge f_{13}, e_{123}, e_{ij}\wedge f_1, e_{12}\wedge
f_i\}.\] The matrix of $\gamma|_{T_pM}$ with respect to the above
basis is ${\rm diag}(\mathbbm{1}_5,-\mathbbm{1}_5)$.  Therefore $\gamma$ has signature
$(5,5)$ on $T_pM$. We have $\gamma (e_{12}\wedge f_1,e_{12}\wedge
f_1)=-1$ and the signature of the induced pseudo-Hermitian form
$\bar{\gamma}$ on $T_p\bar{M}$ is $(4,5)$.

The linear transformation which sends the vectors $e_i$ to $f_i$ and
$f_i$ to $e_i$ induces a linear map $\varphi : V \rightarrow V$, which
interchanges the cone $C(X_{(s,t)})$ with $C(X_{(t,s)})$ and maps
$\gamma$ to $-\gamma$. This shows that $\gamma$ has signature $(9,1)$
and $(5,5)$ on $C(X_{(0,3)})$ and $C(X_{(1,2)})$, respectively.  As a
consequence, the induced pseudo-Hermitian form $\bar{\gamma}$ on
$\bar{M}=X_{(0,3)}$, $X_{(1,2)}$, has still signature $(0,9)$ and
$(4,5)$, respectively.

It follows from these calculations that the  
inclusion $M=C(X_{(s,t)})\subset V$ is a holomorphic conical nondegenerate 
Lagrangian immersion. By Theorem \ref{ACDThm}, it induces 
a conical special \mbox{(pseudo-)}K\"ahler structure $(J,g,\n ,\xi)$ on $M$ and
$\bar{M} =X_{(s,t)}\subset P(V)$ is a projective special K\"ahler manifold.

The above linear anti-isometry $\varphi : (V,\gamma ) \rightarrow (V,\gamma )$ 
maps $\Omega$ to $-\Omega$,  and, 
hence, preserves the real structure $\tau= \sqrt{-1}\gamma^{-1}\circ \Omega$. 
As a result, it 
maps the special K\"ahler structure $(J,g,\n )$ of $C(X_{(s,t)})$ 
to $(J',-g',\n')$, where $(J',g',\n')$ is the special K\"ahler structure 
of $C(X_{(t,s)})$. In particular, it induces a holomorphic isometry 
$X_{(s,t)}\cong X_{(t,s)}$. 
\end{proof} 

\subsubsection{ The case $G=\SU(5,1)$} 
Let $h=\langle \cdot ,\cdot \rangle$ be the standard pseudo-Hermitian form 
of signature $(5,1)$ on $\bC^6$, which is invariant under $G=\SU(5,1)$. 
Let us fix a unitary basis  $(e_1,\ldots ,e_5,f)$ of 
$(\bC^6,h)$, such that $\langle e_i,e_i \rangle = -\langle f,f \rangle =1$.  
As in the previous subsection, $X=Gr_3(\bC^6)$ is decomposed into the 
$G$-invariant real algebraic subvarieties  
$X_{(k)} = X_{(k)}(h)$, of which $X_{(3)} \subset X$ is Zariski open.  $X_{(3)}$ 
is now the union of the two open $G$-orbits $X_{(3,0)}$ and $X_{(2,1)}$.  
$X_{(3,0)}$ is the orbit of the line $\bC e_{123}\in P(V)$ and
$X_{(2,1)}$ is the orbit of $\bC e_{45}\wedge f\in P(V)$. 

\bt The group $\SU(5,1)$ has precisely two open orbits on the highest
weight variety $P(V) \supset X\cong Gr_3(\bC^6)$ of $V=\wedge^3\bC^6$,
namely $X_{(3,0)}$ and $X_ {(2,1)}$.  In both cases the cone
$M=C(X_{(s,t)})\subset V$ carries an $\SU(5,1)$-invariant special
K\"ahler structure.
\[ C(X_{(3,0)})
\cong \bR^*\cdot \SU(5,1)/\,
\SU(3)\times \SU(2,1)\] 
has Hermitian signature $(7,3)$.  
\[ C(X_{(2,1)})\cong \bR^*\cdot \SU(5,1)/\,
\SU(3)\times \SU(2,1)\]  
has Hermitian signature $(3,7)$. 
In both cases, the conical special K\"ahler manifold $M=C(X_{(s,t)})$ 
induces on $\bar{M}=X_{(s,t)}$ the structure of a homogeneous projective
special K\"ahler manifold $\bar{M}$.  
\[  \bar{M}=X_{(3,0)} \cong 
\SU(5,1)/\, {\rm S} (\U(3)\times \U(2,1))\]  
and 
\[ \bar{M}=X_ {(2,1)}
\cong \SU(5,1)/\, {\rm S} (\U(3)\times \U(2,1))\]  
have both Hermitian signature $(6,3)$. 
The special K\"ahler manifolds $C(X_{(3,0)})$ and $C(X_{(2,1)})$ are
equivalent. In fact, they are 
related by a holomorphic $\n$-affine anti-isometry, which induces 
a holomorphic isometry between the corresponding
projective special K\"ahler manifolds. 
\et 

\pf For $p=e_{123}$ and $M=C(X_{(3,0)})$, we have
\[ T_pM= {\rm span}\{ e_{123},e_{234},e_{235},e_{134},
e_{135},e_{124},e_{125}, e_{12}\wedge f,e_{13}\wedge f,e_{23}\wedge
f\}\] and the restriction of $\gamma$ to $T_pM$ is represented by the
matrix ${\rm diag}(\mathbbm{1}_7,-\mathbbm{1}_3)$ with respect to the
above basis.  This shows that the inclusion $M=C(X_{(3,0)})\subset V$
is a holomorphic conical nondegenerate Lagrangian immersion. By
Theorem \ref{ACDThm}, it induces a conical special (pseudo-)K\"ahler
structure $(J,g,\n ,\xi)$ of Hermitian signature $(7,3)$ on $M$ and
$\bar{M} =X_{(3,0)}\subset P(V)$ is a projective special K\"ahler
manifold of Hermitian signature $(6,3)$.  The anti-isometry relating
$C(X_{(3,0)})$ and $C(X_{(2,1)})$ is induced by the linear map
$\varphi : V \rightarrow V$ which has the matrix
\[ \left( \begin{array}{cc} 0&\mathbbm{1}_{10}\\
\mathbbm{1}_{10}&0\end{array}\right), \]
with respect to the basis \re{51Equ}. 
\end{proof} 
\subsection{The homogeneous projective 
special para-K\"ahler manifold \\
$\SL(6,\bR )/\, {\rm S} (\GL(3,\bR )\times \GL(3,\bR ))$}
Let us first briefly recall
the necessary definitions and constructions from special 
para-K\"ahler geometry, see  
\cite{CMMS} for more details. 
\subsubsection{Basic facts about special para-K\"ahler manifolds} 
\bd A {\em para-K\"ahler manifold} $(M,J,g)$ is a pseudo-Riemannian manifold
$(M,g)$ endowed with a parallel skew-symmetric endomorphism field $J\in 
\Gamma ({\rm End}\,TM)$ such that $J^2={\rm Id}$. 
A {\em special para-K\"ahler manifold} $(M,J,g,\n)$ is a 
para-K\"ahler manifold $(M,J,g)$ endowed with a
flat torsion-free connection $\n$ such that $\n \o =0$
and $d^{\n} J=0$, where $\o = g(\cdot,J \cdot)$. 

A {\em conical special para-K\"ahler manifold} $(M,J,g,\n, \xi )$ is a 
special para-K\"ahler manifold \\ $(M,J,g,\n)$ endowed with a timelike or a 
spacelike vector field $\xi$ such that $\n \xi =D \xi = {\rm Id}$,
where $D$ is the Levi-Civita connection. 
\ed     

It follows from the definition of a para-K\"ahler manifold that the
eigenspaces of $J$ are of the same dimension and involutive.  An
endomorphism field $J\in \Gamma ({\rm End}\,TM)$ with these properties
is called a {\em para-complex structure} on $M$. The pair $(M,J)$ is
then called a {\em para-complex manifold}. A smooth map $f: (M,J_M)
\rightarrow (N,J_N)$ between para-complex manifolds is called {\em
  para-holomorphic} if $df \circ J_M = J_N \circ df$. The
skew-symmetry of $J$ in the definition of a para-K\"ahler manifold
implies that the eigenspaces of $J$ are totally isotropic of dimension
$n=\frac{1}{2}\dim M$. In particular, $M$ is of even dimension $2n$
and $g$ is of signature $(n,n)$.

On any conical special para-K\"ahler manifold, the vector fields $\xi$
and $J\xi$ generate a free para-holomorphic action of a
two-dimensional Abelian Lie algebra. If the action can be integrated
to a free para-holomorphic action of a Lie group $A$ such that the
quotient map $M \rightarrow \bar{M}:=M/A$ is a para-holomorphic
submersion, then $\bar{M}$ is called a {\em projective special
  para-K\"ahler manifold}.  The quotient $\bar{M}$ carries a canonical
para-K\"ahler metric $\bar{g}$ compatible with the induced
para-complex structure $J$ on $\bar{M}=M/A$.

Next we explain the extrinsic construction of special para-K\"ahler
manifolds.  Recall that a {\em para-complex vector space} $V$ of
dimension $n$ is simply a free module $V\cong C^n$ over the ring
$C=\bR[e]$, $e^2=1$, of {\em para-complex numbers}. Notice that $C^n$
is a para-complex manifold with the para-complex structure $v\mapsto
ev$ and any para-complex manifold of real dimension $2n$ is locally
isomorphic to $C^n$. An $\bR$-linear map $\tau : V \rightarrow V$ on a
para-complex vector space is called {\em anti-linear} if $\tau (ev)=
-e\tau (v)$ for all $v\in V$. An example is the para-complex
conjugation $C^n \rightarrow C^n$, $z=x+ey \mapsto \bar{z}:=x-ey$.
Let $(V,\Omega )$ be a para-complex symplectic vector space of
dimension $2n$ endowed with a real structure (i.e.\ an anti-linear
involution) $\tau$ such that \re{tauOEqu} holds true. Then
\begin{equation}
\label{gammaEqpara}
\gamma := e\Omega (\cdot ,\tau \cdot )  
\end{equation}
is a  para-Hermitian form and $g_V:= {\rm Re}\, \gamma$ is a flat 
para-K\"ahler metric on $V$.
\bd A para-holomorphic immersion $\phi : M \rightarrow V$ from 
para-complex 
manifold $(M,J)$ of real dimension 2n into $V$ is called 
\begin{itemize}
\item[(i)] 
{\em nondegenerate} if $\phi^*\gamma$  is nondegenerate, 
\item[(ii)] {\em Lagrangian} if $\phi^*\Omega =0$ and 
\item[(iii)] {\em conical} if $\phi (p) \in d\phi (T_pM)$ 
and $\gamma (\phi (p),\phi (p))\neq 0$  for all $p\in M$. 
\end{itemize} 
\ed  

\bt \cite{CMMS} \label{CMMSThm} \begin{itemize}
\item[(i)] Any nondegenerate para-holomorphic Lagrangian immersion
  $\phi : M\rightarrow V$ induces on the para-complex manifold $(M,J)$
  the structure of a special para-K\"ahler manifold $(M,J,g,\n )$,
  where $g={\rm Re}\,\phi^*\gamma$ and $\n$ is determined by the
  condition $\n \phi^* \a =0$ for all $\a \in V^*$ which are real
  valued on $V^\tau$.
\item[(ii)]
Any conical nondegenerate para-holomorphic 
Lagrangian immersion $\phi : M\rightarrow V$ induces
on $(M,J)$ the structure of a  conical special para-K\"ahler manifold 
$(M,J,g,\n ,\xi )$. The vector field $\xi$ is determined 
by the condition $d\phi\,\xi (p) = \phi (p)$.  
\end{itemize} 
\et 

\subsubsection{The (affine) special para-K\"ahler manifold as a para-complex
 Lagrangian cone}  
Now we consider the real symplectic module $V_0=\wedge^3\bR^6$ 
of $G=\SL(6,\bR )$. For convenience, the standard basis of $\bR^6$ is denoted
by $(e_1,e_2,e_3,f_1,f_2,f_3)$.  The para-complexification 
$V:=V_0\ot C= \wedge^3\bC^6\cong C^{20}$ of $V_0$ is a 
para-complex symplectic vector space 
endowed with a real structure $\tau$ such
that $V^\t =V_0$ and \re{tauOEqu}. We put $u_i:=e_i+ef_i$ and consider
the orbit 
\[ V\subset M=\GL^+(6,\bR)p\cong \GL^+(6,\bR)/\,\SL(3,\bR)\times \SL(3,\bR)\] 
of the element $p=u_1\wedge u_2\wedge u_3$.
\bt \label{SL6paraThm} $M=\GL^+(6,\bR)p\subset V$ 
is a nondegenerate para-complex Lagrangian cone. The inclusion $M\subset V$
induces on $M$ an $\SL(6,\bR )$-invariant special
para-K\"ahler structure. The image  
$\bar{M} =\pi (M)\cong  \SL(6,\bR )/\,{\rm S} (\GL(3,\bR)\times \GL(3,\bR))$ under 
the projection $\pi : V' \rightarrow P(V')$ is 
a homogeneous projective special para-K\"ahler manifold of real dimension 
18. Here $V'\subset V$ stands for the subset of nonisotropic vectors. 
\et 

\pf Using the formulas \re{OEqu} and \re{gammaEqpara} 
with $\nu = eu_{123}\wedge \bar{u}_{123}=-8e_{123}\wedge f_{123}$ 
we compute:  $\gamma (p,p)=1$. This shows that $M=\GL^+(6,\bR)p\subset V'$
consists of spacelike vectors. The tangent $T_pM\subset V$ has the 
following basis:
\[ (u_{123},\bar{u}_1\wedge u_{23},\bar{u}_2\wedge u_{31},\bar{u}_3\wedge u_{12},
\bar{u}_2\wedge u_{23},\bar{u}_2\wedge u_{12},\bar{u}_3\wedge u_{23},
\bar{u}_1\wedge u_{13},\bar{u}_3\wedge u_{13},-\bar{u}_1\wedge u_{12}).\] 
The restriction of $\Omega$ to $T_pM$ is zero in view of \re{OEqu}. 
The para-Hermitian form $\gamma|_{T_pM}$ is represented by the 
matrix \re{MatEqu}. This shows that $\gamma|_{T_pM}$ is nondegenerate. 
Hence, the inclusion $M\subset V$ is a conical  
para-holomorphic nondegenerate Lagrangian immersion. In virtue
of Theorem \ref{CMMSThm} it induces an $\SL(6,\bR )$-invariant 
conical special para-K\"ahler structure $(J,g,\n ,\xi )$ on $M$,
which in turn induces a homogeneous projective special para-K\"ahler structure
on $\bar{M}=\pi (M)\subset P(V')$. 
\end{proof} 
\subsubsection{The special para-K\"ahler manifold as an open orbit
of $\GL^+(6,\bR )$ on $\wedge^3\bR^6$}  
The conical special K\"ahler manifold $M=\GL^+(6,\bR)/\,\SL(3,\bC)$ 
described in Theorem \ref{SL6Thm} as a complex Lagrangian cone 
$M\subset V_0\ot \bC$ 
can be identified with the open $\GL^+(6,\bR )$-orbit 
 $\{ \l <0\}\subset V_0=\wedge^3\bR^6$, where $\l$ stands for the quartic 
$\SL(6,\bR )$-invariant \eqref{lambdarho}: 
\bp The projection $\rho : V_0\ot \bC \rightarrow V_0$, $v\mapsto {\rm Re}\, v$,  
induces a $\GL^+(6,\bR )$-equivariant
diffeomorphism from the Lagrangian cone 
$C(X_{(0)})\subset V_0\ot \bC$ described in Theorem \ref{SL6Thm} onto 
$\{ \l <0\}\subset V_0$: 
\[ C(X_{(0)})\cong \{ \l <0\} \cong \GL^+(6,\bR)/\,\SL(3,\bC). \]
\ep  

\pf This follows from the fact that $\l$ is negative on 
the real part of a non-zero 
decomposable $(3,0)$-vector, since $\{ \l <0\}\cong \GL^+(6,\bR)/\,\SL(3,\bC)$ 
is connected, see Proposition \ref{orbitsof3forms}. 
\end{proof}

In that picture the complex structure is less obvious than in the 
complex Lagrangian picture but the flat connection and symplectic 
(K\"ahler) form are simply 
the given structures of the symplectic vector space $V_0$. The 
complex structure is then obtained from the metric, which is the 
Hessian of the function $f=\sqrt{ |\lambda |}$. (We consider $\l$ as a 
scalar invariant by choosing a generator of $\wedge^6\bR^6$.)  
This route was followed by Hitchin in \cite{H1}. 

The other open $\GL^+(6,\bR )$-orbit $\{ \l >0\}\subset V_0$ cannot
be obtained as the real image of a $\GL^+(6,\bR )$-orbit on the 
complex Lagrangian cone $C(X) \subset V_0\ot \bC$  over the
highest weight variety $X\subset P(V_0\ot \bC)$. In fact,
$\GL^+(6,\bR )$ has only one open orbit on $X$, see Theorem \ref{SL6Thm},
and that orbit maps to $\{ \l <0\}\subset V_0$ under the projection 
$V_0\ot \bC \rightarrow V_0$.  Instead we have: 
\bp  \label{paraProp} The projection $\rho : V_0\ot C \rightarrow
V_0$, 
$v\mapsto {\rm Re}\, v = \frac{v+\bar{v}}{2}$,  
induces a $\GL^+(6,\bR )$-equivariant diffeomorphism  
from the para-complex Lagrangian cone $M=\GL^+(6,\bR )p
 \subset V_0\ot C$, $p=u_{123}$,  
described in Theorem \ref{SL6paraThm} onto the open $\GL^+(6,\bR )$-orbit 
$\{ \l >0\}\subset V_0$: 
\[ M\cong \{ \l >0\} \cong \GL^+(6,\bR)/\,(\SL(3,\bR)\times \SL(3,\bR)). \]
\ep 

\pf It suffices to check that ${\rm Re}\, u_{123}\in  \{ \l >0\}$.
This follows from the expression
\begin{eqnarray*}
2 {\rm Re}\, u_{123} &=& 
((e_1+ e f_1)\wedge (e_2+ e f_2)\wedge (e_3+ e f_3)  
+ (e_1- e f_1)\wedge (e_2- e f_2)\wedge (e_3- e f_3)) \\
&=& 
((e_1+f_1)\wedge (e_2+f_2)\wedge (e_3+f_3)  
+ (e_1-f_1)\wedge (e_2-f_2)\wedge (e_3-f_3)), 
\end{eqnarray*}
since a three-vector belongs to $\{ \l >0\}$ if and only if it can be written
as the sum of two decomposable three-vectors which have a non-trivial
wedge product, see Proposition \ref{orbitsof3forms}. 
\end{proof} 

Let us denote by $\n$ the standard flat connection of the vector space
$V_0$, by $\xi$ the position vector field, by $\o$ its $\SL(6,\bR
)$-invariant symplectic form and by $X_f$ the Hamiltonian vector field
associated to the function $f=\sqrt{\l }$. Then we have: \bt The data
($J=\n X_f$, $g=\o \circ J$, $\n$, $\xi$) define on $U= \{ \l
>0\}\subset V_0$ an $\SL(6,\bR )$-invariant conical special
para-K\"ahler structure.  \et

\pf Any three-vector $\psi\in U$ can be written uniquely as
$\psi^++\psi^-$ with decomposable three-vectors $\psi^\pm$ such that
$\psi^+\wedge \psi^-=f(\psi ) \nu$, cf. \eqref{phieuler} and Corollary
\ref{hat}. Differentiation at $\psi$ in direction of a vector $\xi \in
V_0$ yields
\[ (df_\psi\xi)\nu  = (\psi^+-\psi^-)\wedge \xi = \o (\psi^+-\psi^-,\xi )\nu,\]
that is 
\begin{equation} 
  X_f(\psi ) =\psi^+-\psi^-. 
\end{equation}
{}Using this equation, we can calculate $J=\n X_f$ by ordinary
differentiation in the vector space $V_0$. The result is that $J$ acts
as identity on the subspace $\wedge^3 E_+ \oplus \wedge^2E_+\wedge
E_-\subset V_0=\wedge^3\bR^6$ and as minus identity on the subspace
$\wedge^3 E_- \oplus \wedge^2E_-\wedge E_+$ where $E_{\pm} =
\mathrm{span}\{\a\,\inter\,\psi^\pm\;|\;\a \in \wedge^2(\bR^6)^* \}$
denotes the support of the three-vectors $\psi^+$ and $\psi^-$. This
shows that $J^2 = {\rm Id}$ and that $J$ is skew-symmetric with
respect to $\o$. To prove that the data ($J$, $g=\o \circ J$, $\n$,
$\xi$) define on $U= \{ \l >0\}\subset V_0$ an $\SL(6,\bR )$-invariant
conical special para-K\"ahler structure, it suffices to show that
under the map $\rho : V_0 \ot C\rightarrow V_0$ these data correspond
to the conical special para-K\"ahler structure on $M=\GL^+(6,\bR )p
\subset V_0\ot C$, $p=u_{123}$, described in Theorem \ref{SL6paraThm}.
It follows from Proposition \ref{paraProp} and the definition of the
structures on $M$ that the data ($\o$, $\n$, $\xi$) on $U$ correspond
to the symplectic structure, flat connection and Euler vector field of
the conical special para-K\"ahler manifold $M$. One can check by a
simple direct calculation that the endomorphism $J$ on $T_{\rho (p)}U$
corresponds to multiplication by $e\in C$ on $T_pM\subset V_0\ot
C$. This proves the theorem.
 
Alternatively, we give now a direct argument which avoids the use of 
Theorem \ref{SL6paraThm}. The structure $J$ on $U$ satisfies
\[ d^\n J = d^\n \n X_f = (d^\n)^2 X_f =0,\]
since $\n$ is flat. This easily implies the integrability of
$J$ by expanding the brackets in the Nijenhuis tensor 
using that $\n$ has zero torsion. In view of the fact that 
$J$ is skew-symmetric for $\o$, we conclude that $(U,J,g = \o \circ J)$ is
para-K\"ahler. Finally, the flat torsion-free connection $\n$ 
satisfies not only $d^\n J =0$ but also $\n \o =0$, 
since the two-form $\o$ on $V_0$ is constant. This proves 
that $(U,J,g,\n)$ is special para-K\"ahler. Now we check that
$(U,J,g,\n ,\xi )$ is a conical special para-K\"ahler manifold, that is
$\n \xi =D\xi ={\rm Id}$. 
It is clear that $\n \xi = {\rm Id}$, since $\xi$ is the 
position vector field in $V_0$. To prove the second equation, we first 
remark that the Levi-Civita connection is given by
\[ D = \n +\frac{1}{2}J\n J.\] 
(It suffices to check that $D$ is metric and torsion-free.) 
Therefore, the equation $D\xi ={\rm Id}$ is reduced to 
$\n_\xi J =0$. Let us first prove that 
$\xi$ is para-holomorphic, that is $L_\xi J =0$. 
{}By homogeneity of $f$ and $\o$, we have the Lie derivatives 
\[ L_\xi f = 2 f, \quad L_\xi df = 2 df,\quad L_\xi \o = 2\o,\quad 
L_\xi \o^{-1} = -2\o^{-1}\]
and, hence,  
\[ L_\xi X_f = 0.\]
The latter equation implies
\[ L_\xi J=L_\xi (\n X_f) = 0,\]
since $\xi$ is an affine (and even linear) vector field. 
Using $\n_\xi -L_\xi = \n \xi = {\rm Id}$ we get that 
\[ \n_\xi J = L_\xi J + [{\rm Id}, J] = 0.\]  
\end{proof}

\end{document}